\documentclass[12pt]{report}
\usepackage[latin1]{inputenc}
\usepackage[brazil]{babel}
\usepackage{latexsym}
\usepackage{amsfonts}
\usepackage{amssymb}
\usepackage{amsmath}
\usepackage[mathcal]{euscript}

\input amssym.def
\newtheorem{Def}{Definição}[chapter]
\newtheorem{Teo}[Def]{Teorema}

\newtheorem{Lema}[Def]{Lema}
\newtheorem{Prop}[Def]{Proposição}

\newcommand{\dps}{\displaystyle}
\newcommand{\cqd}{
\hspace*{\fill}
\rule{7pt}{7pt}
}

\begin{document}
\vspace*{2cm}

\thispagestyle{empty}
\centerline{\large{\bf K-teoria de operadores pseudodiferenciais}}
\centerline{\large{\bf na reta com símbolos semiperiódicos}} 

\normalsize
\vspace{1.3cm}
\centerline{\bf Cintia Cristina da Silva}

\vspace{2 cm}
\centerline{\bf TESE  APRESENTADA}
\centerline{\bf AO}
\centerline{\bf INSTITUTO DE MATEMÁTICA E ESTATÍSTICA}
\centerline{\bf DA}
\centerline{\bf UNIVERSIDADE DE SÃO PAULO}
\centerline{\bf PARA}
\centerline{\bf OBTENÇÃO DO GRAU DE DOUTORA}
\centerline{\bf EM}
\centerline{\bf MATEMÁTICA }

{
\vspace{1.0 cm}
\begin{center}
Área de Concentração: {\bf Análise}\\
Orientador: {\bf Prof. Dr. Severino Toscano do Rego Melo} 
\end{center}


\vspace{0.5cm}
\centerline{\bf -São Paulo, março de 2005-}

\vspace{0.5cm}
\normalsize 
\centerline{\it \small{Durante a elaboração deste trabalho a autora
  recebeu auxílio financeiro da CAPES.}}

\pagebreak

\thispagestyle{empty}
\vspace*{2cm}

\large
\centerline{\bf K-teoria de operadores pseudodiferenciais }
\centerline{\bf na reta com símbolos semiperiódicos }
\vspace*{0.2cm}
\normalsize

\normalsize
\vspace{2cm}
\begin{flushright}
Este exemplar corresponde à redação \\
final da tese devidamente corrigida\\
e defendida por Cintia Cristina da Silva\\
e aprovada pela comissão julgadora.


\end{flushright}
\vspace{3cm}
Banca examinadora:

\begin{itemize}
\item Prof.~Dr.~Severino Toscano do Rego Melo (orientador) - IME-USP
\item Profa.~Dra.~Cristina Cerri - IME-USP
\item Prof.~Dr.~Ruy Exel Filho - UFSC
\item Prof.~Dr.~Antonio Roberto da Silva - UFRJ
\item Profa.~Dra.~Beatriz Abadie - URU  
\end{itemize}

\vspace{2cm}

\centerline{São Paulo, março de 2005.}

\tableofcontents
 \chapter*{Notações}
\addcontentsline{toc}{chapter}{Notações}
\hspace{8mm} Durante este trabalho adotamos a convenção de que um ideal é
sempre um ideal bilateral.
\begin{itemize}
\item $\mathcal{L}(E) := \{$ operadores limitados sobre o espaço de Banach $E \}$.

\item $\mathcal{K}(E) := \{$ operadores compactos sobre o espaço de Banach $E \}$.

\item $\mathcal{K}_{E} := \mathcal{K}(L^2(E))$.

\item $\mathcal{L}_{E} := \mathcal{L}(L^2(E))$.

\item $S^1 :=$ círculo unitário em $\mathbb{C}$.

\item $S^1_{+} = S^1 \times \{+ \infty\}$.

\item $S^1_{-} = S^1 \times \{- \infty\}$.

\item $C(X) := \{$ funções complexas contínuas sobre $X \}$.

\item $C^{\infty}(X) := \{$ funções complexas infinitamente diferenciáveis sobre $X \}$.

\item $C^{\infty}_{c}(X) := \{$ funções de $C^{\infty}(X)$ com suporte compacto $\}$.

\item $CS(\mathbb{R}) := \{$ funções contínuas de $\mathbb{R}$ em $\mathbb{C}$
  com limites em $+ \infty$ e $- \infty \} \cong C([- \infty, + \infty])$, $[-
  \infty, + \infty]$ denotando a compactificação de $\mathbb{R}$ que se obtém
  acrescentando os pontos $- \infty$ e $+ \infty$.

\item $CS(\mathbb{Z}) := \{$ sequências complexas indexadas por $\mathbb{Z}$ com limites em $+
\infty$ e $- \infty \}$.

\item $P_{2 \pi} := \{$ funções contínuas de $\mathbb{R}$ em $\mathbb{C}$ $2 \pi$-periódicas
$\}$.

\item $C_{0}(\mathbb{R}) := \{$ funções contínuas de $\mathbb{R}$ em
  $\mathbb{C}$ com limite zero em $\pm \infty \}$.

\item $\mathcal{S}(\mathbb{R}) := \{ \phi \in C^{\infty}(\mathbb{R}) :
\dps{\sup_{x}} | x^{\alpha} \phi^{\beta}(x)| < \infty, \, \alpha, \, \beta \in
\mathbb{N} \}$.

\item $F$ denota a Transformada de Fourier em $\mathbb{R}$: $$Fu(x) =
\frac{1}{\sqrt{2 \pi}} \int_{\mathbb{R}} e^{- ixt}u(t) \, dt.$$

\item $Id$ denota o operador identidade em um espaço de Hilbert.

\item $ind$ denota o índice de um operador de Fredholm.

\item $w$ denota o número de rotação de uma curva fechada em $\mathbb{C}
  \setminus \{0\}$.  

\item $\delta_{0}$ é a aplicação exponencial da sequência exata cíclica de seis
termos de K-teoria.

\item $\delta_{1}$ é a aplicação  do índice da sequência exata cíclica de seis
termos de K-teoria.

\item $*$ denota a operação convolução.

\item $T^*$ denota o operador adjunto de $T$.

\item $\hat{T} = FT$.

\item $< \cdot, \cdot >$ denota produto interno.

\item $a(M)$ denota o operador de multiplicação por $a$. Neste trabalho utilizamos

$a(M) \in \mathcal{L}_{\mathbb{R}}$, isto é $a(M)u(x) = a(x)u(x)$, $a \in L^{\infty}(\mathbb{R})$.

$a(M) \in \mathcal{L}_{\mathbb{Z}}$, isto é $a(M)(u_{j})_{j} =
(a(j)u_{j})_{j}$, $a: \mathbb{Z} \rightarrow \mathbb{C}$ limitada.

$a(M) \in \mathcal{L}_{S^1}$, isto é $a(M)u(x) = a(x)u(x)$, $a \in L^{\infty}(S^1)$.

\item $Aut(A)$ denota o conjunto dos automorfismos de $A$, $A$ uma $C^*$- álgebra.

\item $GL(A)$ denota o grupo dos elementos inversíveis de $A$, $A$ uma $C^*$-
  álgebra com unidade.

\item $\mathcal{U}(A)$ denota o grupo dos elementos unitários de $A$, $A$ uma
  $C^*$- álgebra com unidade.

\item $\tilde{A}$ denota a unitização da $C^*$- álgebra $A$.

\item $\bar{\varphi}$: dado o $*$-homomorfismo $\varphi: A \rightarrow B$
  temos $\bar{\varphi}: \tilde{A} \rightarrow \tilde{B}$ onde $\bar{\varphi}(a
  + \alpha 1_{\tilde{A}}) = \varphi(a) + \alpha 1_{\tilde{B}}$, $a \in A$ e
  $\alpha \in \mathbb{C}$.
\end{itemize}
\chapter*{Introdução}
\addcontentsline{toc}{chapter}{Introdução}
\hspace{8mm}Dado o espaço de Hilbert $L^2(\mathbb{R})$, seja
$\mathcal{L}(L^2(\mathbb{R}))$ o conjunto dos operadores limitados de
$L^2(\mathbb{R})$. É conhecido que $\mathcal{L}(L^2(\mathbb{R}))$ é uma $C^*$-
álgebra (veja, por exemplo, \cite{[Murphy]}, exemplo 2.1.3).

Denotemos por $\mathcal{A}$ a $C^*$- subálgebra de
$\mathcal{L}(L^2(\mathbb{R}))$ gerada pelas seguintes classes de operadores:

1. As multiplicações $b(M)$ por funções $b \in CS(\mathbb{R})$.

2. As multiplicações $e^{ijM}$, $j \in \mathbb{Z}$.

3. Os operadores da forma $b(D) := F^{-1}b(M)F$, $F$ denotando a transformada
   de Fourier e $b$ como no item anterior.

Denotemos por $SL$ a $C^*$- subálgebra de $\mathcal{L}(L^2(S^1))$ gerada pelos
operadores $a(M)$ de multiplicação por funções $a \in C^{\infty}(S^1)$ e por
todos os operadores da forma $b(D_{\theta}) := F^{-1}_{d}b(M)F_{d}$, onde
$F_{d}: L^2(S^1) \longrightarrow L^2(\mathbb{Z})$ denota a transformada de
Fourier discreta e $b(M) \in \mathcal{L}(L^2(\mathbb{Z}))$ o operador de
multiplicação por uma sequência $(b_{j})_{j \in \mathbb{Z}}$ que possua
limites quando $j$ tende a $+ \infty$ e $- \infty$. É sabido que $SL$ contém o
ideal $\mathcal{K}_{S^1}$ dos operadores compactos de $L^2(S^1)$ e que o
quociente $SL / \mathcal{K}_{S^1}$ é isomorfo a álgebra das funções contínuas
em $S^1 \times \{ -1, 1\}$, com isomorfismo determinado por 
$$a(M) \longmapsto ((z, \pm 1) \mapsto a(z)) \qquad \mbox{e} \qquad b(D)
\longmapsto ((z, \pm 1) \mapsto \lim_{j \rightarrow \pm \infty} b_{j}).$$

Seja $\mathcal{E}$ o menor ideal fechado contendo todos os comutadores de
$\mathcal{A}$. Então $\mathcal{E}$ contém o ideal $\mathcal{K}_{\mathbb{R}}$
dos operadores compactos em $L^2(\mathbb{R})$ e é isomorfo \cite{[SM]} ao produto
tensorial de $C^*$- álgebras $SL \otimes \mathcal{K}_{\mathbb{Z}}$, onde 
$\mathcal{K}_{\mathbb{Z}}$ é o ideal  dos operadores compactos em
$L^2(\mathbb{Z})$ (a construção explícita de um isomorfismo depende da
escolha de um isomorfismo entre $L^2(\mathbb{R})$ e $L^2(S^1) \otimes
L^2(\mathbb{Z})$). Combinando-se esse isomorfismo com o do parágrafo anterior
e usando que $\mathcal{K}_{\mathbb{R}}$ é isomorfo a $\mathcal{K}_{S^1}
\otimes  \mathcal{K}_{\mathbb{Z}}$, obtém-se:
\begin{equation}
\mathcal{E} / \mathcal{K}_{\mathbb{R}} \cong C(S^1 \times \{-1, 1\}) \otimes
\mathcal{K}_{\mathbb{Z}}.
\end{equation}

Decorre de resultados gerais sobre ``álgebras de comparação'' , que são certas $C^*$-
álgebras geradas por operadores pseudodiferenciais em variedades, estudadas
por Cordes e colaboradores \cite{[Cd]}, que $\mathcal{K}_{\mathbb{R}} \subset
\mathcal{E}$ e que existe um isomorfismo
\begin{equation}
\mathcal{A} / \mathcal{E} \cong C(X \times \{ - \infty, + \infty \}),
\end{equation}
(onde $X = \{ (x, e^{ix}) \, : \, x \in \mathbb{R} \} \cup ( \{ - \infty, +
\infty \} \times S^1)$  é um subconjunto fechado de $[ - \infty, + \infty ]
\times S^1$ ) determinado por:

1. $[a(M)]_{\mathcal{E}} \mapsto ((x, e^{iy}), \pm \infty) \mapsto a(y), \, \forall \, (x,
   e^{iy}) \in X)$, se $a$ é $2 \pi $-periódica.

2. $[b(M)]_{\mathcal{E}}  \mapsto ((x, e^{iy}), \pm \infty) \mapsto b(x), \, \forall \, (x,
   e^{iy}) \in X)$, se $b \in CS(\mathbb{R})$.

3. $[b(D)]_{\mathcal{E}}  \mapsto((m, \pm \infty) \mapsto b(\pm \infty), \, \forall \, m \in
   X)$, se $b \in CS(\mathbb{R})$.

Denotamos, acima, por $[A]_{\mathcal{E}}$ a classe de $A \in \mathcal{A}$ no quociente
$\mathcal{A} / \mathcal{E}$.

O resultado central deste trabalho é o cálculo da K-teoria de $\mathcal{A}$, e
isto é feito nos Capítulos 3 e 4. Nos Capítulos 1 e 2 demonstramos, ou citamos
sem demonstrar, os resultados preliminares necessários para esse trabalho. A seguir detalhamos o conteúdo
de cada capítulo.

No Capítulo 1, provamos que $$\mathcal{K}_{\mathbb{R}} \subset \mathcal{E} \subset
\mathcal{A}$$ e estabelecemos o isomorfismo (2). O que difere esse
capítulo do trabalho de \cite{[SM]} é que para demonstrarmos que
$\mathcal{K}_{\mathbb{R}} \subset \mathcal{E}$ e a existência do isomorfismo
(2) não invocamos resultados sobre álgebras de comparação, embora utilizemos
muitas das idéias e técnicas de \cite{[Cd], [Cordes]}. Para descrevermos o
isomorfismo (1) utilizaremos o Teorema 2.6 de \cite{[SM]}, que será enunciado
  mas não demonstrado.

Denotemos por $\mathcal{C}$ a $C^*$- subálgebra de
$\mathcal{L}(L^2(\mathbb{R}))$ gerada pelas seguintes classes de operadores:

1. As multiplicações $a(M)$ por funções $a \in CS(\mathbb{R})$.

2. Os operadores da forma $b(D) := F^{-1}b(M)F$, $F$ denotando a transformada
   de Fourier e $b$ como no item anterior.

No Capítulo 2 provamos que
$\mathcal{K}_{\mathbb{R}} \subset \mathcal{C}$ e estabelecemos o isomorfismo
$$\mathcal{C} / \mathcal{K}_{\mathbb{R}} \cong C(M_{C})$$ (onde $M_{C} = \{
(x,y) \in [- \infty, + \infty] \times [- \infty, + \infty] \, : \, |x| + |y| =
\infty \}$). Introduzimos, também, nesse capítulo um pouco sobre a teoria de
espaços de Hardy (na reta e no círculo). Essa teoria foi utilizada na
demonstração de que, se $T \in \mathcal{C}$ é um operador de Fredholm, temos $ind(T) = w( \sigma_{T})$ onde $ind$ é o índice de
Fredholm de $T$, $w$ é o número de rotação  e $\sigma$ é a composição do
isomorfismo de Gelfand com a projeção canônica (observamos que, como $M_{C}$ é
homeomorfo ao círculo, faz sentido calcular $w(\sigma_{T})$).
Utilizaremos os resultados obtidos nesse capítulo nos capítulos que seguem.

A álgebra $\mathcal{C}$ foi estudada por Cordes e Herman \cite{[CH]} e Power \cite{[Pw]}. 
No trabalho de Power foi demonstrado que se $T \in \mathcal{C}$ é um operador
de Fredholm então temos
$ind(T) = w( \sigma_{T})$, este resultado estava apenas implícito em
\cite{[CH]}. Damos aqui uma apresentação detalhada e auto-contida dos
resultados de \cite{[CH],[Pw]}, usando a linguagem da K-teoria de álgebra de operadores.

Os isomorfismos (1) e (2) induzem uma sequência exata curta de $C^*$- álgebras
\begin{equation}
0 \longrightarrow C(S^1 \times \{-1, 1\}) \otimes \mathcal{K}_{\mathbb{Z}}
\longrightarrow \mathcal{A} / \mathcal{K}_{\mathbb{R}} \longrightarrow C(X
\times \{ - \infty, + \infty \}) \longrightarrow 0.
\end{equation}

No Capítulo 3 calculamos explicitamente os ``connecting mappings''  da
sequência exata cíclica de seis termos associada a (3). 

A composição do isomorfismo (1) com a projeção canônica (no quociente) pode ser estendida a um $*$-homomorfismo, que chamamos
$\gamma$, da álgebra toda para o conjunto das funções contínuas de $S^1
\times \{-1, 1\}$, tomando valores em $\mathcal{L}(L^2(\mathbb{Z}))$.
No Capítulo 4, calculamos precisamente a imagem de $\gamma$ e verificamos que
a mesma tem estrutura de produto cruzado. Sua K-teoria pôde então ser calculada pela
sequência exata cíclica de Pimsner-Voiculescu \cite{[PV]}. Obtivemos então uma nova
sequência exata curta de $C^*$- álgebras,  
\begin{equation}
0 \longrightarrow J_{0} / \mathcal{K}_{\mathbb{R}}
\longrightarrow \mathcal{A} / \mathcal{K}_{\mathbb{R}} \longrightarrow
\mbox{Im}\, \gamma \longrightarrow 0
\end{equation}
onde $J_{0}$ é a $C^*$- subálgebra de $\mathcal{L}(L^2(\mathbb{R}))$ gerada
por todos os operadores da forma $a(M)b(D)$, onde $a \in C_{0}(\mathbb{R})$ e
$b \in CS(\mathbb{R})$.
Calculamos explicitamente os ``connecting mappings''  da
sequência exata cíclica de seis termos associada a (4). 

Os resultados obtidos com os cálculos dos ``connecting mappings''  das
sequências exatas cíclicas de seis termos associadas a (3) e (4) nos levaram
ao cálculo da K-teoria de
$\mathcal{A}$ e a uma descrição explícita da fórmula do índice de Fredholm $ind:
K_{1}(\mathcal{A} / \mathcal{K}_{\mathbb{R}}) \longrightarrow \mathbb{Z}$.

Seja
\begin{equation}
A = \dps{\sum_{j=1}^{N}}a_{j}(M)b_{j}(D),
\end{equation}
onde $a_{1}, \dots, a_{N}$ pertencem à álgebra finitamente gerada pelas funções infinitamente diferenciáveis de
período $2 \pi$ e por $s(x) = x (1 + x^2)^{-1/2}$ e $b_{1}, \dots,
b_{N}$ estão na álgebra finitamente gerada por $\{ s, 1\}$. Temos que $A$ é um
operador 
pseudodiferencial clássico (isto é, com símbolo poli-homogêneo, veja \cite{[Hormander]}). As funções $a_{1},
\dots, a_{N}$ são chamadas de semiperiódicas. A álgebra gerada pelos
operadores do tipo (5) é densa em $\mathcal{A}$. Logo $\mathcal{A}$ também
pode ser vista como uma $C^*$- álgebra gerada por uma certa classe de operadores
pseudodiferenciais com símbolos semiperiódicos.

Os resultados descritos acima sobre a estrutura de $\mathcal{A}$ foram
estendidos em \cite{[M]} para uma $C^*$- álgebra, $\mathcal{B}$, gerada por operadores
pseudodiferenciais em variedades riemannianas da forma $B \times \mathbb{R}$,
onde $B$ é uma variedade compacta n-dimensional.

Aparentemente para a álgebra $\mathcal{B}$ a  $\mbox{Im} \, \gamma$ não
apresenta estrutura de produto cruzado. Sendo assim, não basta utilizar os
recursos descritos no Capítulo 4 para calcular sua K-teoria. Para utilizarmos
a estratégia do Capítulo 3, precisamos conhecer explicitamente alguns
elementos ``básicos'' dos grupos $K_{i}(\mathcal{B} / \mathcal{E})$ e
$K_{i}(\mathcal{E} / \mathcal{K})$, $i=0$, $1$, como por exemplo
seus geradores. Estas informações parecem ser bem mais difíceis de obter
nessa $C^*$- álgebra.

Para as álgebras $\mathcal{A}$ e $\mathcal{B}$ temos a seguinte sequência de composição (no
sentido de \cite{[Dx]}, 4.3.2)
$$0 \subset \mathcal{K}_{\mathbb{R}} \subset \mathcal{E} \subset \mathcal{A},$$
$$0 \subset \mathcal{K}_{\mathbb{R}} \otimes \mathcal{K}_{B} \subset \mathcal{E} \subset \mathcal{B}.$$
Muitas outras $C^*$- álgebras geradas por operadores pseudodiferenciais
apresentam sequência de composição como acima. Por exemplo, para a $C^*$- álgebra $\mathcal{U}_{M}$ gerada
pelos operadores b-pseudodiferenciais de Melrose numa variedade compacta com bordo
, Lauter \cite{[L]} achou a sequência de composição $$0 \subset \mathcal{K}_{\mathbb{R}} \subset
\mathcal{E} \subset \mathcal{U}_{M}$$ com $\mathcal{U}_{M} / \mathcal{E}$
isomorfo (via símbolo principal) a uma álgebra de funções contínuas e $\mathcal{E} / \mathcal{K}_{\mathbb{R}}$ isomorfo à soma
direta (indexada pelas componentes conexas do bordo) da álgebra de funções
contínuas de $\mathbb{R}$ com valores em operadores compactos. Para mais exemplos ver \cite{[Cordes], [MNS], [M1], [M2]}.

Em 1963, M. Atiyah e I. Singer em \cite{[AS]} publicaram a primeira prova do
teorema, que posteriormente passou a chamar-se \emph{Teorema de
  Atiyah-Singer}, que expressa o índice de Fredholm de um operador elíptico
somente em termos topológicos. Outras
provas podem ser encontradas em \cite{[AS1], [AS2], [BB], [PG]}.

Em 1995, R. Nest e B .Tsygan em \cite{[NT], [NT1], [NT2]} publicaram
generalizações do \emph{Teorema de Atiyah-Singer} em um contexto mais
algébrico. Para essas generalizações eles utilizaram invariantes algébricos
como cohomologia cíclica, K-teoria e caráter de Chern.

V. Nistor em 1997, apresentou em \cite{[NV]} uma nova prova do \emph{Teorema
  do índice de Connes-Moscovici para recobrimento} \cite{[CM]} usando apenas invariantes
  algébricos.

Desde então o interesse pelo estudo dos invariantes algébricos para álgebras
de operadores pseudodiferenciais tem crescido, veja por
exemplo, \cite{[LM], [LM1], [MNS], [MN]}. 

Temos também que álgebras de operadores pseudodiferenciais são álgebras não-comu-\\tativas
naturalmente associadas a situações geométricas singulares, e existem boas razões para esperar que o estudo
de invariantes dessas álgebras não-comutativas revelem informações sobre
a geometria e ajudem a entender o que vem a ser uma variedade não-comutativa \cite{[A1]}.
\chapter{A álgebra $\mathcal{A}$}
\hspace{8mm}Neste capítulo estudaremos uma subálgebra de
$\mathcal{L}(L^2(\mathbb{R}))$, que chamamos $\mathcal{A}$, obtida como o
fecho da álgebra gerada por:\\
i) operadores de multiplicação por funções em $CS(\mathbb{R})$.\\
ii) operadores de multiplicação por $e^{ijt}$, $j \in \mathbb{Z}$.\\
iii) operadores da forma $b(D) := F^{-1}b(M)F$, $b \in CS(\mathbb{R})$ e $F$
denotando a transformada de Fourier.

É fácil verificar que $\mathcal{A}$ é invariante por adjunto, isto é, dado $T
\in \mathcal{A}$ temos $T^{*} \in \mathcal{A}$. Portanto $\mathcal{A}$ é uma $*$- subálgebra de
$\mathcal{L}(L^2(\mathbb{R}))$.
Como $\mathcal{A}$ é fechada e tem unidade temos que $\mathcal{A}$ é uma
$C^*$- álgebra com unidade.

Seja $\mathcal{E}$ o ideal comutador de $\mathcal{A}$, isto é, $\mathcal{E}$ é
o ideal fechado gerado pelos comutadores $[a, b] = ab - ba$, $a, b \in
\mathcal{A}$. É fácil verificar que $\mathcal{A}/ \mathcal{E}$ é abeliano.
Seja $M_{A}$ o espectro da $C^*$- álgebra $\mathcal{A}/ \mathcal{E}$, isto é,
o conjunto dos homomorfismos não nulos de $\mathcal{A} / \mathcal{E}$ em
$\mathbb{C}$ munido da topologia fraca-$*$. Pelo Teorema 2.1.10, pág. 41 de
\cite{[Murphy]} temos que $\varphi : \mathcal{A} / \mathcal{E} \longrightarrow
C(M_{A})$ onde $a \longmapsto \hat{a}$ e $\hat{a}(\phi) = \phi(a)$, para todo
$a \in \mathcal{A} / \mathcal{E}$ e $\phi \in M_{A}$, é um *-isomorfismo
isométrico. O isomorfismo $\varphi$ é chamado de isomorfismo de Gelfand. 

No Teorema 1.17 apresentaremos uma fórmula explícita para $ \varphi: \mathcal{A}/ \mathcal{E} \longrightarrow
C(M_{A})$. Verificaremos no Lema 1.6 que $\mathcal{K}_{\mathbb{R}} \subset
\mathcal{E}$ e na Proposição 1.24 estabeleceremos um isomorfismo $ \psi:
\mathcal{E} / \mathcal{K_{\mathbb{R}}} \longrightarrow C(S^1 \times \{-1, 1
\}, \mathcal{K}_{\mathbb{Z}})$. Estes isomorfismos nos auxiliarão, no Capítulo
3, no cálculo da K- teoria de $\mathcal{A}$. 

Vamos considerar também a $C^*$- álgebra conjugada $\hat{\mathcal{A}} = F^{-1}\mathcal{A}F$. Os geradores de $\hat{\mathcal{A}}$ correspondentes a i), ii) e iii) são dados por:\\
î) operadores da forma $a(D)$, $ a \in CS(\mathbb{R})$.\\
îî) operadores translação $T_{j}$, onde $T_{j}(u)(\tau) := u( \tau +j)$ com $u
\in L^2(\mathbb{R})$, $j \in \mathbb{Z}$.\\
îîî) operadores de multiplicação por funções em $CS(\mathbb{R})$.
\section{Ideal comutador}

\hspace{8mm}O principal resultado desta seção é obtido na Proposição 1.7 onde mostraremos
que o ideal comutador, $\hat{\mathcal{E}} = F^{-1} \mathcal{E} F$, da álgebra
$\hat{\mathcal{A}}$ coincide com o fecho do conjunto:
 $$ \{ \dps{ \sum_{j=
       -N}^{N}}b_{j}(D)a_{j}(M)T_{j} + K\,; \, N \in \mathbb{N}\, , \,b_{j} \in CS(
   \mathbb{R})\, , \,a_{j} \in C_{0}(\mathbb{R})\, , \,K \in
   \mathcal{K}_{\mathbb{R}} \}.$$

Logo podemos concluir que $\mathcal{E}$ coincide com o fecho do conjunto:
$$ \{ \dps{ \sum_{j=
      -N}^{N}}b_{j}(M)a_{j}(D)e^{ijM} + K\,; \, N \in \mathbb{N}\,, \,b_{j} \in
  CS( \mathbb{R})\,, \,a_{j} \in C_{0}(\mathbb{R})\,, \,K \in
  \mathcal{K}_{\mathbb{R}} \}.$$

Temos que o comutador $[b(D), T_{j}]$ é igual a zero e é conhecido que 
\begin{equation}
[a(M),b(D)] \in \mathcal{K}_{\mathbb{R}},\qquad \mbox{ para} \quad a, b \in CS(\mathbb{R})
\end{equation}
(veja por exemplo \cite{[Cordes]}, capítulo III, Lema 9.4). E temos que
\begin{equation}
[T_{j}, a(M)] = (a(M+j) - a(M))T_{j} \qquad \mbox{para} \quad a \in
CS(\mathbb{R}) \quad
\mbox{e} \quad j \in \mathbb{Z}.
\end{equation}
\begin{Def}
Seja $H$ um espaço de Hilbert. Uma $C^*$- subálgebra $A$ de $\mathcal{L}(H)$ é irredutível se os únicos subespaços vetoriais fechados de $H$ que são invariantes por $A$ são $\{0\}$ e $H$.
\end{Def}
No Lema abaixo damos uma condição suficiente para que uma álgebra $A \subset \mathcal{L}(H)$ seja irredutível.
\begin{Lema}
Se $A \subset \mathcal{L}(H)$ é uma álgebra tal que $\{ au: a \in A\}$ é denso em $H$, para todo $u \in H$ não nulo, então $A$ é irredutível.
\end{Lema}
\textbf{Demonstração:}

Seja $M \subset H$ um subespaço vetorial fechado e invariante por $A$, onde $M \neq \{0\}$. Vamos provar que $M^{\bot} = \{0\}$, pois como $H = M \oplus M^{ \bot}$ teremos $M = H$.

Sejam $x \in M^{\bot}$ e $0 \neq y \in M$. Como $\{ay: a \in A \}$ é denso em $H$ temos que existem $a_{j} \in A$ tais que $a_{j}y \rightarrow x $ e $a_{j}y\; \bot\; x$. Então $<x,x> = < \dps{\lim_{j \rightarrow \infty}}a_{j}y,x> = \dps{\lim_{j \rightarrow \infty}}<a_{j}y,x> = 0 $.

Logo $x = 0$ o que implica que $M^{\bot} = \{ 0\}$. Portanto $A$ é irredutível.$\cqd$

Vamos verificar que o ideal comutador $\mathcal{E}$ da $C^*$- álgebra
$\mathcal{A}$ é irredutível e então utilizando o Teorema abaixo teremos  $\mathcal{K_{\mathbb{R}}} \subset \mathcal{E}$.
\begin{Teo}
Seja $A$ uma $C^*$- subálgebra de $\mathcal{L}(H)$ irredutível que tem interseção não nula com $K(H)$. Então $K(H) \subset A$.
\end{Teo}
\textbf{Demonstração:} Ver \cite{[Murphy]}, pág. 58.$\cqd$

Enunciaremos agora o Teorema de Stone -Weierstrass (caso complexo) que será
utilizado na demonstração do próximo resultado.
\begin{Teo}
Seja $B$ uma subálgebra de $C(X)$, $X$ Hausdorff compacto, com a propriedade
de que se $f \in B$ então $\bar{f} \in B$. Se $B$ é fechada na norma do supremo,
$||\cdot||_{\infty}$, e separa pontos (isto é, dados $x \neq y$ em $X$, existe $f
\in B$ tal que $f(x) \neq f(y)$), então $B = C(X)$ ou $B = \{ f \in C(X) : f(x)
=0 \}$ para algum $x$ fixado.
\end{Teo}
\textbf{Demonstração:} Ver \cite{[RS]}, pág. 266. $\cqd$
\begin{Lema}
Seja $B$ a $*$-álgebra gerada por $\{ a(x+j) - a(x) : a \in C_{0}( \mathbb{R}), \; j \in \mathbb{Z} \}$ então $B$ é densa em $C_{0}( \mathbb{R})$.
\end{Lema}
\textbf{Demonstração:}

Como $\mathbb{R}$ não é compacto, para utilizarmos o Teorema 1.4 vamos
considerar $X$ a compactificação de $\mathbb{R}$ por um ponto, isto é, $X =
\mathbb{R} \cup \{\infty\}$. Temos que $B$ é uma subálgebra de $C(X)$ e é
evidente que se $f \in B$ então $\bar{f} \in B$.

Sejam $a(x) = \frac{1}{1+ x^2}$ e $z, y \in \mathbb{R}$ com $ z \neq y$. Então
para algum $j \in \mathbb{Z}$, $a(z + j) - a(z) \neq a(y + j) - a(y)$. Consideremos a mesma $a$, $z = \infty$
e $y \in \mathbb{R}$ temos que $a(z + 1) - a(z) = 0$ e $ a(y + 1) - a(y) \neq
0$ se $y \neq - 1/2$ e $a(y + 2) - a(y) \neq 0$ se $y = - 1/2$. Logo $B$ separa pontos de $X$, então pelo Teorema de Stone-Weierstrass temos que $\bar{B} = \{ a \in C(X): a (+ \infty) = 0 \} = C_{0}( \mathbb{R})$. $ \cqd$

Seja $T: CS( \mathbb{R} ) \longrightarrow \mathcal{L}(L^2(\mathbb{R}))$ definida por $T(a) = a(M)$. $T$ é um $*$-homomorfismo injetivo, logo pelo Teorema 3.1.5 de \cite{[Murphy]} temos que $T$ é uma $*$-isometria. Utilizando esta $*$-isometria temos que a $*$-álgebra gerada por $\{a(M+j) - a(M): a \in C_{0}( \mathbb{R}), \; j \in \mathbb{Z} \}$ é densa em $\{b(M): b \in C_{0}( \mathbb{R}) \}$.

Durante este trabalho estaremos muitas vezes identificando a função $b \in CS(\mathbb{R})$ com o operador $b(M)$ através desta $*$-isometria.
\begin{Lema}
O ideal comutador $\mathcal{E}$ da $C^*$- álgebra $\mathcal{A}$ é irredutível.
\end{Lema}
\textbf{Demonstração:}

Consideremos $\eta(x) = \frac{1}{1+ x^2}$, $x \in \mathbb{R}$. É conhecido que
$$ F^{-1} \eta(x) = \frac{1}{\sqrt{4 \pi}} \int_{0}^{+ \infty} \frac{e^{-
    \delta}}{\sqrt{\delta}} e^{\frac{-|x|^2}{4 \delta}} \; d \delta$$
(veja, por exemplo, \cite{[RS]} pág. 128).

Temos que $\left| F^{-1} \eta (x) \right| \leq \frac{1}{ \sqrt{4 \pi}} \int_{0}^{+ \infty} \frac{e^{ - \delta}}{\sqrt{\delta}} \; d \delta$ e sabemos que 
\begin{equation}
 \int_{0}^{1} \frac{e^{ - \delta}}{\sqrt{\delta}} \; d \delta < \int_{0}^{1}
 \frac{1}{\sqrt{\delta}} \; d \delta < \infty \qquad \mbox{e}
\end{equation} 
\begin{equation}
\int_{1}^{+ \infty} \frac{e^{ - \delta}}{\sqrt{\delta}} \; d \delta \leq \int_{1}^{+ \infty} e^{ - \delta}\; d \delta < \infty.
\end{equation}

Logo por (1.3) e (1.4) temos que $F^{-1} \eta $  é limitada. Fazendo $u =
\frac{x}{2 \sqrt{\delta}}$ temos $du = \frac{dx}{2 \sqrt{\delta}}$, logo
$$\int_{\mathbb{R}} F^{-1} \eta (x)\, dx = \frac{2}{\sqrt{4 \pi}} \int_{0}^{+
  \infty} e^{- \delta} \int_{- \infty}^{ + \infty} e^{-u^2} du\, d \delta <
\infty$$ o que implica que $F^{-1} \eta \in L^1(\mathbb{R})$.

Sejam $\varphi \in \mathcal{S}(\mathbb{R})$ e $T = F^{-1} \eta$. Temos que
$\eta(D)(\varphi) = F^{-1}(\hat{\varphi} \hat{T})$ e como $T$ é uma
distribuição temperada pelo Teorema IX.4 de \cite{[RS]}
sabemos que $F^{-1}(\hat{\varphi} \hat{T}) = \frac{1}{\sqrt{2 \pi}} T *
\varphi$. Como $T \in L^1(\mathbb{R})$ e  $\varphi \in
\mathcal{S}(\mathbb{R}) \subset L^1(\mathbb{R})$ podemos escrever
$$\frac{1}{\sqrt{2 \pi}} T * \varphi(x) = \frac{1}{\sqrt{2 \pi}}
\int_{\mathbb{R}} T(x-y) \varphi (y) \, dy \, \mbox{ para todo} \,x \in \mathbb{R}.$$ 

Logo $\eta(D)$ é um operador integral com núcleo $T(x-y) > 0$. Tomemos $ 0
\neq u \in L^2(\mathbb{R})$ e $ \psi \in C^{\infty}_{c} ( \mathbb{R})$ tal que
$\psi u \neq 0$. Temos que $ \eta(D) \psi u \neq 0$, já que $\eta(M)$ é injetor, e vamos verificar agora que $\eta(D) \psi u$ é contínua.

Sabemos que para todo $x \in \mathbb{R}$, $\eta(D) \psi u (x) =
\int_{\mathbb{R}} T(x - y)( \psi u)(y) \; dy$. Já verificamos que $T \in L^1(
\mathbb{R})$ e é limitada. Por Holder temos que $\psi u \in L^1( \mathbb{R})$
e como $| T(x - y)( \psi u)(y) | \leq C |( \psi u)(y)|$ segue que $T(x - y)(
\psi u)(y)$ é limitada por uma função em $L^1( \mathbb{R})$, independente de
$x$. Então pelo Teorema da Convergência Dominada temos que para cada $x_{0}\in
\mathbb{R}$ $\dps{\lim_{x \rightarrow x_{0}}} \int_{\mathbb{R}} T(x -y)( \psi u)(y) \; dy = \int_{\mathbb{R}} \dps{\lim_{x \rightarrow x_{0}}}T(x -y)( \psi u)(y) \; dy$.

Usando novamente o Teorema da Convergência Dominada temos que $\dps{\lim_{v \rightarrow v_{0}}}T(v) = T(v_{0})$, logo $\dps{\int_{ \mathbb{R}}} \dps{\lim_{x \rightarrow x_{0}}} T(x -y)( \psi u)(y) \; dy = \int_{ \mathbb{R}} T(x_{0} -y)( \psi u)(y) \; dy$, ou seja, $\eta(D) \psi u$ é contínua.

Como $ \eta(D) \psi u \neq 0$ temos que existe $x_{0} \in \mathbb{R}$ tal que $ \eta(D) \psi u(x_{0}) \neq 0 $. Então temos que $\mbox{Re}\,(\eta(D) \psi u(x_{0})) \neq 0$ ou $\mbox{Im}\,(\eta(D) \psi u(x_{0})) \neq 0$ ou ambas o são.

Suponhamos que $\mbox{Re}\,(\eta(D) \psi u(x_{0})) \neq 0$. Então temos que $\mbox{Re}\,(\eta(D)
\psi u(x_{0})) > 0$ ou $\mbox{Re}\,(\eta(D) \psi u(x_{0})) < 0$. Se $\mbox{Re}\,(\eta(D) \psi
u(x_{0})) > 0$ então pela continuidade de $\eta(D) \psi u$ temos que $\mbox{Re}\,(\eta(D) \psi u) > 0$ em
um aberto $B \subset \mathbb{R}$ contendo $x_{0}$. Pelo Corolário I.2.2 de
\cite{[Ho]} existe $\varphi \in
C^{\infty}_{c}(B)$ tal que $ 0 \leq \varphi \leq 1$ e $ \varphi \equiv 1$ em
um aberto contido em $B$.

Portanto $\mbox{Re}\,( \varphi \eta(D) \psi u)(x) \geq 0$ para todo $x \in \mathbb{R}$. E
como $\eta(D)$ tem núcleo positivo $\mbox{Re}\,(\eta(D) \varphi \eta(D) \psi u)(x) > 0$ para todo $x \in \mathbb{R}$.

Se $\mbox{Re}\,(\eta(D) \psi u(x_{0})) < 0$ devemos proceder como no caso da
$\mbox{Re}\,(\eta(D) \psi u(x_{0})) > 0$. Se $\mbox{Im}\,(\eta(D) \psi u(x_{0})) \neq 0$ é análogo ao caso da parte
real ser diferente de zero. Portanto existe $\varphi \in C^{\infty}_{c}(B)$ tal que $[\eta(D) \varphi \eta(D) \psi u](x) \neq 0$ para todo $x \in
\mathbb{R}$.

Dado $v \in C^{\infty}_{c} ( \mathbb{R})$ seja $w = \frac{v}{ \eta(D)
  \varphi(M) \eta(D) \psi(M)u }$. Temos que 
$$ v = w(M) \eta(D) \varphi(M) \eta(D) \psi(M)u.$$

Como $w(M)$, $\varphi(M)$,
$\psi(M) \in \mathcal{A}$ e no Lema 1.5 vimos que $\eta(D) \in \mathcal{E}$
temos que $ v \in \{Eu : E \in \mathcal{E} \}$ o que implica que $C^{ \infty}_{c}(\mathbb{R}) \subset \{Eu : E \in \mathcal{E} \}$.

Usando que $C^{\infty}_{c} ( \mathbb{R})$ é denso em $L^2(\mathbb{R})$ temos que $\{Eu : E \in \mathcal{E} \}$ é denso em $L^2(\mathbb{R})$,
para todo $u \neq 0$. Então pelo Lema 1.2 podemos concluir que $\mathcal{E}$ é
irredutível.$\cqd$

Temos que $[a(M), b(D)] \in \mathcal{K}_{\mathbb{R}}\, \cap \, \mathcal{E}$, para $a,b
\in CS( \mathbb{R})$, e vamos verificar agora que $[a(M), b(D)] \neq 0$ para
algum $a,b \in CS( \mathbb{R})$. Vamos verificar mais, que $a$ e $b$ podem ser
tomados em $\mathcal{S}(\mathbb{R})$.

\textbf{Observação 1.6 1/2:} $[a(M), b(D)] \neq 0$ para algum $a, b \in \mathcal{S}(\mathbb{R})$.\\
De fato, sejam $a, b \in \mathcal{S}(\mathbb{R})$, é fácil verificar que $$[a(M),
b(D)]u(x) = \sqrt{2 \pi}\int_{\mathbb{R}}( (a(x) - a(y))F^{-1}b(x-y))u(y)dy.$$
Logo $[a(M), b(D)]$ é um operador integral com núcleo $K(x,y) = (a(x) -
a(y))F^{-1}b(x-y)$. Seja $T_{K} = [a(M), b(D)]$, afirmamos que $||K||_{L^2}
\neq 0$ implica $T_{K} \neq 0$.\\
De fato, pelo Teorema VI.23 de \cite{ [RS]} temos que a norma Hilbert-Schmidt $||T_{K}||_{HS} =
||K||_{L^2}$, logo $K = 0$ se, e só se,  $T_{K} = 0$.

Consideremos então $b(x) = e^{-x^2}$ e $a$ não constante. Temos que
$F^{-1}b(x-y) = \frac{\sqrt{2}}{2}e^{-(x-y)^2/4} > 0$ para todo $(x,y) \in \mathbb{R}^2$ logo $||K(x,y)||_{L^2}
\neq 0$. Segue da afirmação que $T_{K} \neq 0$.

Portanto pelo Teorema 1.3 temos $\mathcal{K}_{\mathbb{R}} \subset \mathcal{E}$.
\begin{Prop}
O ideal comutador $\hat{\mathcal{E}}$ da $C^*$- álgebra conjugada $\hat{\mathcal{A}}$ coincide com o fecho de 
$$\hat{\mathcal{E}}_{A} := \{ \dps{ \sum_{j= -N}^{N}}b_{j}(D)a_{j}(M)T_{j} + K\,; \, N \in \mathbb{N}\,, \,b_{j} \in CS( \mathbb{R})\,, \,a_{j} \in C_{0}(\mathbb{R})\,, \,K \in \mathcal{K}_{\mathbb{R}} \}.$$
\end{Prop}
\textbf{Demonstração:}

Como já foi provado temos que $\mathcal{K}_{\mathbb{R}} \subset
\mathcal{E}$, logo
\begin{equation}  
\mathcal{K}_{\mathbb{R}} \subset \hat{\mathcal{E}}.
\end{equation}

Por (1.2) temos que todos operadores da forma $b(D)(a(M+j) - a(M))T_{j}$ estão em $\hat{\mathcal{E}}$ para $b \in CS( \mathbb{R})$, $j \in \mathbb{Z}$ e $a \in C_{0}(\mathbb{R})$. Assim, pelo Lema 1.5 temos $\hat{\mathcal{E}}_{A} \subset \hat{\mathcal{E}}$.

Por outro lado, usando (1.1), (1.2), e (1.5), temos que
$\hat{\mathcal{E}}_{A}$ é uma subálgebra de $\hat{\mathcal{A}}$ que contém os
comutadores de todos os geradores de $\hat{\mathcal{A}}$ listados no começo do
capítulo. Além disso, $\hat{\mathcal{E}}_{A}$ é invariante a direita ou a
esquerda por multiplicações por estes operadores. Tomando limite, segue que o
fecho de $\hat{\mathcal{E}}_{A}$ é um ideal de $\hat{\mathcal{A}}$ que contém
todos os seus comutadores. Logo, $\hat{\mathcal{E}}$ está contido no fecho de $ \hat{\mathcal{E}}_{A}$.$\cqd$

\section{A álgebra $\mathcal{A}^{\sharp}$}

\hspace{8mm}Nesta seção estudaremos $\mathcal{A}^{\sharp}$ a $C^*$- álgebra abeliana gerada pelos operadores
do tipo i) e ii) da página 17. No Lema 1.10 explicitaremos o espectro,
$M_{\sharp}$, de $\mathcal{A}^{\sharp}$.
 
Seja $\mathcal{A}^{\sharp}_{0}$ a álgebra finitamente gerada por $CS(
\mathbb{R})$ e $P_{2 \pi}$. Sejam $\chi_{+}$ e $\chi_{-} \in CS( \mathbb{R})$ tais que
$\chi_{\pm}( \pm \infty) = 1$, $\chi_{+} + \chi_{-} = 1$ e $\chi_{\pm}(t) = 0$
se $\mp t > 1$. É fácil verificar que todo $a \in \mathcal{A}^{\sharp}_{0}$ é
da forma\\[0.5cm]
\centerline{$a = a_{+} \chi_{+} + a_{-} \chi_{-} + a_{0}$, $a_{\pm} \in P_{2 \pi}$, $a_{0} \in C_{0}( \mathbb{R})$,}\\[0.5cm]
onde a escolha de $a_{+}$, $a_{-}$ e $a_{0}$ é única. Vamos provar que $\mathcal{A}^{\sharp}_{0}$ é fechada e assim teremos $\mathcal{A}^{\sharp} = \mathcal{A}^{\sharp}_{0}$.

Seja $f_{n} = a_{n} \chi_{+} + b_{n} \chi_{-} + c_{n} \in
\mathcal{A}^{\sharp}_{0}$ e consideremos a sequência $(f_{n})_{n}$ convergindo
uniformemente para $d$. Vamos provar que $d \in\mathcal{A}^{\sharp}_{0}$.  

Para cada $y > 10$ temos que $||a_{n}||_{\infty} = \dps{\sup_{x \in [y, y + 2 \pi]}} |a_{n}(x)| \leq \dps{\sup_{x \in [y, y + 2 \pi]}} |a_{n}(x) + c_{n}(x)| + \dps{\sup_{x \in [y, y + 2 \pi]}} |c_{n}(x)| \leq \dps{\sup_{x \geq 1}} |a_{n}(x) + c_{n}(x)| + \dps{\sup_{x \in [y, y + 2 \pi]}} |c_{n}(x)|$.
Fazendo $y \rightarrow \infty$ temos que $||a_{n}||_{\infty} \leq \dps{\sup_{x \geq 1}} |a_{n}(x) + c_{n}(x)| \leq ||f_{n}||_{\infty}$. 

Com isto temos que $(a_{n})_{n}$ é Cauchy e pelo mesmo argumento $(b_{n})_{n}$
é Cauchy. Digamos que $a_{n} \rightarrow a$ e $b_{n} \rightarrow b$, $a,\; b
\in P_{2 \pi}$. Temos então que $d - a \chi_{+} - b \chi_{-} = \lim c_{n} \in C_{0}( \mathbb{R})$, logo $d \in \mathcal{A}^{\sharp}_{0}$. Portanto $\mathcal{A}^{\sharp}_{0}$ é fechada. Provamos então que:
\begin{Prop}
Dado $a \in \mathcal{A}^{\sharp}$ e fixados $\chi_{+}$ e $\chi_{-}$ como acima, estão unicamente determinados $a_{0} \in C_{0}( \mathbb{R})$ e $a_{\pm} \in P_{2 \pi}$ tais que $a = a_{+} \chi_{+} + a_{-} \chi_{-} + a_{0}$.
\end{Prop}
\hspace{8mm}Vamos achar o espectro $M_{ \sharp}$ da $C^*$- álgebra
$\mathcal{A}^{\sharp}$. Para isto vamos definir duas classes de funcionais lineares multiplicativos.

1) A cada $t \in \mathbb{R}$ vamos definir\\
\centerline{ $w_{t}:\mathcal{A}^{\sharp} \longrightarrow \mathbb{C}$}\\
\centerline{ $\quad \quad f \longmapsto f(t)$} 

2) Para $\theta \in \mathbb{R}$,\\
\centerline{ $w_{\theta,+}: \mathcal{A}^{\sharp} \longrightarrow \mathbb{C}$}\\
\centerline{ $\quad \quad \quad \qquad \quad \quad \quad \quad \quad \quad f \longmapsto
  \dps{\lim_{k \rightarrow \infty}}f( \theta + 2 \pi k) := f^{+}_{\theta}$}\\
\centerline{ $w_{\theta,-}: \mathcal{A}^{\sharp} \longrightarrow \mathbb{C}$}\\
\centerline{ $\quad \quad \qquad \quad \quad \quad \quad \quad \quad \quad f \longmapsto \dps{\lim_{k \rightarrow
      \infty}}f( \theta - 2 \pi k) := f^{-}_{\theta}$}

Usando que para $f \in \mathcal{A}^{\sharp}$ temos $f = f_{+} \chi_{+} + f_{-}
\chi_{-} + f_{0}$ com $f_{\pm} \in P_{2 \pi}$ e $f_{0} \in C_{0} (
\mathbb{R})$, é fácil verificar que os limites existem e que $ \mathbb{R} \ni \theta
\longmapsto f^{\pm}_{\theta} \in \mathbb{C}$ são funções $2 \pi$-periódicas.
\begin{Prop}
Os funcionais $w_{t}$ e $w_{ \theta, \pm}$, $t, \theta \in \mathbb{R}$ são todos os funcionais lineares multiplicativos de $\mathcal{A}^{\sharp}$. E temos que $ \mathbb{R}$ é denso em $M_{ \sharp}$.
\end{Prop}
\textbf{Demonstração:}

Temos que $\mathbb{R} \subset M_{ \sharp}$, pois basta considerarmos a
identificação
\begin{displaymath}
 \begin{array}{cccc}
   i: & \mathbb{R} & \longrightarrow & M_{ \sharp}\\
      & t          & \longmapsto     & w_{t}.\\
 \end{array}
\end{displaymath}
Vamos verificar que $i$ é contínua e injetora.

$\triangleright$ injetividade\\
Sejam $t,\; s \in \mathbb{R}$ com $t \neq s$ e $\varphi \in CS( \mathbb{R})$ tal que $\varphi(t) \neq 0$ e $ \varphi(s) = 0$. Temos que $\varphi \in \mathcal{A}^{\sharp}$ e $w_{t}( \varphi) \neq w_{s}( \varphi)$, isto é, $i(s) \neq i(t)$.   

$\triangleright$ continuidade\\
Seja $\{t_{ \alpha} \} \in \mathbb{R}$ uma rede com $t_{\alpha} \rightarrow t$. Queremos mostrar que $i(t_{ \alpha}) \rightarrow i(t)$.

Temos que $\varphi(t_{ \alpha}) \rightarrow \varphi(t)$ para todo $\varphi \in \mathcal{A}^{\sharp}$. Portanto $w_{t_{ \alpha}}( \varphi) \rightarrow w_{t}( \varphi)$ para todo $\varphi \in \mathcal{A}^{ \sharp}$, logo $i$ é contínua.

Suponhamos que $w \in M_{ \sharp}$ e seja $ s \in CS( \mathbb{R})$ a função $s(t) =
\frac{t}{\sqrt{1 + t^2}}$.

Como $ -1 \leq s \leq 1$ temos que $ -1 \leq w(s) \leq 1$. Vamos provar que se $|w(s)| < 1$ então $w = w_{a}$ onde $a \in \mathbb{R}$ resolve $s(a) = w(s)$.

Por Stone-Weierstrass temos que os polinômios em $s$ são densos em $CS(
\mathbb{R}) =$\\$=  C([ - \infty, + \infty])$, logo $w$ e $w_{a}$ coincidem sobre $CS( \mathbb{R})$.

Seja $ \chi \in C_{0}( \mathbb{R})$ satisfazendo $\chi(a) = 1$, temos $\chi(M) e^{iM} \in CS( \mathbb{R})$ e então $w_{a}(e^{iM}) = e^{ia} = w_{a}( \chi(M) e^{iM}) = w( \chi(M) e^{iM}) = w( \chi(M)) w( e^{iM}) = w(e^{iM})$. Provando que $w$ e $w_{a}$ coincidem sobre o conjunto de geradores de $\mathcal{A}^{ \sharp}$.

Para o caso $w(s) = \pm 1$, seja $\theta$ tal que $w(e^{iM}) = e^{i \theta}$. Temos que $w_{\theta, \pm}(s) = \dps{\lim_{k \rightarrow \infty}} s( \theta \pm 2 \pi k) = \pm 1$. Portanto $w_{\theta, \pm}(s) = w(s)$.

Usando novamente o fato de que os polinômios em $s$ são densos em $CS( \mathbb{R})$ temos que $w$ e $w_{\theta, \pm}$ coincidem em $CS( \mathbb{R})$.

Temos que $w_{\theta, \pm}(e^{iM}) = \dps{\lim_{k \rightarrow \infty}}e^{i(
  \theta \pm 2 \pi k)} = e^{i \theta}$. Portanto $w(e^{iM}) = w_{\theta,
  \pm}(e^{iM})$, logo $w$ e $w_{\theta, \pm}$ coincidem sobre os geradores
de $\mathcal{A}^{\sharp}$.

Agora só nos falta provar que $ \mathbb{R}$, ou melhor, $i( \mathbb{R})$ é
denso em $M _{ \sharp}$.

Seja $K = \overline{i( \mathbb{R})} \subset M_{ \sharp}$ e suponhamos que $K \neq 
M_{ \sharp}$. Então existe $w_{0} \in M_{ \sharp} \setminus K$, e pelo Lema de Urysohn temos que existe $\hat{f} \in C(M_{ \sharp})$ tal que $\hat{f}( w_{0}) = 1$ e $\hat{f} \equiv 0$ em $K$.

Temos que $\mathcal{A}^{ \sharp}$ é isomorfo a $C(M_{ \sharp})$, pelo Teorema
de Gelfand. Logo existe $f \in \mathcal{A}^{ \sharp}$ tal que $\hat{f}$ é a
imagem de $f$ por este isomorfismo  e temos que $ f(x) = w_{x}(f) = \hat{f}(w_{x}) = 0$ para todo $x \in \mathbb{R}$, ou seja, $f$ é nula. Absurdo, pois $\hat{f}$ não é nula.
Portanto $\mathbb{R}$ é denso em $M_{ \sharp}$.$\cqd$
\begin{Lema}
$M_{\sharp}$ é homeomorfo a um subconjunto de $[ - \infty, + \infty] \times
S^1$ a saber, $M_{\sharp} \approx \{(x, e^{ix}) : x \in \mathbb{R} \} \cup ( \{ -
\infty, + \infty \} \times S^1)$. 
\end{Lema}
\textbf{Demonstração:}

Pela Proposição acima  temos que $M_{\sharp} = \{ w_{t}: t \in \mathbb{R} \}
\; \cup \; \{ w_{ \theta, \pm} :  \theta \in \mathbb{R} \}$ com a topologia fraca-$*$. Usando as notações acima temos que $f^{ \pm}_{\theta + 2 \pi j} = f^{ \pm}_{ \theta}$, logo $f^{ \pm}_{\theta}$ é um par de funções em $\mathbb{R}/2 \pi \mathbb{Z} = S^1$.

Consideremos as inclusões
\begin{displaymath}
 \begin{array}{cccc}
 i_{1}: & CS( \mathbb{R}) & \longrightarrow & \mathcal{A}^{\sharp}\\
        & a               & \longmapsto      & a(M)\\
\end{array}
\end{displaymath}
\centerline{e}
\begin{displaymath}
 \begin{array}{cccc}
 i_{2}: & P_{2 \pi} & \longrightarrow & \mathcal{A}^{\sharp}\\
        & a               & \longmapsto      & a(M).\\
\end{array}
\end{displaymath}

Temos que o espectro de $CS( \mathbb{R})$ é o
conjunto $[- \infty, + \infty]$ e que o espectro de $P_{2 \pi}$ é $S^1$.
Consideremos $$i^{*}_{1}: M_{ \sharp} \longrightarrow [
- \infty, + \infty ], $$ onde $i_{1}^*(w) = w \circ i_{1}$, e também
$$i^{*}_{2}: M_{ \sharp} \longrightarrow S^1,$$ onde $i_{2}^*(w) = w \circ
i_{2}$. As aplicações $i^{*}_{1}$ e $i^{*}_{2}$ são chamadas de aplicações
duais de $i_{1}$ e $i_{2}$ respectivamente.
 
Seja $j : M_{\sharp} \longrightarrow  [ - \infty, + \infty ] \times S^1$ dada
por $j(w) = (i^{*}_{1}(w), i^{*}_{2}(w))$ para $w \in M_{\sharp}$.

Como $\mbox{Im} \;i_{1}$ e $\mbox{Im} \;i_{2}$ geram $\mathcal{A}^{\sharp}$ é
fácil ver que $j$ é injetora. Como $M_{\sharp}$ é compacto (Teorema 1.3.5,
\cite{[Murphy]}), $[- \infty, + \infty] \times S^1$ é Hausdorff e $j$ é contínua temos que $j$ é um
homeomorfismo sobre a imagem. Podemos então identificar topologicamente $M_{
  \sharp}$ com um subconjunto de $[- \infty, + \infty] \times S^1$. Pela
definição da aplicação dual, podemos concluir que $M_{\sharp} = \{ (x, e^{ix}): \; x \in \mathbb{R} \} \; \cup \; ( \{ - \infty, + \infty \} \times S^1)$.$\cqd$

Pelo Teorema de Gelfand, temos que $ \varphi : \mathcal{A}^{ \sharp} \rightarrow C(M_{
  \sharp})$ definido por 
$\varphi( a ) = \hat{a}$
é um $*$-isomorfismo isométrico. Mas sabemos quem são os funcionais de $M_{ \sharp}$, logo temos\\[0.5cm]
\centerline{$ \hat{a}(w_{t}) = w_{t}(a) = a(t)$, $t \in \mathbb{R}$\hspace{.2cm} e} \\[0.4cm]
\centerline{$\quad \hat{a}(w_{ \theta, \pm}) = w_{ \theta, \pm}(a) = a^{
  \pm}_{ \theta}$, $e^{i \theta}_{ \pm} \in S^1$.}

\section{A álgebra $J_{0}$}

\hspace{8mm}Nesta seção estudaremos $J_{0} \subset \mathcal{L}(L^2( \mathbb{R}))$ a $C^*$-
subálgebra gerada por $a(M)b(D)$, $a \in C_{0}( \mathbb{R})$ e $b \in CS(
\mathbb{R})$. No Lema 1.15 mostraremos que $J_{0} /
\mathcal{K}_{\mathbb{R}} \cong C_{0}(\mathbb{R} \times \{ - \infty, + \infty
\})$ e daremos uma fórmula explícita para o isomorfismo de Gelfand. E na Proposição 1.16 vamos verificar que $J_{0} \cap \mathcal{E} =
\mathcal{K}_{\mathbb{R}}$. Os resultados desta seção serão utilizados na próxima seção e também no
Capítulo 4.
\begin{Def}
O espaço símbolo $M_{A}$ de uma $C^*$- álgebra $A$ é o espectro da álgebra $A/E$, onde
$E$ é o ideal comutador de $A$.
\end{Def}
\begin{Def}
O $\sigma$-símbolo de uma $C^*$- álgebra $A$ é a composição da projeção
canônica $\pi$ com a transformada de Gelfand da $C^*$- álgebra $A / E$, onde
$E$ é o ideal comutador. \textbf{Notação:} $\sigma: A \rightarrow C(M_{A})$ com
$a \mapsto \sigma_{a} = \widehat{\pi(a)}$.
\end{Def}
\begin{Lema}
A $C^*$- álgebra $J_{0}$ é irredutível.
\end{Lema}
\textbf{Demonstração:}

Vamos verificar que $\{ au : a \in J_{0} \}$ é denso em $L^2(\mathbb{R})$, para todo $u$ não nulo. Então pelo Lema 1.2 teremos que $J_{0}$ é irredutível.

Na demonstração do Lema 1.6, vimos que, dados $v \in C^{\infty}_{c} (
\mathbb{R})$ existem $\eta, \varphi, \psi, w, u$ tais que 
$$ v = w(M)  \eta(D) \varphi(M) \eta(D) \psi(M)
u.$$
Temos que $w(M) \eta(D)$,
$\varphi(M) \eta(D)$ e $\psi(M) \in J_{0}$, logo $v \in \{ au : a \in J_{0}
\}$. Como $v \in C^{\infty}_{c}(\mathbb{R})$ é qualquer temos que
$C^{\infty}_{c}(\mathbb{R}) \subset \{ au : a \in J_{0} \}$. Portanto $\{ au :
a \in J_{0} \}$ é denso em $L^2( \mathbb{R})$, para todo $u$ não nulo.  $\cqd$
\begin{Lema}
O comutador $[a(M), b(D)]$ é compacto para $a, b \in C_{0}(\mathbb{R})$.
\end{Lema}
\textbf{Demonstração:}

Consideremos $a,b \in \mathcal{S}(\mathbb{R})$. É fácil verificar que
$$a(M)b(D)u(x) = \frac{1}{2\pi} \int_{\mathbb{R}} \Big( \int _{\mathbb{R}}
e^{it(x-y)}b(t)a(x) dt \Big) u(y) dy$$ e também que $$b(D)a(M)u(x) = \frac{1}{2\pi} \int_{\mathbb{R}} \Big( \int _{\mathbb{R}}
e^{it(x-y)}b(t)a(y) dt \Big) u(y) dy.$$
Logo $a(M)b(D)$ é um operador integral com núcleo $\hat{b}(y-x)a(x) \in
L^2(\mathbb{R}^2)$ e $b(D)a(M)$ é um operador integral com núcleo
$\hat{b}(y-x)a(y) \in L^2(\mathbb{R}^2)$. Portanto pelo Teorema VI.23 de \cite{[RS]}
temos que $a(M)b(D)$ e $b(D)a(M)$ são compactos e assim $[a(M), b(D)]$ é
compacto para $a, b \in \mathcal{S}(\mathbb{R})$.

O fim da demonstração segue do fato de $\mathcal{S}(\mathbb{R})$ ser denso em
$C_{0}(\mathbb{R})$.$\cqd$

\textbf{Observação 1.14 1/2:} Também temos que $a(M)b(D) \in \mathcal{K}_{\mathbb{R}}$
se $a \in C_{0}(\mathbb{R})$ e $b \in L^{\infty}(\mathbb{R}) \cap
L^2(\mathbb{R})$, a demonstração é análoga a do Lema acima.

Pelo Lema 1.14 e pela Observação 1.6 1/2 da página 22 temos que $J_{0} \cap \mathcal{K}_{\mathbb{R}} \neq 0$. Logo
pelo Teorema 1.3 temos que $\mathcal{K}_{\mathbb{R}} \subset J_{0}$ e assim $
\mathcal{K}_{\mathbb{R}} \subset J_{0} \cap \mathcal{E}$.

Seja $G$ o ideal comutador de $J_{0}$. Por (1.1) temos que $G \subset
\mathcal{K}_{\mathbb{R}}$ e que $G$ é não nulo pela Observação 1.6 1/2 da página 22. Como $ \mathcal{K}_{\mathbb{R}} \subset J_{0}$
segue que $G = \mathcal{K}_{\mathbb{R}}$, pois o único ideal não nulo de
$\mathcal{K}_{\mathbb{R}}$ é $\mathcal{K}_{\mathbb{R}}$.
\begin{Lema}
O espaço símbolo $M_{J}$ da $C^*$- álgebra $J_{0}$ é dado como:\\
\centerline{$M_{J} = \mathbb{R} \times \{- \infty, + \infty \}$.}\\
O $ \sigma$-símbolo é dado por:\\
\centerline{$ \sigma_{a}(m, \pm \infty) = a(m)$, $m \in \mathbb{R}$, $
  a \in C_{0}(\mathbb{R})$ \hspace{.2cm} e}
\centerline{$ \sigma_{b(D)}(m, \pm \infty) = b( \pm \infty )$, $m \in \mathbb{R}$, $b \in CS( \mathbb{R})$.}
\end{Lema}
\textbf{Demonstração:}

Consideremos a aplicação
\begin{displaymath}
 \begin{array}{cccc}
   i: &  C_{0}(\mathbb{R}) & \longrightarrow  &
   J_{0}/\mathcal{K}_{\mathbb{R}}\\    
      &    a               & \longmapsto      &
   [a(M)]_{\mathcal{K}_{\mathbb{R}}}.\\
\end{array}
\end{displaymath}
Temos que $i$ é injetora já que $a(M) \in \mathcal{K}_{\mathbb{R}}$ se, e só
se, $a \equiv 0$, logo sua aplicação dual $i^{*}: M_{J} \longrightarrow
\mathbb{R}$ onde $\mathbb{R}$ é o espectro de $C_{0}(\mathbb{R})$ é
sobrejetora (veja, por exemplo, \cite{[Cordes]}, pág. 314).

Para $w \in M_{J}$ sejam $t_{0} = i^{*}(w)$, $\phi \in C^{
  \infty}_{c}(\mathbb{R})$, $\phi$ real, com $\phi$ igual a um em uma
  vizinhança de $t_{0}$.

Sejam $A = \phi(M)s(D)$ onde $s(x) = \frac{x}{\sqrt{1 + x^2}}$ e $\xi_{0} =
\sigma_{A}(w) = w([\phi(M)s(D)]_{\mathcal{K}_{\mathbb{R}}})$. Temos que
$\xi_{0}$ é real, pois $(\phi(M)s(D))^* = s(D)\phi(M)$ que é côngruo módulo
compacto a $\phi(M)s(D)$. Logo $[\phi(M)s(D)]^*_{\mathcal{K}_{\mathbb{R}}} =
[\phi(M)s(D)]_{\mathcal{K}_{\mathbb{R}}}$, consequentemente $\xi_{0} =
\bar{\xi_{0}}$.

Definimos assim o par $(t_{0}, \xi_{0})$. Vamos verificar que $\xi_{0}$ não
depende da $\phi$ considerada.

Sejam $\phi$ e $\psi \in C^{\infty}_{c}(\mathbb{R})$ com $\phi$ igual a um em
uma vizinhança U de $t_{0}$ e $\psi$ igual a um em uma vizinhança V de
$t_{0}$. Seja $B = \psi(M)s(D)$, seja $\chi \in C^{\infty}_{c}(\mathbb{R})$ com suporte contido em $U \cap V$ e igual a um em uma vizinhança de
$t_{0}$. Como $(1 - \chi)( \phi - \psi) = \phi - \psi$ temos que $A - B = (1 -
\chi)(M)(A - B)$. Logo $\sigma_{A-B}(w) = \sigma_{A-B}(w) -
\sigma_{ \chi(M)(A-B)}(w) = \sigma_{A-B}(w) - \sigma_{\chi(M)}(w)\sigma_{A-B}(w)
  = 0$, pois $\sigma_{\chi(M)}(w) = \chi(t_{0}) = 1$. Portanto $\sigma_{A}(w) = \sigma_{B}(w)$.

Temos que $\phi(M)s(D)\phi(M)s(D)$ é côngruo módulo compacto a
$\phi^2(M)s^2(D)$ e que $\phi^2(M)(1 - s^2)(D)$ é compacto pelo Lema
1.14. Logo $A^2$ é côngruo módulo compacto a $\phi^2(M)$, com isto temos que
$\xi^2_{0} = w([A]^2_{\mathcal{K}_{\mathbb{R}}}) =
w([\phi(M)]^2_{\mathcal{K}_{\mathbb{R}}}) = \phi^2(i^*(w)) = 1$. Portanto
$\xi_{0} = \pm 1$.

Definimos então 
\begin{displaymath}
 \begin{array}{cccc}
 \nu: & M_{J} & \longrightarrow & \mathbb{R} \times \{-1, 1\}\\
      &  w    & \longmapsto     & (i^*(w), \xi),\\
\end{array}
\end{displaymath}
onde $\xi = w([\phi(M)s(D)]_{\mathcal{K}_{\mathbb{R}}})$, $\phi \in
C^{\infty}_{c}(\mathbb{R})$ com $\phi$ igual a um em uma vizinhança de
$i^*(w)$. 

Vamos provar agora que se $a \in C_{0}(\mathbb{R})$ e $\nu(w) =(t, \xi)$ então
$w([a(M)s(D)]_{\mathcal{K}_{\mathbb{R}}}) = a(t)\xi$. De fato, seja $\chi_{j}
\in C^{\infty}_{c}(\mathbb{R})$ com $\chi_{j}$ igual a um em $(-j, j)$ então
$\xi = w([\chi_{j}(M)s(D)]_{\mathcal{K}_{\mathbb{R}}})$. Temos que
$w([\chi_{j}(M)a(M)s(D)]_{\mathcal{K}_{\mathbb{R}}}) = a(t) \xi$ para todo
$j$ suficientemente grande, mas $\chi_{j}(M)a(M)s(D) \rightarrow a(M)s(D)$. Logo temos o que
queríamos.

Vamos verificar que $\nu$ é um homeomorfismo.

$\triangleright$ injetora

Sejam $w$, $w^{,} \in M_{J}$ tais que $(t_{0}, \xi_{0}) = \nu(w) = \nu(w^{,})
= (t, \xi)$. Temos então que $ t_{0} = i^*(w) = i^*(w^{,}) = t$ e $\sigma_{\phi}(w) =
\sigma_{\phi}(w^{,})$ para toda $\phi \in C^{\infty}_{c}(\mathbb{R})$. Como
$\xi = \xi_{0}$ temos também $\sigma_{\phi(M)s(D)}(w) =
\sigma_{\phi(M)s(D)}(w^{,})$, para todo $\phi \in C^{\infty}_{c}(\mathbb{R})$.

O conjunto de todos os $\phi(M)s(D)$ como acima geram $J_{0}$, módulo
$\mathcal{K}_{\mathbb{R}}$, logo $\sigma_{A}(w) = \sigma_{A}(w^{,})$ para todo
$A \in J_{0}$. Portanto $w = w^{,}$.

$\triangleright$ sobrejetora

Consideremos o operador paridade $P$, isto é $(Pf)(x) = f(-x)$, $f \in
L^2(\mathbb{R})$. É fácil verificar que$Pa(M)P = a(-M)$ e $Pb(D)P = b(-D)$
para $a \in C_{0}(\mathbb{R})$ e $b \in CS(\mathbb{R})$.

Dado $w \in M_{J}$ vamos definir $\lambda([A]_{\mathcal{K}_{\mathbb{R}}}) =
w([PAP]_{\mathcal{K}_{\mathbb{R}}})$, $A \in J_{0}$. Seja $\nu(w) = (t_{0},
\xi_{0}) = (i^*(w), \xi_{0})$ como $\lambda([a(M)]_{\mathcal{K}_{\mathbb{R}}})
= w([a(-M)]_{\mathcal{K}_{\mathbb{R}}})$ temos que $i^*(\lambda) = - t_{0}$.

Seja $\psi \in C^{\infty}_{c}(\mathbb{R})$ com $\psi$ igual a um em uma
vizinhança de $-t_{0}$. Consideremos $\tilde{\psi}(x) = \psi(-x)$ então
$\tilde{\psi}$ é igual a um em uma vizinhança de $t_{0}$. Temos que
$\lambda([\psi(M)s(D)]_{\mathcal{K}_{\mathbb{R}}}) =
w([\psi(-M)s(-D)]_{\mathcal{K}_{\mathbb{R}}}) =
-w([\tilde{\psi}(M)s(D)]_{\mathcal{K}_{\mathbb{R}}}) = - \xi_{0}$ já que $s$ é
uma função ímpar.

Para $t \in \mathbb{R}$ vamos definir $(T_{t}f)(x) = f(x + t)$, $f \in
L^2(\mathbb{R})$. É fácil verificar que $T_{t}a(M)T_{-t} = (T_{t}a)(M)$ e
$T_{t}b(D)T_{-t} = b(D)$ para $a \in C_{0}(\mathbb{R})$ e $b \in
CS(\mathbb{R})$. 

Dado $w \in M_{J}$ vamos definir $\tau([A]_{\mathcal{K}_{\mathbb{R}}}) =
w([T_{t}AT_{-t}]_{\mathcal{K}_{\mathbb{R}}})$, $A \in J_{0}$. Seja $\nu(w) =
(t_{0}, 1)$. Como $\tau([a(M)]_{\mathcal{K}_{\mathbb{R}}})
= w([(T_{t}a)(M)]_{\mathcal{K}_{\mathbb{R}}})$ temos que $i^*(\tau) = t_{0}
+ t$.

Seja $\psi \in C^{\infty}_{c}(\mathbb{R})$ com $\psi$ igual a um em uma
vizinhança de $t_{0} + t$. Consideremos $\tilde{\psi}(x) = \psi(x + t)$ então
$\tilde{\psi}$ é igual a um em uma vizinhança de $t_{0}$. Temos que
$\tau([\psi(M)s(D)]_{\mathcal{K}_{\mathbb{R}}}) =
w([(T_{t} \psi)(M)s(D)]_{\mathcal{K}_{\mathbb{R}}}) =
w([\tilde{\psi}(M)s(D)]_{\mathcal{K}_{\mathbb{R}}}) = 1$.

Portanto se $(t_{0}, 1) \in \mbox{Im} \, \nu$ então $(t_{0} + t, 1) \in
\mbox{Im} \, \nu$ para todo $t \in \mathbb{R}$. Portanto usando a
sobrejetividade de $i^*$ e o que foi visto acima temos que $\nu$ é
sobrejetora. 

$\triangleright$ contínua

Sejam $\{w_{\alpha} \}_{\alpha}$ uma rede em $M_{J}$ e $w \in M_{J}$ com
$w_{\alpha} \rightarrow w$. Seja também $\nu(w) = (t, \xi)$ e $\nu(w_{\alpha})
= (t_{\alpha}, \xi_{\alpha})$.

Como $w_{\alpha} \rightarrow w$ temos que $i^*(w_{\alpha}) \rightarrow i^*(w)$
logo $t_{\alpha} \rightarrow t$. Seja $\chi \in C^{\infty}_{c}(\mathbb{R})$
com $\chi$ igual a um em uma vizinhança de $t$, temos que $\xi =
w([\chi(M)s(D)]_{\mathcal{K}_{\mathbb{R}}})$. Seja $\alpha_{0}$ tal que
para todo $\alpha \geq \alpha_{0}$ temos $t_{\alpha} \in U$, $U$ aberto de
$\mathbb{R}$ onde $\chi$ é igual a um. Então para todo $\alpha \geq
\alpha_{0}$ temos que $\xi_{\alpha} =
w_{\alpha}([\chi(M)s(D)]_{\mathcal{K}_{\mathbb{R}}}) \rightarrow
w([\chi(M)s(D)]_{\mathcal{K}_{\mathbb{R}}}) = \xi$. Portanto $\nu$ é contínua.

$\triangleright$ a inversa é contínua

Consideremos $w$, $w_{\alpha}$, $(t, \xi)$, $(t_{\alpha}$, $\xi_{\alpha})$ como
antes e seja $(t_{\alpha}$, $\xi_{\alpha}) \rightarrow (t, \xi)$. Temos que
$t_{\alpha} = i^*(w_{\alpha}) \rightarrow i^*(w) = t$ logo
$w_{\alpha}([a(M)]_{\mathcal{K}_{\mathbb{R}}}) \rightarrow
w([a(M)]_{\mathcal{K}_{\mathbb{R}}})$ para todo $a \in C_{0}(\mathbb{R})$.

Então $w_{\alpha}([a(M)s(D)]_{\mathcal{K}_{\mathbb{R}}}) = a(t_{\alpha})
\xi_{\alpha} \rightarrow a(t) \xi =
w([a(M)s(D)]_{\mathcal{K}_{\mathbb{R}}})$. Como\\ $\{
[a(M)s(D)]_{\mathcal{K}_{\mathbb{R}}} : a \in C_{0}(\mathbb{R}) \}$ gera $J_{0}
    / \mathcal{K}_{\mathbb{R}}$, segue que $w_{\alpha} \rightarrow w$. Portanto $\nu$ é um homeomorfismo.   

Para concluir a demonstração basta usarmos o isomorfismo de $C_{0}(
\mathbb{R} \times \{- 1, 1 \})$ em $C_{0}( \mathbb{R} \times \{- \infty, +
\infty\})$. Observar que $s(\pm \infty) = \pm 1$ e usar que os polinômios em
$s$ são densos em $CS(\mathbb{R})$. $\cqd$
\begin{Prop}
Temos que $J_{0} \cap \mathcal{E} = \mathcal{K}_{\mathbb{R}}.$
\end{Prop}
\textbf{Demonstração:}

Já vimos que $\mathcal{K}_{\mathbb{R}} \subset J_{0} \cap \mathcal{E}$, vamos
verificar então que $J_{0} \cap \mathcal{E} \subset \mathcal{K}_{\mathbb{R}}$.

Decorre da Proposição 1.7 que $\mathcal{E}$ é o fecho do conjunto
$\Big \{ \dps{ \sum_{j= -N}^{N}}b_{j}(M)a_{j}(D)e^{ijM} + K$: $ N \in \mathbb{N}$, $b_{j} \in CS( \mathbb{R})$, $a_{j} \in C_{0}(\mathbb{R})$, $K \in \mathcal{K}_{\mathbb{R}}\Big \}$.

Sabemos que $a(D - j)e^{ijM} = e^{ijM}a(D)$, logo $a_{j}(D)e^{ijM} = e^{ijM}a_{j}(D + j)$. Portanto $\mathcal{E}$ é o fecho do conjunto\\
\centerline{$ \Big\{ \dps{ \sum_{j= -N}^{N}}b_{j}(M)e^{ijM}a_{j}(D + j)\; +\; K$:
  $ N \in \mathbb{N}$, $b_{j} \in CS( \mathbb{R})$, $a_{j} \in
  C_{0}(\mathbb{R})$, $K \in \mathcal{K}_{\mathbb{R}}\Big \} $.}

Seja $A \in (J_{0} \cap \mathcal{E})$ e consideremos $\chi_{j} \in
C^{\infty}_{c}( \mathbb{R})$ com $0 \leq \chi_{j} \leq 1$ e $\chi_{j} \equiv
1$ em $(-j, j)$ com $j \in \mathbb{N}$. Como $A \in J_{0}$ podemos calcular o
$\sigma$ - símbolo de $A$, logo 
$ \sigma_{A} = \dps{\lim_{j \rightarrow \infty}} \sigma_{\chi_{j} A}$.
Portanto $[A]_{\mathcal{K}_{\mathbb{R}}} = \dps{\lim_{j \rightarrow \infty}}
[\chi_{j} A]_{\mathcal{K}_{\mathbb{R}}}$ em $J_{0}/\mathcal{K}_{\mathbb{R}}$,
ou seja 
 $|| [ \chi_{j} A]_{\mathcal{K}_{\mathbb{R}}} - [A]_{\mathcal{K}_{\mathbb{R}}}|| \longrightarrow 0$ quando $j \longrightarrow \infty$. Logo existem $(K_{j})_{j}$ onde $K_{j} \in \mathcal{K}_{\mathbb{R}}$ tais que $\dps{\lim_{j \rightarrow \infty}} \chi_{j} A + K_{j} = A$.

Como $A \in \mathcal{E}$ temos, também, que $A = \dps{\lim_{k \rightarrow \infty}} \dps{ \sum_{l= -N_{k}}^{N_{k}}}b^{k}_{l}(M)e^{ilM}a^{k}_{l}(D + l) + S_{k}$ onde $b^{k}_{l} \in CS(\mathbb{R})$, $a^{k}_{l} \in C_{0}( \mathbb{R})$ e $S_{k} \in \mathcal{K}_{\mathbb{R}}$.

Logo $$A = \dps{\lim_{j \rightarrow \infty}} \dps{\lim_{k \rightarrow \infty}}
\dps{ \sum_{l= -N_{k}}^{N_{k}}} \chi_{j}(M) b^{k}_{l}(M)e^{ilM}a^{k}_{l}(D +
l) + S_{k} + K_{j}.$$ Como $\chi_{j}(M) b^{k}_{l}(M)e^{ilM}a^{k}_{l}(D +
l)$ é compacto, vem que $A$ é o limite de uma soma de compactos, logo compacto.
Portanto $(J_{0} \cap \mathcal{E}) \subset \mathcal{K}_{\mathbb{R}}$,
assim $J_{0} \cap \mathcal{E} = \mathcal{K}_{ \mathbb{R}}$.$\cqd$
 
\section{O espaço símbolo de $\mathcal{A}$}

\hspace{8mm}Nesta seção calcularemos o espaço símbolo, $M_{A}$, de
$\mathcal{A}$ e daremos uma fórmula explicíta para o homomorfismo $\varphi:
\mathcal{A} \longrightarrow C(M_{A})$.
\begin{Teo}
O espaço símbolo $M_{A}$ da $C^*$- álgebra $ \mathcal{A}$ é dado como:\\
\centerline{$M_{A} = M_{ \sharp} \times \{- \infty, + \infty \}$.}\\
O $ \sigma$-símbolo é dado por:\\
\centerline{$ \sigma_{a}(m, \pm \infty) = a(m)$, $m \in M_{ \sharp}$, $
  a \in \mathcal{A}^{ \sharp}$\hspace{.2cm} e}
\centerline{$ \sigma_{b(D)}(m, \pm \infty) = b( \pm \infty )$, $m \in M_{ \sharp}$, $b \in CS( \mathbb{R})$.}
\end{Teo}
\textbf{Demonstração:}

Consideremos as aplicações
\begin{displaymath}
 \begin{array}{cccc}
i_{1}: & \mathcal{A}^{ \sharp} & \longrightarrow &
  \mathcal{A}/\mathcal{E}\\
       &  a(M)                 & \longmapsto     &
[a(M)]_{\mathcal{E}}\\
\end{array}
\end{displaymath}
\centerline{e}
\begin{displaymath}
 \begin{array}{cccc}
i_{2}: & CS( \mathbb{R}) & \longrightarrow &
  \mathcal{A}/\mathcal{E}\\
       &  b          & \longmapsto     &
[b(D)]_{\mathcal{E}}.\\
\end{array}
\end{displaymath}

E suas aplicações duais 
\begin{displaymath}
 \begin{array}{cccc}
i^{*}_{1}: & M_{A} & \longrightarrow & M_{ \sharp}\\
\end{array}
\end{displaymath}
\centerline{e}
\begin{displaymath}
 \begin{array}{cccc}
 i^{*}_{2}: & M_{A} & \longrightarrow & [ - \infty, + \infty ],\\
\end{array}
\end{displaymath}
onde dado $w \in M_{A}$ temos $i^{*}_{j}(w) = w \circ i_{j}$, $j = 1$ e $2$.

Consideremos $i$ o produto das duais:\\
\centerline{$i : M_{A} \longrightarrow M_{ \sharp} \times [ - \infty, + \infty ]$,}\\
com $w \in M_{A}$ tendo imagem $i(w) =
( w \circ i_{1}, w \circ i_{2} )$.

Como $\mbox{Im}\; i_{1}$ e $\mbox{Im}\; i_{2}$ geram $\mathcal{A}/
\mathcal{E}$ é fácil ver que $i$ é injetora. Como $M_{A}$ é compacto
(Teorema 1.3.5, \cite{[Murphy]}), $M_{\sharp} \times [ - \infty, + \infty ]$ é
Hausdorff e $i$ é contínua temos que $i$ é um homeomorfismo sobre a imagem.

Usando este homeomorfismo como identificação temos que:\\[0.5cm]
\centerline{$ \sigma_{a}(m,t) = \hat{a}(m)$\hspace{3mm} e \hspace{3mm}$ \sigma_{b(D)}(m,t) = b(t)$,}\\[0.5cm]
 para $ a \in \mathcal{A}^{ \sharp}$, $b \in CS( \mathbb{R})$ e $(m,t) \in M_{ \sharp} \times [ - \infty, + \infty ]$.

Até aqui temos que $M_{A} \subset M_{ \sharp} \times [ - \infty, + \infty ]$.
Mas demonstraremos nos Lemas 1.20, 1.21 e 1.22 respectivamente que $M_{A} \cap
( \mathbb{R} \times  \{ - \infty, + \infty \}) \neq \emptyset$; $(x, + \infty)
\in M_{A} \cap ( \mathbb{R} \times  \{ - \infty, + \infty \})$ se, e só se, $
( - x, - \infty ) \in M_{A} \cap ( \mathbb{R} \times  \{ -
\infty, + \infty \})$; $(x, + \infty) \in M_{A} \cap ( \mathbb{R} \times  \{ -
\infty, + \infty \})$ se, e só se, $
( y, + \infty ) \in M_{A} \cap ( \mathbb{R} \times  \{ -
\infty, + \infty \}) $, para todo $y \in \mathbb{R}$. Logo $\mathbb{R} \times \{ - \infty, + \infty \} \subset
M_{A}$ e como $\mathbb{R}$ é denso em $M_{\sharp}$ temos que $M_{ \sharp}
\times \{ - \infty, + \infty \} \subset M_{A}$. No Lema 1.19 demonstraremos que $(m,t) \notin M_{A}$ se $|t| < \infty$, logo $M_{A} = M_{ \sharp}
\times \{ - \infty, + \infty \}$. $\cqd$
\begin{Lema}
Dado qualquer $x \in M_{ \sharp}$ existe $y \in [ - \infty, + \infty ]$ tal que $(x,y) \in M_{A}$.
\end{Lema}
\textbf{Demonstração:}

Dados $a(M) \in \mathcal{E} \cap \mathcal{A}^{ \sharp}$ e $ \chi \in C^{
  \infty}_{c}( \mathbb{R})$ temos que $ \chi(M) a(M) \in (J_{0} \cap \mathcal{E})
= \mathcal{K}_{ \mathbb{R}}$, já que $b(M) \in J_{0}$ se $b \in
  C_{0}(\mathbb{R})$ e $\mathcal{E}$ é ideal de $\mathcal{A}$. Temos então que
  $ \chi a(M) = 0$, logo $\chi a \equiv 0$ para todo $ \chi \in  C^{ \infty}_{c}( \mathbb{R})$ o que implica que $a \equiv 0$.
Portanto $\mathcal{E} \cap \mathcal{A}^{ \sharp} = \{ 0 \}$.

Consideremos a aplicação $i_{1}$ do teorema acima, como $\mathcal{E} \cap \mathcal{A}^{ \sharp} = \{ 0 \}$ temos
que $i_{1}$ é um isomorfismo sobre a imagem, logo $i^{*}_{1}$ é sobrejetora
  (veja, por exemplo, pág. 314 de \cite{[Cordes]}).

Então dado $x \in M_{ \sharp}$ existe $x' \in M_{A}$ tal que $i^{*}_{1}(x') =
x$ e $i(x') = (i^{*}_{1}(x'), i^{*}_{2}(x')) = (x, i^{*}_{2}(x')) \in
M_{A}$.$\cqd$
\begin{Lema}
$(m,t) \notin M_{A}$ se $|t| < \infty$.
\end{Lema}
\textbf{Demonstração:}

Dado $t_{0} \in \mathbb{R}$, podemos escolher $a \in C_{0}( \mathbb{R})$ tal que $a(t_{0}) \neq 0$. Temos que $a(D) \in \mathcal{E} = Ker \sigma$,
logo $ \sigma_{a(D)}(m,t) = a(t) = 0$, para todo $(m,t) \in M_{ A}$. Assim
$(m, t_{0}) \notin M_{A}$, para qualquer $ m \in M_{ \sharp}$ já que $a(t_{0})
\neq 0$.$\cqd$
\begin{Lema}
$M_{A} \cap ( \mathbb{R} \times  \{ - \infty, + \infty \}) \neq \emptyset$.
\end{Lema}
\textbf{Demonstração:}

Dado $x \in \mathbb{R} \subset M_{ \sharp}$ temos pelo Lema 1.18 que existe $y \in [ - \infty, + \infty ]$ tal que $(x,y) \in M_{A}$.
E pelo Lema 1.19 temos que se  $|y| < \infty$ então $(x,y) \notin M_{A}$.
Portanto $(x,y) \in M_{A} \cap \mathbb{R} \times \{ - \infty, + \infty \}$.$\cqd$
\begin{Lema}
$(x, + \infty) \in M_{A} \cap ( \mathbb{R} \times  \{ - \infty, + \infty \}) \Longleftrightarrow ( - x, - \infty ) \in M_{A} \cap ( \mathbb{R} \times  \{ - \infty, + \infty \})$.
\end{Lema}
\textbf{Demonstração:}

Consideremos o operador paridade $P$, isto é, $(Pf)(x) = f(-x)$, $f \in L^2(\mathbb{R})$. É fácil verificar que $Pa(M)P = a(-M)$ e $Pb(D)P = b(-D)$ para $a \in \mathcal{A}^{ \sharp}$ e $b \in CS( \mathbb{R})$.

Seja $g:  \mathcal{A}/ \mathcal{E} \longrightarrow  \mathcal{A}/ \mathcal{E}$, onde $[a]_{\mathcal{E}} \mapsto [PaP]_{\mathcal{E}}$. 
Dado $\lambda \in M_{A}$ definimos $w = \lambda \circ g$. Pode-se verificar
sem dificuldades que $w \in M_{A}$.

Se $i( \lambda) = (x, + \infty)$ então $i(w) = (-x, - \infty)$ e se $i( \lambda) = (-x, - \infty)$ então $i(w) = (x, + \infty)$.
Portanto temos o resultado.$\cqd$
\begin{Lema}
$(x, + \infty) \in M_{A} \cap ( \mathbb{R} \times  \{ - \infty, + \infty \})
\Longleftrightarrow ( y, + \infty ) \in M_{A} \cap ( \mathbb{R} \times  \{ -
\infty, + \infty \}) $, para todo $y \in \mathbb{R}$.
\end{Lema}
\textbf{Demonstração:}

Dado $t \in \mathbb{R}$ vamos definir $(T_{t}f)(x) = f(x+t)$, $ f \in L^2(
\mathbb{R})$. É fácil verificar que $T_{t}a(M)T_{-t} = (T_{t}a)(M)$ e
$T_{t}b(D)T_{-t} = b(D)$ para $a \in \mathcal{A}^{\sharp}$ e $b \in
CS(\mathbb{R})$.

Seja $g:  \mathcal{A}/ \mathcal{E} \longrightarrow  \mathcal{A}/ \mathcal{E}$,
onde $[a]_{\mathcal{E}} \mapsto [T_{t}aT_{-t}]_{\mathcal{E}}$, a partir daqui
a demonstração é análoga a do Lema anterior.$\cqd$ 
\section{O quociente $\mathcal{E}/\mathcal{K}_{ \mathbb{R}}$}

Nesta seção vamos construir o isomorfismo de $\mathcal{E}/\mathcal{K}_{
  \mathbb{R}}$ para $C( S^1 \times \{-1,1\}, \mathcal{K}_{ \mathbb{Z}})$.

Seja $SL$ a $C^*$- subálgebra de $\mathcal{L}(L^2(S^1))$ gerada pelos operadores $a(M)$
de multiplicação por funções $a \in C^{ \infty}(S^1)$, e por todos os
operadores $b(D_{\theta}) := F^{-1}_{d}b(M)F_{d}$, onde
\begin{displaymath}
\begin{array}{cccc}
F_{d}:&  L^2(S^1) & \rightarrow & L^2( \mathbb{Z})\\
      &   f       & \mapsto     &  (f_{j})_{j}    \\
\end{array}
\end{displaymath}
é a transformada de Fourier discreta, com $f_{j} = \frac{1}{\sqrt{2 \pi}}
\dps{\int_{0}^{2 \pi}} e^{- i j \theta} f(e^{i \theta}) d \, \theta$, $b \in CS(\mathbb{Z})$, $b(M)(u_{j})_{j}$\\ $= (b(j)u_{j})_{j}$ e $CS( \mathbb{Z})$ é o
  conjunto das sequências $(b_{j})_{j}$ que possuem limites quando $j$ tende a
  $+ \infty$ e $- \infty$.

Podemos verificar que $SL$ é irredutível (verificação análoga à de
$\mathcal{E}$) e que todos os seus comutadores são compactos (veja, por
exemplo, \cite{[Cordes]}), então pelo Teorema 1.2 temos que $\mathcal{K}_{S^1} \subset
SL$. Usando a mesma estratégia do Teorema 1.17 podemos verificar que o espaço
símbolo de $SL$ é $M_{SL} =
S^1 \times \{ -1, 1 \}$. Logo pelo Teorema de Gelfand temos que
$SL/ \mathcal{K}_{S^1}$ é isomorfo a $C( S^1 \times \{ -1, 1 \})$ e que
$\sigma_{a}(z, \pm 1) = a(z)$\hspace{3mm}
e \hspace{3mm}
$\sigma_{b(D_{\theta})}(z, \pm 1) = b( \pm \infty)$, $a \in C^{
  \infty}(S^1)$ e $b \in CS( \mathbb{Z})$. 

Logo $A \otimes K \longmapsto \sigma_{A} \otimes K$ induz um
$*$-isomorfismo  
\begin{equation}
\frac{SL \hat{\otimes} \mathcal{K}_{ \mathbb{Z}}}{\mathcal{K}_{S^1 \times
    \mathbb{Z}}} \cong \frac{SL \hat{\otimes} \mathcal{K}_{
    \mathbb{Z}}}{\mathcal{K}_{S^1} \hat{\otimes} \mathcal{K}_{ \mathbb{Z}}}
    \cong \frac{SL}{\mathcal{K}_{S^1}} \hat{\otimes} \mathcal{K}_{ \mathbb{Z}}
    \cong C( M_{SL}) \hat{\otimes} \mathcal{K}_{ \mathbb{Z}} \cong C( M_{SL}, \mathcal{K}_{ \mathbb{Z}}).
\end{equation}
(veja \cite{[Breuer]}, para mais detalhes). Acima $\hat{\otimes}$ denota o
produto tensorial de $C^*$- álgebras nucleares. Em particular $SL
\hat{\otimes} \mathcal{K}_{ \mathbb{Z}}$ denota o fecho do produto tensorial algébrico de $SL
\otimes \mathcal{K}_{ \mathbb{Z}}$ em $\mathcal{L}(L^2(S^1) \hat{\otimes} L^2(\mathbb{Z}))$.

Dado $u \in L^2( \mathbb{R})$ denotamos:\\
\centerline{$u^{ \diamond}( \varphi) = ( u( \varphi - j))_{ j \in \mathbb{Z}}$
  para cada $ \varphi \in \mathbb{R}$.}

Temos, pelo Teorema de Fubini, que $u^{ \diamond}( \varphi) \in L^2( \mathbb{Z})$ para quase todo $ \varphi \in \mathbb{R}$.

Para cada $\varphi \in \mathbb{R}$ definimos $Y_{ \varphi} = F_{d} M^{-
  \varphi } F^{-1}_{d}$, onde $M^{- \varphi } \in \mathcal{L}(L^2(S^1))$ é o
  operador de multiplicação pela função $(z \mapsto z^{- \varphi})$, isto é
  $(M^{- \varphi }f)(e^{i \theta}) = e^{- i \varphi \theta} f(e^{i \theta})$. Podemos verificar que $Y_{ \varphi }$, $\varphi \in
  \mathbb{R}$, é uma família de
  operadores unitários de $\mathcal{L}_{\mathbb{Z}}$ com as
  seguintes propriedades:\\[0.5cm]
\centerline{$Y_{ \varphi}Y_{w} = Y_{ \varphi + w}$, $ \varphi, w \in
  \mathbb{R}$ \hspace{5mm} e \hspace{5mm} $(Y_{k}u)_{j} = u_{j+k}$, $k \in \mathbb{Z}$, $ u \in L^2( \mathbb{Z})$.}

Podemos então definir o operador unitário (com $S^1 = \{ e^{2 \pi i \varphi}
\,:\, \varphi \in \mathbb{R} \}$)
\begin{displaymath}
\begin{array}{cccc}
W: & L^2( \mathbb{R}) & \longrightarrow & L^2(S^1, d \varphi; L^2( \mathbb{Z}))\\
   & u                & \longmapsto     & (Wu)( \varphi) = Y_{ \varphi}u^{ \diamond}( \varphi).\\
\end{array}
\end{displaymath}

Usando a primeira propriedade de $Y_{ \varphi}$ verifica-se facilmente que se
$b \in CS(\mathbb{R})$, então $\varphi \mapsto Y_{ \varphi}b( \varphi - M)
Y_{- \varphi}$ é 1-periódica e pode ser vista, portanto como um elemento de
$C(S^1, \mathcal{L}_{\mathbb{Z}})$.

Então pelo Teorema 2.6 de \cite{[SM]}, que enunciamos abaixo, podemos definir o
isomorfismo $\psi$.
\begin{Teo}
Temos que $W \hat{\mathcal{E}}W^{-1} = SL \hat{\otimes} \mathcal{K}_{ \mathbb{Z}}$.\\ Além disso, para $b \in CS( \mathbb{R})$, $a \in C_{0}( \mathbb{R})$ e $j \in \mathbb{Z}$, temos:\\
\centerline{$Y_{ \varphi}a( \varphi -M)Y_{- \varphi} \in C(S^1, \mathcal{K}_{ \mathbb{Z}})$ \hspace{.2cm}e}
\begin{equation}
W(b(D)aT_{j})W^{-1} = b(D_{ \theta})Y_{ \varphi}a( \varphi
  -M)Y_{- \varphi-j} + K,\hspace{3mm} K \in \mathcal{K}_{S^1 \times
    \mathbb{Z}},
\end{equation} 
onde denotamos por $f(\varphi)$ funções $S^1 \ni e^{2 \pi i \varphi} \mapsto f(\varphi) \in \mathcal{K}_{\mathbb{Z}}$. 
\end{Teo}
\begin{Prop}
Existe um $*$-isomorfismo\\[0.5cm]
\centerline{$ \psi: \mathcal{E}/\mathcal{K}_{ \mathbb{R}} \longrightarrow C(M_{SL}, \mathcal{K}_{ \mathbb{Z}})$}\\[0.5cm]
tal que se $\tilde{ \gamma}$ denota a composição de $ \psi$ com a projeção
canônica $\mathcal{E} \longrightarrow \mathcal{E}/\mathcal{K}_{ \mathbb{R}}$,
e se $E \in \mathcal{E}$ satisfaz $F^{-1}EF = b(D)aT_{j}$, $b \in CS(
\mathbb{R})$, $a \in C_{0}( \mathbb{R})$, $ j \in \mathbb{Z}$ então temos:\\
\centerline{$\tilde{ \gamma}_{E}( \varphi, \pm 1) = b( \pm \infty)Y_{ \varphi}a( \varphi -M)Y_{- \varphi-j}$, \hspace{3mm}$(e^{2 \pi i \varphi}, \pm 1) \in M_{SL}$.}
\end{Prop}
\textbf{Demonstração:}

Vamos definir $\psi$  como sendo a composição das aplicações:\\[1cm]
\centerline{$ \mathcal{E}/ \mathcal{K}_{ \mathbb{R}} \longrightarrow 
    \hat{\mathcal{E}}/ \mathcal{K}_{ \mathbb{R}} \longrightarrow  \frac{SL \hat{
      \otimes} \mathcal{K}_{ \mathbb{Z}}}{\mathcal{K}_{S^1 \times \mathbb{Z}}}
  \longrightarrow C(M_{SL}, \mathcal{K}_{ \mathbb{Z}})$,}\\[1cm]
onde, a primeira aplicação leva $A + \mathcal{K}_{ \mathbb{R}}$ para
$F^{-1}AF + \mathcal{K}_{ \mathbb{R}}$, a segunda para $WF^{-1}AFW^{-1} +
\mathcal{K}_{S^1 \times \mathbb{Z}}$ e a última é o isomorfismo
(1.6). Pelo Teorema 1.23 e pela equação (1.7), temos o resultado para $\tilde{ \gamma}_{E}( \varphi, \pm 1)$.$\cqd$

Estenderemos $\tilde{ \gamma}$, definida sobre $\mathcal{E}$, para a álgebra
inteira $\mathcal{A}$. Como $\mathcal{E}/ \mathcal{K}_{ \mathbb{R}}$ é um
ideal de $\mathcal{A}/ \mathcal{K}_{ \mathbb{R}}$, todo $A \in \mathcal{A}$
define um operador $T_{A}$ em $ \mathcal{L}(\mathcal{E}/ \mathcal{K}_{ \mathbb{R}})$ por\\
\centerline{$ T_{A}(E + \mathcal{K}_{ \mathbb{R}}) = AE + \mathcal{K}_{ \mathbb{R}}$.}\\
Assim define-se:\\
\centerline{$T: \mathcal{A} \longrightarrow  \mathcal{L}(\mathcal{E}/ \mathcal{K}_{ \mathbb{R}})$.}

Está claro que $||T_{A}|| \leq ||A||$. Vamos definir
\begin{equation}
\begin{array}{cccc}
 \gamma: & \mathcal{A} & \longrightarrow & \mathcal{L}(C(M_{SL},
  \mathcal{K}_{ \mathbb{Z}})) \\
& A & \longmapsto & \psi T_{A} \psi^{-1}.\\
\end{array}
\end{equation}

Para $E \in \mathcal{E}$, $\gamma_{E}$ é uma multiplicação por
$\tilde{\gamma}_{E} \in C(M_{SL}, \mathcal{K}_{ \mathbb{Z}})$ da Proposição 1.24. Identificando funções de $C(M_{SL}, \mathcal{K}_{ \mathbb{Z}})$ com operadores de multiplicação correspondentes em $\mathcal{L}(C(M_{SL}, \mathcal{K}_{ \mathbb{Z}}))$, nós podemos então dizer que $\gamma$ estende $\tilde{ \gamma}$.
\begin{Prop}
Existe um $*$-homomorfismo\\
\centerline{ $\gamma : \mathcal{A} \longrightarrow C(M_{SL}, \mathcal{L}_{ \mathbb{Z}})$}\\
estendendo $\tilde{ \gamma}$ da Proposição 1.24. Sobre os geradores de $\mathcal{A}$, $ \gamma$ é dada por:\\[0.8cm]
\centerline{$ \gamma_{a}( \varphi, \pm 1) = a( \pm \infty)$, $ a \in CS( \mathbb{R})$,}\\[0.8cm]
\centerline{$\gamma_{b(D)}( \varphi, \pm 1) = Y_{ \varphi}b(M - \varphi )Y_{- \varphi}$, $ b \in CS( \mathbb{R})$ e} \\[0.8cm]
\centerline{$\gamma_{e^{ijM}}( \varphi, \pm 1) = Y_{-j }$, $ j \in \mathbb{Z}$.}
\end{Prop}
\textbf{Demonstração:}

É suficiente provar as equações acima para a aplicação $\gamma$ definida em (1.8). Por continuidade, a imagem de $ \gamma$ estará então contida em $C(M_{SL}, \mathcal{L}_{ \mathbb{Z}})$ considerada como uma subálgebra fechada de $\mathcal{L}(C(M_{SL}, \mathcal{K}_{ \mathbb{Z}}))$.

Dado $ a \in CS( \mathbb{R})$, precisamos calcular $ \tilde{ \gamma}_{aE}$ em
termos de $\tilde{ \gamma}_{E}$, para $ E \in \mathcal{E}$. Pela Proposição
1.7, temos que basta considerar $E$ tal que

\centerline{$F^{-1}EF = d(D)cT_{k}$, $d \in CS( \mathbb{R})$, $c \in C_{0}( \mathbb{R})$, e $k \in \mathbb{Z}$.}

Obtemos então $F^{-1}(aE)F = (ad)(D)cT_{j}$ e portanto, pela Proposição 1.24,

$\quad \quad \quad \tilde{\gamma}_{aE}( \varphi, \pm 1) = a( \pm \infty)d( \pm \infty)Y_{ \varphi}c( \varphi -M )Y_{- \varphi - j} = a( \pm \infty) \tilde{ \gamma}_{E}( \varphi, \pm 1)$.

Para $E$ como acima e $b \in CS( \mathbb{R})$,

$ \quad \quad \quad F^{-1}(b(D)E)F = b(-M)d(D)cT_{j} = d(D)b(-M)cT_{j} + K$, $K \in \mathcal{K}_{ \mathbb{R}}$,\\
e, assim,

$ \quad \quad \quad \tilde{\gamma}_{b(D)E}( \varphi, \pm 1) = d( \pm \infty)Y_{ \varphi}b( - \varphi + M)c( \varphi -M )Y_{- \varphi - j} = Y_{ \varphi}b(M - \varphi ) Y_{- \varphi} \tilde{ \gamma}_{E}( \varphi, \pm 1)$.

Para o mesmo $E$ temos

$\quad \quad F^{-1}e^{ijM}EF = T_{j}d(D)cT_{k} = d(D)c(M + j)T_{j+k}$\\[0.5cm]
e, então, usando que $Y_{ -j}c( \varphi -M ) = c( \varphi - M + j)Y_{-j}$,

$\quad \quad \quad  \tilde{\gamma}_{e^{ijM}E}( \varphi, \pm 1) = d( \pm
\infty)Y_{ \varphi}c( \varphi -M + j)Y_{- \varphi - j - k} = Y_{-j}d( \pm
\infty) Y_{ \varphi}c( \varphi -M)Y_{- \varphi - k} = Y_{-j} \tilde{
  \gamma}_{E}( \varphi, \pm 1)$. A aplicação preserva $*$, pois preserva $*$
nos geradores.$\cqd$
\begin{Prop}
Seja $\mathcal{A}^{\diamond} \subset \mathcal{A}$ a subálgebra fechada gerada
pelos operadores do tipo ii) e iii) da página 17. Se $A \in
\mathcal{A}^{\diamond}$ então $f = (WF^{-1})A(WF^{-1})^{-1} \in C(S^1,
\mathcal{L}_{\mathbb{Z}})$ e $\gamma_{A}(z, 1) = f(z) = \gamma_{A}(z, -1)$
para todo $z \in S^1$. Em particular, $||\gamma_{A}|| = ||A||$ para todo $A
\in \mathcal{A}^{\diamond}$.
\end{Prop}
\textbf{Demonstração:}

Dado $b \in CS(\mathbb{R})$, para cada $\varphi \in \mathbb{R}$, $(b(\varphi -
j))_{j \in \mathbb{Z}} \in CS(\mathbb{Z})$; denotamos por $b(\varphi -M) \in
\mathcal{L}_{\mathbb{Z}}$ o operador de multiplicação por esta sequência.

Para $ j \in \mathbb{Z}$, $ \varphi \in \mathbb{R}$ e $u \in L^2(\mathbb{R})$ temos,

\centerline{$(Y_{ \varphi}b( \varphi -M )Y_{-j}Y_{- \varphi})(Y_{ \varphi} u^{ \diamond}(
\varphi)) = Y_{ \varphi}(b( \varphi - k)u( \varphi - k + j))_{k \in \mathbb{Z}}$,}
logo $WbT_{j}W^{-1} = Y_{ \varphi}b( \varphi -M )Y_{- \varphi - j}$, para $ b \in
CS(\mathbb{R})$ e $j \in \mathbb{Z}$. 

Por outro lado segue das equações da Proposição 1.25 que

$ \gamma_{b(D)e^{ijM}}(\varphi, \pm 1)= Y_{ \varphi}b( M - \varphi)Y_{- \varphi - j} =
Wb(-M)T_{j}W^{-1} = WF^{-1}b(D)e^{ijM}FW^{-1}$,
o que prova que $\gamma_{A}(z, 1) = \gamma_{A}(z, -1)$, para todo $A \in
\mathcal{A}^{ \diamond}$. Como $WF^{-1}$ é um operador unitário segue que
$||\gamma_{A}|| = ||A||$ para todo $A \in \mathcal{A}^{ \diamond}$.$\cqd$
\begin{Prop}
Seja $J_{0}$ a $C^*$- álgebra definida na seção 1.3. Temos que $\mbox{Ker} \,
\gamma = J_{0}$.
\end{Prop}
\textbf{Demonstração:}

Se $A \in J_{0}$, vimos na demonstração da Proposição 1.16 que $\dps{\lim_{j
    \rightarrow \infty}} \chi_{j} A + K_{j} = A$ onde $\chi_{j} \in
C^{\infty}_{c}( \mathbb{R})$ com $0 \leq \chi_{j} \leq 1$ e $\chi_{j} \equiv
1$ em $(-j, j)$ com $j \in \mathbb{N}$; $K_{j} \in \mathcal{K}_{\mathbb{R}}$.

Agora se $A \in \mathcal{A}$ e $\dps{\lim_{j
    \rightarrow \infty}} \chi_{j} A + K_{j} = A$ temos que $A \in J_{0}$, já
     que\\ $A = \dps{\lim_{k \rightarrow \infty}} \dps{ \sum_{l=
    -N_{k}}^{N_{k}}}b_{l}^{k}(M)a_{l}^{k}(D)e^{ilM} + C$; $ N_{k} \in
    \mathbb{N}$, $b_{l} \in CS( \mathbb{R})$, $a_{l} \in CS(\mathbb{R})$, $C
    \in \mathcal{K}_{\mathbb{R}}$ e $a(D - j)e^{ijM} = e^{ijM}a(D)$.

Com isto provamos que $A \in \mathcal{A}$ pertence a $J_{0}$ se, e só se,
$\dps{\lim_{j \rightarrow \infty}} \chi_{j} A + K_{j} = A$, $K_{j} \in
\mathcal{K}_{\mathbb{R}}$.

Podemos então concluir que $A \in \mathcal{A}$ pertence a $J_{0}$ se, e só se,
$\sigma_{A}(m, \pm \infty) = 0$ para todo $m \in S^1_{+} \cup S^1_{-}$, onde
$S^1_{+} = S^1 \times \{ + \infty\}$ e $S^1_{-} = S^1 \times \{ -
\infty\}$. Pela Proposição 3.4 de \cite{[SM]} temos que $$ \sup \{ | \sigma_{A} (m , \xi)| : x
\in S^1_{+} \cup S^1_{-}, \: \xi = \pm \infty \}  \leq || \gamma_{A}||.$$ Logo
se $A \in \mbox{Ker} \, \gamma$ então
$\sigma_{A}(m, \pm \infty) = 0$ para $m \in S^1_{+} \cup S^1_{-}$. E vimos
acima que isto ocorre se, e só se, $A \in J_{0}$, portanto $\mbox{Ker} \, \gamma
= J_{0}$.$\cqd$ 
\chapter{A álgebra $\mathcal{C}$}
\hspace{8mm} Neste capítulo vamos estudar a $C^*$- subálgebra de $\mathcal{L}(L^2(\mathbb{R}))$, que chamamos $\mathcal{C}$, gerada por operadores de multiplicação $a(M)$ e $b(D) := F^{-1}b(M)F$ com $a, b \in CS(\mathbb{R})$.

Calcularemos sua K - teoria e obteremos uma fórmula para o índice de Fredholm.

\begin{Def}
Seja $\mathcal{C}$ a $C^*$- subálgebra de $\mathcal{L}(L^2(\mathbb{R}))$ obtida como o fecho da $C^*$- álgebra gerada por:\\
i) operadores de multiplicação $a(M)$, $a \in CS(\mathbb{R})$.\\
ii) operadores da forma $b(D) := F^{-1}b(M)F$, $b \in CS(\mathbb{R})$.
\end{Def}

É óbvio que $\mathcal{C}$ é uma $C^*$- subálgebra de $\mathcal{A}$,
introduzida no Capítulo 1.
\section{O espaço símbolo de $\mathcal{C}$}
\hspace{8mm} Nesta seção calcularemos o espaço símbolo, $M_{C}$, de $\mathcal{C}$ e daremos
uma fórmula explícita para o homomorfismo $\varphi: \mathcal{C}
\longrightarrow C(M_{C})$. 
\begin{Prop}
A $C^*$- álgebra $\mathcal{C}$ é irredutível.
\end{Prop}
\textbf{Demonstração:}

Temos que $\{ au : a \in \mathcal{C} \}$ é denso em
$L^2(\mathbb{R})$, para todo $u$ não nulo. Pois do Lema 1.13 temos que $\{ au : a \in
J_{0} \}$ é denso em $L^2(\mathbb{R})$, para todo $u$ não nulo e $J_{0}
\subset \mathcal{C}$. Então pelo Lema 1.2 temos que
$\mathcal{C}$ é irredutível$\cqd$

Usando o Teorema 1.3, (1.1) e a Observação 1.6 1/2 temos que
\begin{equation}
\mathcal{K}_{\mathbb{R}} \subset \mathcal{C}.
\end{equation}
Seja $J$ o ideal comutador de $\mathcal{C}$, temos 
\begin{equation}
J \subset \mathcal{K}_{\mathbb{R}}
\end{equation}
já que $[a(M), b(D)] \in \mathcal{K}_{\mathbb{R}}$, $a, b \in CS(
\mathbb{R})$. Por (2.1), (2.2) e pela Observação 1.6 1/2 temos $J = \mathcal{K}_{\mathbb{R}}.$
\begin{Teo}
O espaço símbolo $M_{C}$ da $C^*$- álgebra $\mathcal{C}$ é dado como \\
\centerline{$M_{C} = \{ (x, \xi) \in [ - \infty, + \infty] \times [ - \infty, + \infty] : |x| + |\xi| = \infty \}$.}

O $\sigma$- símbolo é dado por:\\
\centerline{$ \sigma_{a}(x,\xi) = a(x)$, $(x,\xi) \in M_{C}$,}\\[.2cm]
\centerline{$ \sigma_{b(D)}(x,\xi) = b(\xi)$, $(x,\xi) \in M_{C}$.}
\end{Teo}
\textbf{Demonstração:}

Consideremos as aplicações
\begin{displaymath}
\begin{array}{cccc}
 i_{1}: & CS( \mathbb{R}) & \longrightarrow &
\mathcal{C}/ \mathcal{K}_{\mathbb{R}}\\
        &  a              & \longmapsto     &
[a(M)]_{\mathcal{K}_{\mathbb{R}}}\\
\end{array}
\end{displaymath}
\centerline{e}
\begin{displaymath}
\begin{array}{cccc}
 i_{2}: & CS( \mathbb{R}) & \longrightarrow & \mathcal{C}/
 \mathcal{K}_{\mathbb{R}}\\
        &   b             & \longmapsto     &
 [b(D)]_{\mathcal{K}_{\mathbb{R}}}.\\
\end{array}
\end{displaymath}
Note que $i_{j}$, $j = 1, 2$, é injetora já que $a(M) \in
\mathcal{K}_{\mathbb{R}}$ se, e só se, $a \equiv 0$.

Temos que o espectro de $CS(\mathbb{R})$
é o conjunto $[- \infty, + \infty]$. Consideremos as
 aplicações duais  de $i_{1}$ e $i_{2}$: $$i^{*}_{1}: M_{C} \longrightarrow [
- \infty, + \infty ] \qquad \mbox{e}$$ $$i^{*}_{2}: M_{ C} \longrightarrow [
- \infty, + \infty ],$$
onde dado $w \in M_{C}$ temos $i^{*}_{j}(w) = w \circ i_{j}$, $j = 1$, $2$. Como $i_{j}$, $j = 1, 2$ é injetora temos que $i^{*}_{j}$, $j = 1, 2$ é
sobrejetora (veja, por exemplo, [3], pág. 314).

Seja $i$ o produto das duais $$i : M_{C} \longrightarrow  [ - \infty, + \infty ] \times [
- \infty, + \infty ],$$ onde $w \in M_{C}$ tem imagem $i(w) = ( w \circ i_{1}, w \circ i_{2})$.

Como $\mbox{Im} \;i_{1}$ e $\mbox{Im} \;i_{2}$ geram $\mathcal{C}/
\mathcal{K}_{\mathbb{R}}$ é fácil ver que $i$ é injetora. Como $M_{C}$ é compacto, $[ - \infty, + \infty ] \times [
- \infty, + \infty ]$ é Hausdorff e $i$ é contínua temos que $i$ é um homeomorfismo sobre a imagem.

Usando este homeomorfismo como identificação temos 
$$\sigma_{a}(x,\xi) = a(x) \qquad \quad \mbox{e} \qquad \quad \sigma_{b(D)}(x,\xi) = b(\xi).$$

Sabemos que $A = a(M)b(D) \in \mathcal{K}_{\mathbb{R}} = \mbox{Ker}\, \sigma$ para
$a,b \in C_{0}(\mathbb{R})$ e podemos escolher $a$ e $b$ que nunca se anulam em
$\mathbb{R}$. Como $$\sigma_{A}(x,\xi) = a(x)b(\xi) = 0 \quad \forall \quad (x,\xi) \in M_{C},$$
e $a(x)b(\xi) = 0$ se, e somente se, $|x| = \infty$ ou $|\xi| = \infty$ ou ambos,
segue que $M_{C} \subset \{ (x,\xi) \in [ - \infty, + \infty ] \times [ - \infty, + \infty] : |x| + |\xi| = \infty \}$.

Seja $w = (x, \xi) \in M_{C}$ e vamos definir
$\lambda([A]_{\mathcal{K}_{\mathbb{R}}}) =
w([F^{-1}AF]_{\mathcal{K}_{\mathbb{R}}})$. Temos que
$\lambda([a(M)]_{\mathcal{K}_{\mathbb{R}}}) =
w([a(D)]_{\mathcal{K}_{\mathbb{R}}}) = a( \xi)$ e $\lambda([b(D)]_{\mathcal{K}_{\mathbb{R}}}) =
w([b(-M)]_{\mathcal{K}_{\mathbb{R}}}) = b(-x)$, logo $\lambda = (\xi, x)$.

Usando o operador translação do Lema 1.22 temos que dado $x \in
\mathbb{R}$, $(x, + \infty) \in M_{C}$
se, e só se, $(y, + \infty) \in M_{C}$ para todo $y \in \mathbb{R}$. E usando o operador paridade do Lema
1.21 temos que dado $x \in \mathbb{R}$, $(x, + \infty) \in M_{C}$ se, e só se, $(-x, - \infty) \in
M_{C}$. Usando a sobrejetividade das aplicações duais, temos  $M_{C} = \{ (x,\xi) \in [ - \infty, + \infty ] \times [ - \infty, + \infty] : |x| + |\xi| = \infty \}$.$\cqd$ 

\section{Espaços de Hardy e operadores de Toeplitz}
\hspace{8mm} Nesta seção apresentaremos a definição de espaços de Hardy no
círculo, $H^2(S^1)$, e na
reta, $H^2(\mathbb{R})$, e provaremos que $U(H^2(S^1)) = H^2( \mathbb{R})$ (com $S^1 = \mathbb{R}
/ 2 \pi \mathbb{Z}$), onde $U$ é um operador unitário. Esta teoria será necessária na próxima seção.

Sejam $D = \{ z \in \mathbb{C}: |z| < 1 \}$ e $H(D)$ a classe de todas as funções analíticas em $D$.\\
Para $f \in H(D)$ definimos $M_{2}(f;r) = \Big\{ \frac{1}{2 \pi} \int^{\pi}_{ - \pi} |f(r e^{i t})|^2 \, d t \Big\}^{1/2}$, pelo Teorema 17.6 de \cite{[Rud]} temos que $M_{2}$ é uma função monótona crescente de $r \in [0,1)$, logo podemos definir para $f \in H(D)$
$$||f||_{2} = \dps{\lim_{r \rightarrow 1}}\: M_{2}(f;r).$$
\begin{Def}
$H^2(S^1) = \{ f \in H(D) : ||f||_{2} < \infty \}$.
\end{Def}
Aplicando a desigualdade de Minkowski para $M_{2}(f;r)$ é fácil verificar que
$||f||_{2}$ satisfaz a desigualdade triangular, logo $H^2(S^1)$ é um espaço
linear normado. Mais ainda, temos que $H^2(S^1)$ é um espaço de Banach, veja
por exemplo na página 331 de \cite{[Rud]} uma demonstração de que $H^2(S^1)$ é
completo.

Veremos a seguir que $H^2(S^1)$ pode ser identificado com um subespaço de $L^2(S^1)$, $S^1 = \{e^{i
  \theta} : \theta \in \mathbb{R} \}$ munido da medida $\frac{1}{2 \pi} d \, \theta$.
\begin{Teo}
a) Uma função $f \in H(D)$, da forma $$f(z) = \dps{\sum_{n=0}^{\infty}}a_{n}z^n
\qquad (z \in D),$$ está em  $H^2(S^1)$ se, e somente se,
$\dps{\sum_{n=0}^{\infty}}|a_{n}|^2 < \infty$; neste caso, $$||f||_{2} =
\{\dps{\sum_{n=0}^{\infty}}|a_{n}|^2 \}^2.$$
b) Se $f \in H^2(S^1)$, então $f$ tem limite radial $f^*(e^{it}) =
\dps{\lim_{r \rightarrow 1}} f(r e^{it})$ em quase
todo ponto de $S^1$; $f^* \in L^2(S^1)$; o n-ésimo coeficiente de Fourier de
$f^*$ é $a_{n}$ se $n \geq 0$ e zero se $n < 0$; $$\dps{\lim_{r \rightarrow
    1}} \, \frac{1}{2 \pi} \dps{\int_{- \pi}^{\pi}} |f^*(e^{it}) - f(re^{it})|^2 \,
dt = 0.$$ Além disso $f$ é a integral de Poisson de $f^*$, isto é, se $z = re^{i \theta}$ então
$$f(z) = \frac{1}{2 \pi} \dps{\int_{- \pi}^{\pi}}P_{r}(\theta - t)f^*(e^{i
  t}) \, dt$$ onde $P_{r}(\theta -t) $ é o núcleo de Poisson.\\
c) A aplicação $f \mapsto f^*$ é uma isometria de $H^2(S^1)$ sobre o subespaço
  de $L^2(S^1)$ que consiste naquelas $g \in L^2(S^1)$ que tem $\hat{g} (n) =
  0$ para todo $n < 0$.
\end{Teo}
\textbf{Demonstração:} Ver \cite{[Rud]}, Teorema 17.10. $\cqd$ 
\begin{Def}
$H^2(\mathbb{R})$ é o subespaço de $L^2( \mathbb{R})$ consistindo de todas as funções $f$ cuja extensão harmônica $F(x + iy) = \frac{1}{\pi} \int_{\mathbb{R}} \frac{y f(t)}{y^2 + (x - t)^2} \, dt$, é analítica no semi - plano superior.
\end{Def}
Enunciaremos agora um Teorema de Paley e Wiener que nos mostra que $F H^2(\mathbb{R}) = L^2( \mathbb{R}_{+})$.
\begin{Teo}
Uma função $f \in L^2( \mathbb{R})$ pertence a $H^2( \mathbb{R})$ se, e somente se, $\hat{f}$ é nula em quase todo o ponto no eixo real negativo.
\end{Teo}
\textbf{Demonstração:} Ver \cite{[He]}, pág. 101.$\cqd$

Provaremos agora um resultado que relaciona $H^2(S^1)$ com $H^2( \mathbb{R})$.
\begin{Prop}
Seja $U : L^2(S^1) \longrightarrow L^2( \mathbb{R})$ onde $g \longmapsto (Ug)(t) = \frac{1}{\sqrt{\pi} (1 - it)} g ( \frac{1 +it}{1 - it} )$, $t \in \mathbb{R}$. Então $U$ é um operador unitário de $L^2(S^1)$ para $L^2( \mathbb{R})$ tal que $U(H^2(S^1)) = H^2( \mathbb{R})$.
\end{Prop}
\textbf{Demonstração:}

Fazendo alguns cálculos, concluímos que $U$ é inversível e  $U^{-1}: L^2(
\mathbb{R}) \longrightarrow L^2(S^1)$ é dada por $h  \longmapsto (U^{-1}h)(z) = \frac{2 \sqrt{\pi}}{z + 1} h ( \frac{i (z - 1)}{-z -1} )$.

Vamos verificar que $U$ é uma isometria. Temos que $$||Ug||^{2}_{L^2( \mathbb{R})} = \frac{1}{ \pi} \int_{\mathbb{R}} \frac{1}{|1 - it|^2} \Big|g ( \frac{1 + it}{1 - it} ) \Big|^2 \, dt.$$

Seja $e^{i \theta} = \frac{1 + it}{1 - it}$ então $e^{i \theta} \,d \theta = \frac{2}{(1 -it )^2} dt$ sendo assim $d \theta = \frac{2}{|1 -it|^2} \,dt$. Substituindo esses valores obtemos que $||Ug||^{2}_{L^2( \mathbb{R})} = ||g||^{2}_{L^2(S^1)}$.

Verificaremos agora que $U(H^2(S^1)) = H^2( \mathbb{R})$. A aplicação conforme
$z = \frac{1 + iw}{1 -iw}$ leva o semi - plano superior sobre o disco
unitário, com $i$ sendo levado no zero e o zero sendo levado no um.  Sejam
$z_{1} = re^{i \theta}$ e $z_{2} = re^{it}$, $r \in [0, 1)$, temos que
$P_{r}(\theta -t) = \mbox{Re} \Big[ \frac{e^{it} + z_{1}}{e^{it} - z_{1}} \Big]$. Como a aplicação conforme é bijetora temos que
existem $w_{1} = x_{1} + i y_{1}$ e $w_{2} = x_{2} + i y_{2}$ tais que $z_{1}
= \frac{1 + iw_{1}}{1 -iw_{1}}$ e $z_{2} = \frac{1 + iw_{2}}{1 -iw_{2}}$. Na
fronteira temos $e^{i t} = \frac{1 + ix_{2}}{1 - ix_{2}}$ logo $\frac{1}{2 \pi}d
t = \frac{1}{\pi} \frac{1}{(1 + x_{2}^2)}dx_{2}$. Compondo a fórmula do núcleo
de Poisson com a aplicação conforme obtemos $\frac{y_{1}(1 + x_{2}^2)}{y_{1}^2
  + (x_{1} - x_{2})^2}$.

Dada $f \in H^2(S^1)$ consideremos $\tilde{f}(z) = f(z) \frac{1 + z}{2
  \sqrt{\pi}}$. Temos que $\tilde{f} \in H^2(S^1)$. Então pelo Teorema 2.5 item
  b temos que $\tilde{f}(z_{1}) = \frac{1}{2 \pi} \dps{\int_{-
  \pi}^{\pi}}P_{r}(\theta -t)\tilde{f}^*(e^{i
  t}) \, dt $. Compondo com a aplicação conforme temos $\tilde{f}(\frac{1 +
  iw_{1}}{1 -iw_{1}}) = \frac{1}{\pi}
\int_{\mathbb{R}} \frac{y_{1}}{y_{1}^2 + (x_{1} - x_{2})^2}
  \tilde{f}^*(\frac{1 + ix_{2}}{1 - ix_{2}})\, dx_{2} =$\\$= \frac{1}{\pi}
\int_{\mathbb{R}} \frac{y_{1}}{y_{1}^2 + (x_{1} - x_{2})^2}
  \frac{1}{\sqrt{\pi}(1-ix_{2})}f^*(\frac{1 + ix_{2}}{1 - ix_{2}})\, dx_{2} =
  \frac{1}{\pi} \int_{\mathbb{R}} \frac{y_{1}}{y_{1}^2 + (x_{1} - x_{2})^2}
  (Uf^*)(x_{2}) d \, x_{2}$ que é
  analítica no semi-plano superior. 

Portanto $U(H^2(S^1)) \subset H^2( \mathbb{R})$. A demonstração da volta
utiliza a mesma estratégia.$\cqd$ 
\begin{Def}
Seja $P$ a projeção ortogonal de $L^2(S^1)$ sobre $H^2(S^1)$. Se $\varphi \in
L^{\infty}(S^1)$, o operador $$T_{\varphi}: H^2(S^1) \longrightarrow H^2(S^1),
\qquad \varphi \longmapsto P( \varphi f)$$ é chamado operador de Toeplitz com
símbolo $\varphi$.
\end{Def}
Enunciaremos agora um teorema  que será utilizado na próxima seção.
\begin{Teo}
 Se $\varphi \in C(S^1)$, então $T_{\varphi}$ é um operador de Fredholm se, e
 só se, $\varphi$ nunca se anula. Seja $\varphi$ um elemento inversível em $C(S^1)$. Então o índice de Fredholm
de $T_{\varphi}$ é menos o número de rotação de $\varphi$, isto é,
$$ind(T_{\varphi}) = - w( \varphi).$$
\end{Teo}
\textbf{Demonstração:} Ver \cite{[Murphy]}, pág. 103 e 104.$\cqd$
\section{Índice de Fredholm e a K-teoria de $\mathcal{C}$}
\hspace{8mm} Nesta seção e nos capítulos  que seguem estudaremos a K-teoria de $C^*$-
álgebras. Como referência de um amplo estudo sobre K-teoria para $C^*$-
álgebras recomendamos os livros de B. Blackadar, \emph{K-Theory for Operator Algebras} \cite{[Bk]} e
M. Rordam, F. Larsen e N. J. Laustisen, \emph{An Introduction to K-Theory for $C^*$-
algebras} \cite{[R]}. Para um breve resumo sobre K-teoria para $C^*$- álgebras
recomendamos o Apêndice C de \cite{[Cris]}. Utilizaremos neste trabalho as notações e resultados de
\cite{[Bk]} e \cite{[R]}.

A K-teoria foi desenvolvida por Atiyah e Hirzebruch por volta de 1960, baseada
no trabalho de Grothendieck em geometria algébrica. Ela foi introduzida como
uma ferramenta na teoria de $C^*$- álgebras por volta de 1970. 

Resumidamente, K-teoria (para $C^*$- álgebras) é um par de funtores, chamados $K_{0}$ e
$K_{1}$, que a cada $C^*$- álgebra $A$ associa dois grupos abelianos
$K_{0}(A)$ e $K_{1}(A)$. Temos que $K_{0}$ é um grupo de diferenças formais de
classes de projeções (matrizes) de tamanho arbitrário e $K_{1}$ é o grupo de
classes de homotopias de matrizes unitárias de tamanho arbitrário.
Uma das propriedades dos funtores $K_{0}$ e $K_{1}$ é que dada uma sequência exata curta de $C^*$- álgebras
 \begin{center}
\begin{picture}(280,15)
\put(15,0){0}
\put(30,3){\vector(1,0){25}}
\put(43,0){\makebox[60pt]{$A$}}
\put(90,3){\vector(1,0){35}}
\put(110,0){\makebox[60pt]{$B$}}
\put(152,3){\vector(1,0){35}}
\put(181,0){\makebox[60pt]{$C$}}
\put(230,3){\vector(1,0){25}}
\put(262,0){0}
\put(90,8){\makebox[35pt]{$\varphi $}}
\put(152,8){\makebox[35pt]{$\psi$}}
\end{picture}
\end{center}
podemos associar uma sequência exata cíclica de seis termos de grupos
abelianos

\begin{picture}(350,15)
\put(48,0){\makebox[60pt]{$K_{0}(A)$}}
\put(111,3){\vector(1,0){35}}
\put(149,0){\makebox[60pt]{$K_{0}(B)$}}
\put(212,3){\vector(1,0){35}}
\put(250,0){\makebox[60pt]{$K_{0}(C)$}}
\put(111,8){\makebox[35pt]{$ \varphi_{*}$}}
\put(212,8){\makebox[35pt]{$ \psi_{*} $}}
\end{picture}

\begin{picture}(350,45)
\put(78,0){\vector(0,1){35}}
\put(280,35){\vector(0,-1){35}}
\put(81,17){$\delta_{1}  $}
\put(283,17){$ \delta_{0} $}
\end{picture}

\begin{picture}(350,15)
\put(48,0){\makebox[60pt]{$K_{1}(C)$}}
\put(146,3){\vector(-1,0){35}}
\put(149,0){\makebox[60pt]{$K_{1}(B)$}}
\put(247,3){\vector(-1,0){35}}
\put(250,0){\makebox[60pt]{$K_{1}( A)$}}
\put(111,8){\makebox[35pt]{$ \psi_{*}$}}
\put(212,8){\makebox[35pt]{$ \varphi_{*} $}}
\end{picture}\\
onde as setas horizontais são homomorfismos induzidos funtorialmente e a seta
de conexão de $K_{1}$ para $K_{0}$ é chamada aplicação do índice e a seta de
conexão de $K_{0}$ para $K_{1}$ é chamada aplicação exponencial.

Vimos na seção 1 que $\mathcal{C} / \mathcal{K}_{\mathbb{R}}$ é isomorfo a
$C(M_{C})$. Usando que $M_{C}$ é homeomorfo a $S^1$ temos que $K_{i}(\mathcal{C} /
\mathcal{K}_{\mathbb{R}}) \cong K_{i}(C(S^1))$, $i = 0,1$. Pelo exemplo
11.3.3 de \cite{[R]} sabemos que
$K_{i}(C(S^1)) \cong \mathbb{Z}$, $i = 0,1$, onde a classe $[.]_{1}$ da
função $f: S^1 \longrightarrow \mathbb{C}$ definida por $f(z) = z$ gera
$K_{1}(C(S^1))$ e a classe $[.]_{0}$ da função $g: S^1 \longrightarrow
\mathbb{C}$ definida por $g(z) = 1$ gera $K_{0}(C(S^1))$. Pelo
exemplo 11.3.4 de \cite{[R]} temos que 
\begin{displaymath}
\begin{array}{cccc}
\triangle: &  K_{1}(C(S^1)) & \longrightarrow & \pi^1(S^1)\\
           &  [u]_{1}       & \longmapsto     & <Det(u)>,\\
\end{array}
\end{displaymath}
onde $\pi^1(S^1)$ é o grupo fundamental de $S^1$ e $u$ é um
unitário é um isomorfismo. Pelo exemplo 8.3.2 de \cite{[R]} temos que
\begin{displaymath}
\begin{array}{cccc}
w: &  \pi^1(S^1) & \longrightarrow & \mathbb{Z}\\
           &  < \lambda >       & \longmapsto     & w(\lambda),\\
\end{array}
\end{displaymath}
onde $ \lambda \in C(S^1, S^1)$ é um isomorfismo e satisfaz a seguinte
propriedade $u $ é homotópico a $v$ se, e só se, $w(u) = w(v)$, $u$ e $v \in C(S^1, S^1)$.

Logo 
\begin{displaymath}
\begin{array}{cccc}
w \circ \triangle: &  K_{1}(C(S^1)) & \longrightarrow & \mathbb{Z} \\
           &  [u]_{1}       & \longmapsto     & w(Det(u))\\
\end{array}
\end{displaymath}
é um isomorfismo. Com isto temos que a classe de qualquer $f \in C(S^1, S^1)$ que tenha $w(f)
=1$ ou $-1$ gera $K_{1}(C(S^1))$.
 
Fixando um homeomorfismo $h$ de $M_{C}$ para $S^1$ temos que $[1]_{0}$ gera
$K_{0}(C(M_{C}))$ e $[h]_{1}$ gera $K_{1}(C(M_{C}))$.

Usando  o $\sigma$-símbolo temos que $[[Id]_{\mathcal{K}_{\mathbb{R}}}]_{0}$
gera $K_{0}(\mathcal{C} /
\mathcal{K}_{\mathbb{R}})$ e $[[e^{2 \pi i c(M)} b(D) +
c(D)]_{\mathcal{K}_{\mathbb{R}}}]_{1}$ onde $b, c \in CS( \mathbb{R})$ são tais que $b$ é crescente com $b(x) = 0$, se $x \leq
-1$ e $b(x) = 1$, se $x \geq 1$ e $c$ é decrescente com $c(x) = 0$, se $x \geq
1$ e $c(x) = 1$, se $x \leq -1$, $b + c \equiv 1$ gera $K_{1}(\mathcal{C} /
\mathcal{K}_{\mathbb{R}})$, isto é, $$K_{0}(\mathcal{C} /
\mathcal{K}_{\mathbb{R}}) = \mathbb{Z}[[Id]_{\mathcal{K}_{\mathbb{R}}}]_{0}
\hspace{0.5cm} \mbox{e}$$ $$K_{1}(\mathcal{C} /
\mathcal{K}_{\mathbb{R}}) = \mathbb{Z}[[e^{2 \pi i c(M)} b(D) +
c(D)]_{\mathcal{K}_{\mathbb{R}}}]_{1}.$$
\begin{Lema} 
Seja $T \in \mathcal{C}$. Então $T$ é um operador de Fredholm se, e somente se,
$\sigma_{T}$ é inversível.
\end{Lema}
\textbf{Demonstração:}

Como $\mathcal{C}/ \mathcal{K}_{\mathbb{R}}$ é uma $C^*$- subálgebra de
$\mathcal{L}(L^2(\mathbb{R})) / \mathcal{K}_{\mathbb{R}}$, segue do Teorema
2.1.11 de \cite{[Murphy]} que $[T]_{\mathcal{K}_{\mathbb{R}}} \in \mathcal{C}
/ \mathcal{K}_{\mathbb{R}}$ é inversível se, e somente se,
$[T]_{\mathcal{K}_{\mathbb{R}}} \in \mathcal{L}(L^2(\mathbb{R})) /
\mathcal{K}_{\mathbb{R}}$ é inversível. Pelo Teorema 1.4.16 de \cite{[Murphy]}
temos que  $[T]_{\mathcal{K}_{\mathbb{R}}} \in \mathcal{L}(L^2(\mathbb{R})) /
\mathcal{K}_{\mathbb{R}}$ é inversível se, e somente se, $T$ é de Fredholm. Por
outro lado, segue do Teorema 2.3 que $[T]_{\mathcal{K}_{\mathbb{R}}} \in
\mathcal{C}/ \mathcal{K}_{\mathbb{R}}$ é inversível se, e só se, $\sigma_{T}(x,
\xi) \neq 0$ para todo $(x, \xi) \in M_{C}$.$\cqd$

Sejam $\chi$ a função característica da semi-reta real positiva e $\phi \in CS(\mathbb{R})$ tal que $\phi( + \infty) = \phi( - \infty) = 1$, $\phi$ não nula e $w( \phi) = -1$.
Pelo Teorema 2.7 temos que $\chi(D) \phi(M) \chi(D)$ é equivalente a um
operador de Toeplitz em $H^2(\mathbb{R})$, logo unitariamente equivalente a um
operador de Toeplitz em $H^2(S^1)$, pela Proposição 2.8. Então pelo Teorema 2.10 temos que $\chi(D) \phi(M) \chi(D)$ é um operador de Fredholm com índice 1. Logo $\chi(D) \phi(M) \chi(D) \oplus Id_{H^2(\mathbb{R})^{\bot}} = \chi(D) \phi(M) \chi(D)\, +\, (Id - \chi(D))$ é um operador de Fredholm com índice 1.
\begin{Lema}
O operador $(Id - \chi(D)) \phi(M) \chi(D)$ é compacto.
\end{Lema}
\textbf{Demonstração:}

Pela Proposição 2.8 e pelo Teorema 2.7 temos que $U^{-1} \chi(D) U = P$ onde
$P$ é a projeção ortogonal em $H^2(S^1)$. Temos também que $U^{-1} \phi(M) U =
\tilde{\phi}(M)$ onde $\tilde{\phi}(z) = \phi \big(\frac{z-1}{i(z+1)} \big)$.

Suponha $\phi(t) = 1$ se $|t| > 1$, então $\tilde{\phi}(z) \equiv 1$ em uma
vizinhança de $-1$. Portanto $\tilde{\phi} \in C^{\infty}(S^1)$. Logo,
$U^{-1}( Id - \chi(D)) \phi(M) \chi(D) U = (Id - P) \tilde{\phi}(M)P$ com $P$ a
projeção sobre $H^2(S^1)$ e $\tilde{\phi} \in  C^{\infty}(S^1)$.

Consideremos $(e_{k})_{k \in \mathbb{Z}}$ base ortonormal de $L^2(S^1)$ onde $e_{k}(z) =
z^{k}$. Temos que $$(Id
-P)e_{n}(M)Pe_{k} = 0 \quad \mbox{se} \quad k < 0 \quad \mbox{e} \quad (Id
-P)e_{n}(M)Pe_{k} = (Id - P)e_{k+n} \quad \mbox{se} \quad k \geq 0.$$
Mas $(Id - P)e_{k+n} = 0$ se $k + n \geq 0$. Vamos estudar então o
que acontece para $n \geq 0$ e $n < 0$.

Se $n \geq 0$ então $k + n \geq 0$ para todo $k \geq 0$, logo $ (Id - P)e_{k+n}
= 0$.

Se $n <0$ teremos $k + n \geq 0$ se, e somente se, $k \geq -n$. Então
$(Id - P)e_{k+n} = 0$ para todo $k \geq -n$ e $(Id - P)e_{k+n} = e_{k+n}$ caso
contrário. Logo para n fixo $\mbox{Im}\, [(Id-P)e_{n}P]$ é gerada por  $\{
e_{k} : 0 \leq k \leq -n\}$.

Com isto temos que $(Id-P)p_{n}(M)P$ é um operador de posto finito, logo
compacto.
Logo se $p$ é um polinômio temos que $(Id-P)p(M)P$ é compacto. Pelo Teorema de
Stone-Weierstrass temos que existem  polinômios $p_{n}$ convergindo
uniformemente para $\tilde{\phi}$ então
$$||(Id-P)\tilde{\phi}(M)P - (Id-P)p_{n}(M)P|| \leq ||\tilde{\phi}(M) -
p_{n}(M)|| = ||\tilde{\phi} - p_{n}||_{\infty} \longrightarrow 0.$$

Como $\mathcal{K}(L^2(S^1 ))$ é fechado em $\mathcal{L}(L^2(S^1
))$ temos que $(Id-P)\tilde{\phi}(M)P$ é compacto em $L^2(S^1)$. Logo  $(Id -
\chi(D)) \phi(M) \chi(D)$ é compacto em $L^2(\mathbb{R})$.$\cqd$
\begin{Teo}
Seja $T \in \mathcal{C}$ um operador de Fredholm, temos que $ind(T) =
\delta_{1}([[T]_{\mathcal{K}_{ \mathbb{R}}}]_{1}) = w(\sigma_{T})$, onde
$\delta_{1}$ é a aplicação do índice da sequência de K-teoria associada a
sequência exata curta abaixo.
\begin{center}
\begin{picture}(280,15)
\put(15,0){0}
\put(30,3){\vector(1,0){25}}
\put(43,0){\makebox[60pt]{$\mathcal{K}_{\mathbb{R}}$}}
\put(90,3){\vector(1,0){35}}
\put(110,0){\makebox[60pt]{$\mathcal{C}$}}
\put(152,3){\vector(1,0){35}}
\put(181,0){\makebox[60pt]{$\mathcal{C}/ \mathcal{K}_{\mathbb{R}}$}}
\put(230,3){\vector(1,0){25}}
\put(262,0){0}
\put(90,8){\makebox[35pt]{$i $}}
\put(152,8){\makebox[35pt]{$\pi$}}
\end{picture}
\end{center}
\end{Teo}
\textbf{Demonstração:}

Vamos considerar as seguintes sequências exatas curtas
\begin{center}
\begin{picture}(280,15)
\put(15,0){0}
\put(30,3){\vector(1,0){25}}
\put(43,0){\makebox[60pt]{$\mathcal{K}_{\mathbb{R}}$}}
\put(90,3){\vector(1,0){35}}
\put(110,0){\makebox[60pt]{$\mathcal{C}$}}
\put(152,3){\vector(1,0){35}}
\put(181,0){\makebox[60pt]{$\mathcal{C}/ \mathcal{K}_{\mathbb{R}}$}}
\put(230,3){\vector(1,0){25}}
\put(262,0){0}
\put(90,8){\makebox[35pt]{$i $}}
\put(152,8){\makebox[35pt]{$\pi$}}
\end{picture}
\end{center}
\begin{center}
\begin{picture}(395,15)
\put(10,0){0}
\put(30,3){\vector(1,0){25}}
\put(48,0){\makebox[60pt]{$\mathcal{K}_{\mathbb{R}}$}}
\put(101,3){\vector(1,0){35}}
\put(149,0){\makebox[60pt]{$\mathcal{L}(L^2( \mathbb{R}))$}}
\put(212,3){\vector(1,0){35}}
\put(270,0){\makebox[60pt]{$\mathcal{L}(L^2( \mathbb{R}))/ \mathcal{K}_{\mathbb{R}} $}}
\put(353,3){\vector(1,0){25}}
\put(388,0){0}
\put(101,8){\makebox[35pt]{$ i$}}
\put(212,8){\makebox[35pt]{$\pi $}}
\end{picture}
\end{center}
onde $i$ é a inclusão e $\pi$ a projeção canônica.

Sabemos pela Proposição 9.4.2 de \cite{[R]} que o índice, $\bar{\delta}_{1}$,
para a sequência exata curta 
\begin{center}
\begin{picture}(395,15)
\put(10,0){0}
\put(30,3){\vector(1,0){25}}
\put(48,0){\makebox[60pt]{$\mathcal{K}_{\mathbb{R}}$}}
\put(101,3){\vector(1,0){35}}
\put(149,0){\makebox[60pt]{$\mathcal{L}(L^2( \mathbb{R}))$}}
\put(212,3){\vector(1,0){35}}
\put(270,0){\makebox[60pt]{$\mathcal{L}(L^2( \mathbb{R}))/ \mathcal{K}_{\mathbb{R}} $}}
\put(353,3){\vector(1,0){25}}
\put(388,0){0}
\put(101,8){\makebox[35pt]{$ i$}}
\put(212,8){\makebox[35pt]{$\pi $}}
\end{picture}
\end{center}
é o índice de Fredholm, isto é $ind(T) = \bar{\delta}_{1}([[T]_{\mathcal{K}_{ \mathbb{R}}}]_{1})$.

Consideremos o diagrama comutativo

\begin{picture}(350,15)
\put(10,0){0}
\put(20,3){\vector(1,0){25}}
\put(48,0){\makebox[60pt]{$\mathcal{K}_{\mathbb{R}}$}}
\put(111,3){\vector(1,0){35}}
\put(149,0){\makebox[60pt]{$\mathcal{C}$}}
\put(212,3){\vector(1,0){35}}
\put(256,0){\makebox[60pt]{$\mathcal{C}/ \mathcal{K}_{\mathbb{R}}$}}
\put(325,3){\vector(1,0){25}}
\put(353,0){0}
\put(111,8){\makebox[35pt]{$i $}}
\put(212,8){\makebox[35pt]{$\pi $}}
\end{picture} 

\begin{picture}(350,35)
\put(78,35){\vector(0,-1){35}}
\put(179,35){\vector(0,-1){35}}
\put(286,35){\vector(0,-1){35}}
\put(81,17){$ $}
\put(182,17){$i$}
\put(289,17){$i$}
\end{picture}

\begin{picture}(350,15)
\put(10,0){0}
\put(20,3){\vector(1,0){25}}
\put(48,0){\makebox[60pt]{$\mathcal{K}_{\mathbb{R}}$}}
\put(111,3){\vector(1,0){35}}
\put(149,0){\makebox[60pt]{$\mathcal{L}(L^2( \mathbb{R}))$}}
\put(212,3){\vector(1,0){35}}
\put(256,0){\makebox[60pt]{$\mathcal{L}(L^2( \mathbb{R})) / \mathcal{K}_{\mathbb{R}}$}}
\put(325,3){\vector(1,0){25}}
\put(353,0){0}
\put(111,8){\makebox[35pt]{$i $}}
\put(212,8){\makebox[35pt]{$\pi $}}
\end{picture} 
   
Pela naturalidade da aplicação do índice, veja Proposição 9.1.5 de \cite{[R]}, temos que o diagrama abaixo comuta 
\begin{displaymath}
\begin{array}{ccc}
K_{1}(\mathcal{C} / \mathcal{K}_{\mathbb{R}}) & {\buildrel \delta_{1} \over \longrightarrow} & \mathbb{Z}\\
\downarrow & & \Arrowvert\\
K_{1}(\mathcal{L}(L^2( \mathbb{R})) / \mathcal{K}_{\mathbb{R}}) & {\buildrel
  \bar{\delta_{1}} \over \longrightarrow} & \mathbb{Z}\\
\end{array}
\end{displaymath}

Logo $\delta_{1}$ é o índice de Fredholm, ou seja $ind(T) = \delta_{1}([[T]_{\mathcal{K}_{ \mathbb{R}}}]_{1})$. Vamos verificar agora que $ind(T) = w(\sigma_{T})$.

Consideremos $S = \chi(D) \phi(M) \chi(D) \,+\, (Id - \chi(D)) \,+\, (Id -
\chi(D)) \phi(M) \chi(D) = \phi(M) \chi(D) + (Id - \chi(D))$. Como vimos no
Lema 2.12 que $(Id - \chi(D)) \phi(M) \chi(D)$ é compacto, temos que $S$ é um
operador de Fredholm com índice 1, pela observação antes do Lema 2.12.

Seja $T' = \phi(M) b(D) + (Id - b(D))$ com $b$ contínua crescente com $b(x) =
0$ para $x \leq -1$ e $b(x) = 1$ para $x \geq 1$. Temos $T'\in \mathcal{C}$ e $T' - S =
(\phi -1)(M)(b - \chi)(D)$ logo $T'- S$ é compacto pela Observação 1.14
1/2. Portanto $T' = S + K$, $K$ compacto, logo $T'$ é Fredholm e $ind(T') = 1$ e temos também que $w(\sigma_{T'}) = 1$. Concluimos assim que $ind(T) =
\delta_{1}([[T]_{\mathcal{K}_{ \mathbb{R}}}]_{1}) = w(\sigma_{T})$, $T \in
\mathcal{C}$ um operador de Fredholm.$\cqd$
\begin{Prop}
Seja $\pi: \mathcal{C} \longrightarrow \mathcal{C} / \mathcal{K}_{\mathbb{R}}$
a projeção canônica. Então  $\pi_{*}: K_{0}( \mathcal{C}) \longrightarrow
K_{0}( \mathcal{C} / \mathcal{K}_{\mathbb{R}})$ é um isomorfismo. Além disso, $K_{1}( \mathcal{C}) = 0$.
\end{Prop}
\textbf{Demonstração:}

Consideremos a sequência exata curta 

\begin{center}
\begin{picture}(280,15)
\put(15,0){0}
\put(30,3){\vector(1,0){25}}
\put(43,0){\makebox[60pt]{$\mathcal{K}_{\mathbb{R}}$}}
\put(90,3){\vector(1,0){35}}
\put(110,0){\makebox[60pt]{$\mathcal{C}$}}
\put(152,3){\vector(1,0){35}}
\put(181,0){\makebox[60pt]{$\mathcal{C}/ \mathcal{K}_{\mathbb{R}}$}}
\put(230,3){\vector(1,0){25}}
\put(262,0){0}
\put(90,8){\makebox[35pt]{$i $}}
\put(152,8){\makebox[35pt]{$\pi$}}
\end{picture}
\end{center} 
 e sua sequência exata cíclica de seis termos 

\begin{picture}(350,15)
\put(48,0){\makebox[60pt]{$K_{0}(\mathcal{K}_{\mathbb{R}})$}}
\put(111,3){\vector(1,0){35}}
\put(149,0){\makebox[60pt]{$K_{0}(\mathcal{C})$}}
\put(212,3){\vector(1,0){35}}
\put(250,0){\makebox[60pt]{$K_{0}(\mathcal{C}/ \mathcal{K}_{\mathcal{R}})$}}
\put(111,8){\makebox[35pt]{$ i_{*}$}}
\put(212,8){\makebox[35pt]{$ \pi_{*} $}}
\end{picture}

\begin{picture}(350,45)
\put(78,0){\vector(0,1){35}}
\put(280,35){\vector(0,-1){35}}
\put(81,17){$\delta_{1}  $}
\put(283,17){$ \delta_{0} $}
\end{picture}

\begin{picture}(350,15)
\put(48,0){\makebox[60pt]{$K_{1}(\mathcal{C} / \mathcal{K}_{\mathbb{R}})$}}
\put(146,3){\vector(-1,0){35}}
\put(149,0){\makebox[60pt]{$K_{1}(\mathcal{C})$}}
\put(247,3){\vector(-1,0){35}}
\put(250,0){\makebox[60pt]{$K_{1}( \mathcal{K}_{\mathbb{R}})$}}
\put(111,8){\makebox[35pt]{$ \pi_{*}$}}
\put(212,8){\makebox[35pt]{$ i_{*} $}}
\end{picture}

Pelo corolário 6.4.2 e pelo exemplo 8.2.9 de \cite{[R]} temos que
$K_{0}(\mathcal{K}_{\mathbb{R}}) = \mathbb{Z}$ e
$K_{1}(\mathcal{K}_{\mathbb{R}}) = 0$. Substituindo então as K - teorias conhecidas temos

\begin{picture}(350,15)
\put(48,0){\makebox[60pt]{$\mathbb{Z}$}}
\put(111,3){\vector(1,0){35}}
\put(149,0){\makebox[60pt]{$K_{0}(\mathcal{C})$}}
\put(212,3){\vector(1,0){35}}
\put(250,0){\makebox[60pt]{$\mathbb{Z}$}}
\put(111,8){\makebox[35pt]{$ i_{*}$}}
\put(212,8){\makebox[35pt]{$ \pi_{*} $}}
\end{picture}

\begin{picture}(350,45)
\put(78,0){\vector(0,1){35}}
\put(280,35){\vector(0,-1){35}}
\put(81,17){$\delta_{1}  $}
\put(283,17){$ \delta_{0} $}
\end{picture}

\begin{picture}(350,15)
\put(48,0){\makebox[60pt]{$\mathbb{Z}$}}
\put(146,3){\vector(-1,0){35}}
\put(149,0){\makebox[60pt]{$K_{1}(\mathcal{C})$}}
\put(247,3){\vector(-1,0){35}}
\put(250,0){\makebox[60pt]{$0$}}
\put(111,8){\makebox[35pt]{$ \pi_{*}$}}
\put(212,8){\makebox[35pt]{$ i_{*} $}}
\end{picture}

Seja $T' = \phi(M)b(D) + (Id - b(D))$ o operador do Teorema 2.13 , temos que
$[[T']_{\mathcal{K}_{\mathbb{R}}}]_{1} \in K_{1}(\mathcal{C} /
\mathcal{K}_{\mathbb{R}})$ e já vimos que
$\delta_{1}([[T']_{\mathcal{K}_{\mathbb{R}}}]_{1}) = 1$, logo $\delta_{1}$ é
um isomorfismo. Portanto $\pi_{*}: K_{0}(\mathcal{C}) \longrightarrow
\mathbb{Z}$ é isomorfismo e $K_{1}(\mathcal{C}) = 0$.$\cqd$
\chapter{K-teoria de $\mathcal{A}$}
\hspace{8mm}Neste capítulo calcularemos a K-teoria da $C^*$- álgebra
$\mathcal{A}$ utilizando a seguinte estratégia. Primeiro calcularemos a
K-teoria de $\mathcal{A}/ \mathcal{E}$ e $\mathcal{E} /
\mathcal{K}_{\mathbb{R}}$ e para isto utilizaremos os isomorfismos 
 $\varphi:
\mathcal{A}/ \mathcal{E} \longrightarrow C(M_{A})$ e $ \gamma:
\mathcal{E}/ \mathcal{K}_{\mathbb{R}} \longrightarrow C(M_{SL},
\mathcal{K}_{\mathbb{Z}})$, estudados no
Capítulo 1. De posse dessas K-teorias e utilizando a seguinte sequência exata curta 
\begin{center}
\begin{picture}(350,15)
\put(10,0){0}
\put(20,3){\vector(1,0){25}}
\put(48,0){\makebox[60pt]{$\mathcal{E}/ \mathcal{K}_{\mathbb{R}}$}}
\put(111,3){\vector(1,0){35}}
\put(149,0){\makebox[60pt]{$\mathcal{A}/ \mathcal{K}_{\mathbb{R}}$}}
\put(212,3){\vector(1,0){35}}
\put(250,0){\makebox[60pt]{$ \mathcal{A}/ \mathcal{E}$}}
\put(313,3){\vector(1,0){25}}
\put(341,0){0}
\put(111,8){\makebox[35pt]{}}
\put(212,8){\makebox[35pt]{}}
\end{picture}
\end{center}
calcularemos a K-teoria  de $\mathcal{A}/
\mathcal{K}_{\mathbb{R}}$. Por último utilizamos a sequência exata
curta abaixo para o cálculo da K-teoria de $\mathcal{A}$
\begin{center}
\begin{picture}(280,15)
\put(15,0){0}
\put(30,3){\vector(1,0){25}}
\put(43,0){\makebox[60pt]{$\mathcal{K}_{\mathbb{R}}$}}
\put(90,3){\vector(1,0){35}}
\put(110,0){\makebox[60pt]{$\mathcal{A}$}}
\put(152,3){\vector(1,0){35}}
\put(181,0){\makebox[60pt]{$\mathcal{A}/ \mathcal{K}_{\mathbb{R}}$}}
\put(230,3){\vector(1,0){25}}
\put(262,0){0}
\put(90,8){\makebox[35pt]{}}
\put(152,8){\makebox[35pt]{}}
\end{picture}.
\end{center}
\section{A K-teoria de $\mathcal{A}/ \mathcal{E}$ e $\mathcal{E} /
\mathcal{K}_{\mathbb{R}}$}

\hspace{8mm}Começaremos o nosso trabalho por descrever geradores para $K_{i}(C(M_{\sharp}))$,
$i = 0,1$, já que pelo Teorema 1.17 $M_{A} = M_{ \sharp} \times \{ - \infty, +
\infty\}$ e pelo Teorema de Gelfand temos $\varphi:
\mathcal{A}/ \mathcal{E} \longrightarrow C(M_{A})$ isomorfismo . Pelo Lema 1.10
podemos considerar $M_{ \sharp} = S^{1}_{-} \cup \mathbb{R} \cup S^{1}_{+}$,
onde $S^{1}_{+} = S^{1} \times \{ + \infty\}$ e $S^{1}_{-} = S^{1} \times \{ - \infty\}$.
\begin{Prop}
Seja $\psi: C(M_{\sharp}) \longrightarrow C(S^1 \times \{-\infty, + \infty
\})$ a restrição das $f \in C(M_{\sharp})$ a $S^{1}_{+} \cup S^{1}_{-}$. Para $i = 0,1$
$\psi_{*}: K_{i}(C(M_{\sharp})) \longrightarrow K_{i}(C(S^1 \times \{-\infty,
+ \infty \}))$ é injetora. A imagem de $\psi_{*}: K_{0}(C(M_{\sharp}))
\longrightarrow K_{0}(C(S^1 \times \{-\infty,+ \infty \}))$ é isomorfa a
$\mathbb{Z}$ e a imagem de $\psi_{*}: K_{1}(C(M_{\sharp})) \longrightarrow
K_{1}(C(S^1 \times \{-\infty,+ \infty \}))$ é isomorfa a $\mathbb{Z} \oplus
\mathbb{Z}$. Assim $K_{0}(C(M_{\sharp})) \cong \mathbb{Z}$ e
$K_{1}(C(M_{\sharp})) \cong \mathbb{Z} \oplus \mathbb{Z}$.
\end{Prop}
\textbf{Demonstração:}

Consideremos a sequência exata curta
\begin{center}
\begin{picture}(395,15)
\put(10,0){0}
\put(24,3){\vector(1,0){25}}
\put(48,0){\makebox[60pt]{$C_{0}( \mathbb{R})$}}
\put(107,3){\vector(1,0){35}}
\put(143,0){\makebox[60pt]{$C(M_{ \sharp})$}}
\put(204,3){\vector(1,0){35}}
\put(262,0){\makebox[60pt]{$C(S^1 \times \{-\infty, + \infty \})$}}
\put(345,3){\vector(1,0){25}}
\put(382,0){0}
\put(107,8){\makebox[35pt]{$ \varphi$}}
\put(204,8){\makebox[35pt]{$\psi $}}
\end{picture}
\end{center}
onde $ \varphi$ é a inclusão e $\psi$ é a restrição.

Sabemos que esta sequência exata curta induz a sequência exata cíclica de seis
termos em K-teoria

\begin{picture}(395,15)
\put(21,0){\makebox[60pt]{$K_{0}(C_{0}( \mathbb{R}))$}}
\put(117,3){\vector(1,0){35}}
\put(160,0){\makebox[60pt]{$K_{0}(C(M_{ \sharp}))$}}
\put(228,3){\vector(1,0){35}}
\put(300,0){\makebox[60pt]{$K_{0}(C(S^1 \times \{-\infty, + \infty \}))$}}
\put(117,8){\makebox[35pt]{$ \varphi_{*}  $}}
\put(228,8){\makebox[35pt]{$ \psi_{*} $}}
\end{picture}

\begin{picture}(350,45)
\put(49,0){\vector(0,1){35}}
\put(330,35){\vector(0,-1){35}}
\put(52,17){$\delta_{1}  $}
\put(333,17){$ \delta_{0} $}
\end{picture}

\begin{picture}(350,15)
\put(20,0){\makebox[60pt]{$K_{1}(C(S^1 \times \{- \infty, + \infty \}))$}}
\put(152,3){\vector(-1,0){35}}
\put(160,0){\makebox[60pt]{$K_{1}(C(M_{ \sharp}))$}}
\put(263,3){\vector(-1,0){35}}
\put(300,0){\makebox[60pt]{$K_{1}(C_{0}( \mathbb{R}))$}}
\put(117,8){\makebox[35pt]{$ \psi_{*} $}}
\put(228,8){\makebox[35pt]{$ \varphi_{*} $}}
\end{picture}

Pelo exemplo 11.3.2 de \cite{[R]} sabemos que $K_{0}(C_{0}( \mathbb{R})) = 0$ e
$K_{1}(C_{0}( \mathbb{R})) = \mathbb{Z}$. Assim como pelo exemplo 11.3.3 de
\cite{[R]} e o fato de $C(S^1 \times \{- \infty, + \infty\})$ ser isomorfo a $C(S^1)
\oplus C(S^1)$ temos para $i= 0, 1$ que $K_{i}(C(S^1 \times \{- \infty, +
\infty \})) \cong \mathbb{Z} \oplus \mathbb{Z}$ onde a classe $[.]_{0}$ das funções
$k,\tilde{k} : S^1 \times \{- \infty, + \infty \} \rightarrow \mathbb{C}$ com
$k(z,+ \infty) = 1$;
$k(z, - \infty) = 0$; $\tilde{k}(z,+ \infty) = 0$;  $\tilde{k}(z, - \infty) = 1$ geram
$K_{0}(C(S^1 \times \{- \infty, + \infty \}))$. E a classe $[.]_{1}$ das
funções $l,\tilde{l} : S^1 \times \{- \infty, + \infty \} \rightarrow
\mathbb{C}$ com $l(z,+ \infty) = z $; $l(z, - \infty) = 1$; $\tilde{l}(z, +
\infty) = 1$; $\tilde{l}(z, - \infty) = z $ geram $K_{1}(C(S^1 \times \{- \infty, + \infty \}))$.

Substituindo na sequência exata cíclica de seis termos as informações  acima temos que $\delta_{1} \equiv 0$. Só nos falta calcular $\delta_{0}$ para obter a K-teoria de $C(M_{ \sharp})$.

Pela Proposição 12.2.2 de \cite{[R]} sabemos que precisamos de funções $g$ e $\tilde{g} \in C(M_{ \sharp})$ reais tais que $\psi(g) = k$ e $\psi(\tilde{g}) = \tilde{k}$.

Consideremos $b: M_{ \sharp} \rightarrow \mathbb{R}$ com $b \equiv 0$ em
$S^{1}_{-}$; $b \equiv 1$ em $S^{1}_{+}$; $b(x) = 0$ para $x \leq -1$ e $b(x) = 1$ para $x \geq 1$.\\
E também $c: M_{ \sharp} \rightarrow \mathbb{R}$ com $c \equiv 1$ em
$S^{1}_{-}$; $c \equiv 0$ em $S^{1}_{+}$; $c(x) = 1$ para $x \leq -1$ e $c(x)
= 0$ para $x \geq 1$. As funções $b$ e $c$ são extensões contínuas das funções
$b$ e $c$ definidas na página 47, mas utilizaremos a mesma notação.

Temos que $\psi(b) = k$ e $\psi(c) = \tilde{k}$, então pela Proposição
12.2.2 de \cite{[R]} temos que existem (e são únicos) elementos $u$ e $\tilde{u} \in \mathcal{U}(\widetilde{C_{0}( \mathbb{R})})$ tais que $\bar{\varphi}(u) = e^{2 \pi i b}$ e $\bar{\varphi}(\tilde{u}) = e^{2 \pi i c}$. Mas como $\widetilde{C_{0}( \mathbb{R})} \subset C(M_{ \sharp})$ temos que $ \varphi = \bar{\varphi}$, logo $u = e^{2 \pi i b}$ e $\tilde{u} = e^{2 \pi i c}$.

Portanto $\delta_{0}([k]_{0}) = - [u]_{1}$ e $\delta_{0}([\tilde{k}]_{0}) = -
[\tilde{u}]_{1}$ onde $[u]_{1}$ e $[\tilde{u}]_{1}$ pertencem a
$K_{1}(C_{0}(\mathbb{R}))$. Sabemos que $K_{1}(C_{0}(\mathbb{R}))$ é isomorfo
a $K_{1}(C(S^1))$ e temos que o número de rotação induz um isomorfismo de $K_{1}(C(S^1))$ para $\mathbb{Z}$.

Calculando o número de rotação de $u$ e $\tilde{u}$ temos que $w(u) = 1$ e
$w(\tilde{u}) = -1$. Logo com respeito aos isomorfismos $K_{0}(C(S^1 \times
\{- \infty, + \infty \})) = \mathbb{Z}[k]_{0} \oplus \mathbb{Z}[\tilde{k}]_{0}$
e $K_{1}(C(S^1 \times \{- \infty, + \infty \})) = \mathbb{Z}[l]_{1} \oplus
\mathbb{Z}[\tilde{l}]_{1}$ discutidos acima, podemos escrever $\delta_{0}(1,0)
= - 1$ e $\delta_{0}(0,1) = 1$, consequentemente $\delta_{0}$ é sobrejetora e
$\delta_{0}(x, y) = -x + y$.

Segue da exatidão da sequência exata cíclica de seis termos que $\mbox{Im}\,
\psi_{*} = \mbox{Ker}\, \delta_{0} \cong \mathbb{Z}$, $\mbox{Ker}\,
\psi_{*} = \mbox{Im}\, \varphi_{*} = 0$, $\mbox{Im}\,
\psi_{*} = \mbox{Ker}\, \delta_{1} \cong \mathbb{Z} \oplus \mathbb{Z}$ e
$\mbox{Ker}\, \psi_{*} = \mbox{Im}\, \varphi_{*} = 0$. 
Portanto $K_{0}(C(M_{ \sharp})) \cong \mathbb{Z}$ e $K_{1}(C(M_{ \sharp})) \cong\mathbb{Z} \oplus \mathbb{Z}$.$\cqd$

Usando os isomorfismos da Proposição acima podemos verificar que $$K_{0}(C(M_{
  \sharp})) = \mathbb{Z}[1]_{0}.$$

Sejam $S^1 = \{e^{ix} \: : x \in \mathbb{R} \}$, $L$ e $\tilde{L} \in
C(M_{\sharp})$ tais que:

$L(e^{ix}, + \infty)= e^{ix}$ em $S^{1}_{+}$; $L(e^{ix}, - \infty)= 1$ em $S^{1}_{-}$; $L(x) =
e^{ix}$, se $x \geq 0$ e $L(x) = 1$, se $x < 0$ e,

$\tilde{L}(e^{ix}, + \infty) = 1$ em $S^{1}_{+}$; $\tilde{L}(e^{ix}, - \infty) = e^{ix}$ em $S^{1}_{-}$; $\tilde{L}(x) = 1$, se $x > 0$ e $\tilde{L}(x) = e^{ix}$, se $x \leq 0$.

Consideremos a restrição de $L$ à reta. Então temos que $b(x)e^{ix} + c(x) -
L(x) \in C_{0}(\mathbb{R})$, onde $b, c \in CS( \mathbb{R})$ são tais que $b$ é crescente com $b(x) = 0$, se $x \leq
-1$ e $b(x) = 1$, se $x \geq 1$ e $c$ é decrescente com $c(x) = 0$, se $x \geq
1$ e $c(x) = 1$, se $x \leq -1$, $b + c \equiv 1$. Logo $L(M) = b(M)e^{iM} + c(M) - a_{0}(M)$,
onde $a_{0} \in C_{0}(\mathbb{R})$. O mesmo
vale para $\tilde{L}$, logo $\tilde{L}(M) = b(M) + c(M)e^{iM} - b_{0}(M)$.

Com o que acabamos de provar, concluímos que $L$, $\tilde{L}$ estão em
$\mathcal{A}^{\sharp}$ (veja definição na seção 1.2). Usando novamente o
isomorfismo da Proposição acima concluímos que $$K_{1}(C(M_{\sharp})) =
\mathbb{Z}[L]_{1} \oplus \mathbb{Z}[\tilde{L}]_{1}.$$
\begin{Prop}
$K_{0}(C(M_{A}))$ é isomorfo a $\mathbb{Z} \oplus \mathbb{Z}$ e $K_{1}(C(M_{A}))$ é isomorfo a $\mathbb{Z} \oplus \mathbb{Z} \oplus \mathbb{Z} \oplus \mathbb{Z}$.
\end{Prop}
\textbf{Demonstração:}

Sabemos pelo Teorema 1.17 que $M_{A} = M_{ \sharp} \times \{ - \infty, + \infty \}$, então consideremos o isomorfismo\\
\begin{displaymath}
\begin{array}{clc}
C(M_{ \sharp} \times \{ - \infty, + \infty \}) & \longrightarrow  & C(M_{ \sharp}) \oplus C(M_{ \sharp})\\
f(z , \pm \infty)                              & \longmapsto      & (f(z, +
\infty), f(z, - \infty)).\\
\end{array}
\end{displaymath}

Logo temos que $K_{i}(C(M_{A})) \cong K_{i}(C(M_{
 \sharp}) \oplus C(M_{ \sharp}))$, $i = 0, 1$.

Pelas Proposições 4.3.4 e 8.2.6 de \cite{[R]} temos que $K_{i}(C(M_{\sharp}) \oplus C(M_{\sharp})) \cong K_{i}(C(M_{\sharp})) \oplus K_{i}(C(M_{\sharp}))$,  logo $K_{i}(C(M_{A})) \cong K_{i}(C(M_{\sharp})) \oplus K_{i}(C(M_{\sharp}))$, $i = 0, 1$.
Consideremos as funções: 

\centerline{$ f:  M_{ \sharp} \times \{ - \infty, + \infty \} \longrightarrow \mathbb{C}$ com $f(z, + \infty) = 1$ e $f(z, - \infty) = 0$,}

\centerline{$\tilde{f}:  M_{ \sharp} \times \{ - \infty, + \infty \} \longrightarrow \mathbb{C}$ com $\tilde{f}(z, + \infty) = 0$ e $\tilde{f}(z, - \infty) = 1$,}

\centerline{$w_{i}: M_{\sharp} \times \{ - \infty, + \infty \} \longrightarrow \mathbb{C}$, $i = 1,..,4$, com}

\centerline{$w_{1}(z, + \infty) = L(z)$ e $ w_{1}(z, - \infty) = 1$,}

\centerline{$w_{2}(z, + \infty) = \tilde{L}(z)$ e $ w_{2}(z, - \infty) = 1$,}

\centerline{$w_{3}(z, + \infty) = 1$ e $ w_{3}(z, - \infty) = L(z)$,}

\centerline{$w_{4}(z, + \infty) = 1$ e $ w_{4}(z, - \infty) = \tilde{L}(z)$, $z
\in M_{\sharp}$.}

Utilizando o isomorfismo acima temos que $$K_{0}(C(M_{A})) = \mathbb{Z}[f]_{0}
\oplus \mathbb{Z}[\tilde{f}]_{0}\hspace{0.5cm} \mbox{e}$$ $$K_{1}(C(M_{A})) =
\mathbb{Z}[w_{1}]_{1} \oplus \mathbb{Z}[w_{2}]_{1} \oplus \mathbb{Z}[w_{3}]_{1}
\oplus \mathbb{Z}[w_{4}]_{1}.$$ $\cqd$

Agora utilizando o $\sigma$-símbolo definido no Teorema 1.17 concluímos que
\begin{equation}
K_{0}(\mathcal{A}/ \mathcal{E}) = \mathbb{Z}[[b(D)]_{\mathcal{E}}]_{0}
\oplus \mathbb{Z}[[c(D)]_{\mathcal{E}}]_{0}\hspace{0.5cm} \mbox{e}
\end{equation}
\begin{equation} 
K_{1}(\mathcal{A}/ \mathcal{E}) = \mathbb{Z}[[A_{1}]_{\mathcal{E}}]_{1}
\oplus \mathbb{Z}[[A_{2}]_{\mathcal{E}}]_{1} \oplus
\mathbb{Z}[[A_{3}]_{\mathcal{E}}]_{1} \oplus
\mathbb{Z}[[A_{4}]_{\mathcal{E}}]_{1} 
\end{equation}
onde $A_{1} =
L(M)b(D) + c(D)$, $A_{2} = \tilde{L}(M)b(D) + c(D)$, $A_{3} = b(D) + L(M)c(D)$
e $A_{4} = b(D) + \tilde{L}(M)c(D)$ com $b, c$ as mesmas funções definidas na
página 54.
 \begin{Prop}
O isomorfismo da Proposição 1.24,  $\gamma:\mathcal{E}/ \mathcal{K}_{\mathbb{R}} \longrightarrow C(M_{SL},
\mathcal{K}_{\mathbb{Z}})$, induz o isomorfismo $K_{i}( \mathcal{E}/
\mathcal{K}_{\mathbb{R}}) \cong \mathbb{Z} \oplus \mathbb{Z}$, $i = 0, 1$.
\end{Prop}
\textbf{Demonstração:}

Usando o isomorfismo $\gamma: \mathcal{E}/
\mathcal{K}_{\mathbb{R}} \longrightarrow C(M_{SL}, \mathcal{K}_{\mathbb{Z}})$
temos que
$K_{i}( \mathcal{E}/ \mathcal{K}_{\mathbb{R}})$ é isomorfo a $K_{i}(C(M_{SL},
\mathcal{K}_{\mathbb{Z}}))$, $i = 0, 1$ onde $M_{SL} = S^1 \times \{-1, 1\}$.

Então vamos calcular a K-teoria de $C(M_{SL}, \mathcal{K}_{\mathbb{Z}})$ para
conhecermos a K-teoria de $\mathcal{E}/ \mathcal{K}_{\mathbb{R}}$.
Começaremos por calcular a K-teoria de $C(S^1, \mathcal{K}_{\mathbb{Z}})$, por
isso consideremos a seguinte sequência exata curta cindida
\begin{center}
\begin{picture}(395,15)
\put(10,0){0}
\put(25,3){\vector(1,0){25}}
\put(58,0){\makebox[60pt]{$S \mathcal{K}_{\mathbb{Z}}$}}
\put(128,3){\vector(1,0){35}}
\put(172,0){\makebox[60pt]{$C(S^1, \mathcal{K}_{\mathbb{Z}})$}}
\put(245,3){\vector(1,0){35}}
\put(285,0){\makebox[60pt]{$\mathcal{K}_{\mathbb{Z}}$}}
\put(353,3){\vector(1,0){25}}
\put(390,0){0}
\put(128,8){\makebox[35pt]{$\mu  $}}
\put(245,8){\makebox[35pt]{$\nu $}}
\end{picture}
\end{center}
onde $\nu (f(z)) = f(1)$ e $S \mathcal{K}_{\mathbb{Z}}$ denota a suspensão
de $\mathcal{K}_{\mathbb{Z}}$, isto é $S \mathcal{K}_{\mathbb{Z}} = \{ f: S^1
\rightarrow \mathcal{K}_{\mathbb{Z}}: f(1) = 0\}$.

E sua sequência exata cíclica de seis termos\\
\begin{picture}(395,15)
\put(68,0){\makebox[60pt]{$K_{0}(S \mathcal{K}_{\mathbb{Z}})$}}
\put(134,3){\vector(1,0){35}}
\put(180,0){\makebox[60pt]{$K_{0}(C(S^1, \mathcal{K}_{\mathbb{Z}}))$}}
\put(251,3){\vector(1,0){35}}
\put(290,0){\makebox[60pt]{$K_{0}(\mathcal{K}_{\mathbb{Z}})$}}
\put(134,8){\makebox[35pt]{$ \mu_{*}  $}}
\put(251,8){\makebox[35pt]{$ \nu_{*} $}}
\end{picture}

\begin{picture}(350,45)
\put(78,0){\vector(0,1){35}}
\put(300,35){\vector(0,-1){35}}
\put(81,17){$\delta_{1}  $}
\put(303,17){$ \delta_{0} $}
\end{picture}

\begin{picture}(350,15)
\put(48,0){\makebox[60pt]{$K_{1}(\mathcal{K}_{\mathbb{Z}})$}}
\put(149,3){\vector(-1,0){35}}
\put(160,0){\makebox[60pt]{$K_{1}(C(S^1, \mathcal{K}_{\mathbb{Z}}))$}}
\put(266,3){\vector(-1,0){35}}
\put(270,0){\makebox[60pt]{$K_{1}(S \mathcal{K}_{\mathbb{Z}})$}}
\put(117,8){\makebox[35pt]{$ \nu_{*} $}}
\put(228,8){\makebox[35pt]{$ \mu_{*} $}}
\end{picture} 

Como a sequência é cindida temos que $\delta_{0} = \delta_{1} = 0$. Segue da
exatidão da sequência exata cíclica de seis termos que $\nu_{*}: K_{0}(C(S^1,
\mathcal{K}_{\mathbb{Z}})) \rightarrow K_{0}(\mathcal{K}_{\mathbb{Z}}) =
\mathbb{Z}[E]_{0}$, onde $E$ é qualquer projeção de posto 1 em
$\mathcal{L}(L^2(\mathbb{Z}))$ é isomorfismo e que $\mu_{*}: K_{1}(S
\mathcal{K}_{\mathbb{Z}}) \rightarrow  K_{1}(C(S^1,
\mathcal{K}_{\mathbb{Z}}))$ é isomorfismo.

Usando o isomorfismo da periodicidade de Bott temos que $ K_{1}(S
\mathcal{K}_{\mathbb{Z}}) = \mathbb{Z}[1 + (z -1)E]_{1}$
 logo $K_{1}(C(S^1,\mathcal{K}_{\mathbb{Z}})) = \mathbb{Z}[1 + (z -1)E]_{1}$,
 $z$ denotando a função identidade em $S^1$.

Considerando o isomorfismo 
\begin{displaymath}
\begin{array}{clc}
C(S^1 \times \{ - 1, 1 \}, \mathcal{K}_{\mathbb{Z}} ) & \longrightarrow  &
C(S^1, \mathcal{K}_{\mathbb{Z}} ) \oplus C(S^1, \mathcal{K}_{\mathbb{Z}})\\
f(z , \pm 1)                              & \longmapsto      & (f(z, 1), f(z, - 1)),\\
\end{array}
\end{displaymath}
temos que 
$K_{0}(C(M_{SL}, \mathcal{K}_{\mathbb{Z}})) \cong \mathbb{Z} \oplus
\mathbb{Z}$ e $$K_{1}(C(M_{SL}, \mathcal{K}_{\mathbb{Z}})) =
\mathbb{Z}[u_{1}]_{1} \oplus \mathbb{Z}[u_{2}]_{1},$$ onde $$u_{1}(z, 1) = 1 + (z -1)E \hspace{0.3cm} \mbox{e}
\hspace{0.3cm} u_{2}(z, -1) = 1$$
\centerline{e}
$$u_{2}(z, 1) = 1 \hspace{0.3cm} \mbox{e}
\hspace{0.3cm} u_{2}(z, -1) = 1 + (z -1)E.$$

Pelo isomorfismo $\gamma$ temos que
$K_{0}(\mathcal{E}/ \mathcal{K}_{\mathbb{R}}) \cong \mathbb{Z} \oplus
  \mathbb{Z}$ e $K_{1}(\mathcal{E}/ \mathcal{K}_{\mathbb{R}}) \cong
    \mathbb{Z} \oplus \mathbb{Z}$. $\cqd$
\section{A K-teoria de $\mathcal{A}/ \mathcal{K}_{\mathbb{R}}$ e $\mathcal{A}$}

\hspace{8mm}Vamos calcular agora a K-teoria de $\mathcal{A}/
\mathcal{K}_{\mathbb{R}}$ e então estaremos prontos para calcular a K-teoria
de $\mathcal{A}$. Enunciaremos e demonstraremos antes um Lema que será
utilizado no cálculo da aplicação do índice da sequência exata cíclica de seis termos de
$\mathcal{A}/ \mathcal{K}_{\mathbb{R}}$. Nesta seção e nas que seguem
estaremos utilizando as definições de $K_{0}$ e $K_{1}$ de \cite{[Bk]}, isto é
$K_{0}$ é um grupo de diferenças formais de classes de idempotentes (matrizes)  de
tamanho arbitrário e $K_{1}$ é
o grupo de classes de homotopias de matrizes inversíveis de tamanho arbitrário. A demonstração de que as
definições apresentadas em \cite{[Bk]} e \cite{[R]} (que usa diferenças formais
de classes de projeções ortogonais (matrizes) para $K_{0}$ e de unitários para $K_{1}$)
são equivalentes encontra-se em \cite{[Bk]}.
\begin{Lema}
Sejam $A$ uma $C^*$- álgebra e $J$ um ideal de $A$. Consideremos a sequência
exata curta :
\begin{center}
\begin{picture}(280,15)
\put(15,0){0}
\put(30,3){\vector(1,0){25}}
\put(43,0){\makebox[60pt]{$J$}}
\put(90,3){\vector(1,0){35}}
\put(110,0){\makebox[60pt]{$A$}}
\put(152,3){\vector(1,0){35}}
\put(181,0){\makebox[60pt]{$A/J$}}
\put(230,3){\vector(1,0){25}}
\put(262,0){0}
\put(90,8){\makebox[35pt]{$i $}}
\put(152,8){\makebox[35pt]{$\pi$}}
\end{picture}
\end{center}
onde $i$ é a inclusão e $\pi$ é a projeção canônica.\\
 Se $u \in M_{n}(A/J)$ é inversível, $\pi(a) = u$ e $\pi(b) = u^{-1}$ então:
\begin{displaymath}
 \delta_{1}([u]_{1}) = \left [ \left(
  \begin{array}{cc}
   2ab - (ab)^2 & a(2 - ba)(1 - ba)\\
   (1 - ba)b    & (1 - ba)^2\\
  \end{array}
 \right ) \right ]_{0} - \left [ \left (
  \begin{array}{ll}
   1 & 0\\
   0 & 0\\
  \end{array} 
 \right ) \right ]_{0}.
\end{displaymath}
\end{Lema}

\textbf{Demonstração:}

Seja $w \in M_{2n}(A)$ definido por
$$ w = \left(   
\begin{array}{cc}
   2a - aba & ab - 1\\
   1 - ba   & b\\
\end{array}
 \right )$$
É fácil verificar que $w$ é inversível e
$$ w^{-1} = \left(   
\begin{array}{cc}
   b     & 1 - ba\\
   ab -1 & 2a - aba  \\
\end{array}
 \right ).$$
Temos também que 
$$\pi \left( 
\begin{array}{cc}
   2a - aba & ab - 1\\
   1 - ba   & b\\
\end{array}
 \right ) = \left( 
\begin{array}{cc}
u  & 0\\
0  & u^{-1}\\
\end{array}
\right).$$
Pela definição 8.3.1 de \cite{[Bk]} segue que
\begin{displaymath}
 \delta_{1}([u]_{1}) = \left [ \left(
  \begin{array}{cc}
   2a - aba & ab - 1\\
   1 - ba   & b\\
  \end{array}
 \right )  \left(
\begin{array}{ll}
1 & 0\\
0 & 0\\
\end{array}
\right ) \left ( 
\begin{array}{cc}
b     & 1 - ba\\
ab -1 & 2a - aba\\
\end{array}
\right ) \right ]_{0} - \left [ \left (
  \begin{array}{ll}
   1 & 0\\
   0 & 0\\
  \end{array} 
 \right ) \right ]_{0}
\end{displaymath}
\begin{displaymath}
= \left [ \left(
  \begin{array}{cc}
   2ab - (ab)^2 & a(2 - ba)(1 - ba)\\
   (1 - ba)b    & (1 - ba)^2\\
  \end{array}
 \right ) \right ]_{0} - \left [ \left (
  \begin{array}{ll}
   1 & 0\\
   0 & 0\\
  \end{array} 
 \right ) \right ]_{0}.
\end{displaymath}$\cqd$
\begin{Teo}
Seja $\psi: \mathcal{A}/ \mathcal{K}_{\mathbb{R}} \rightarrow  \mathcal{A}/
\mathcal{E}$ onde $[a]_{\mathcal{K}_{\mathbb{R}}} \mapsto [a]_{\mathcal{E}}$,
$a \in \mathcal{A}$. Temos que $\psi_{*}:K_{0}(\mathcal{A}/
\mathcal{K}_{\mathbb{R}}) \rightarrow K_{0}(\mathcal{A}/
\mathcal{E})$ é injetora e a imagem de $\psi_{*}$ é isomorfa a
$\mathbb{Z}$. Assim,  $K_{0}(\mathcal{A}/ \mathcal{K}_{\mathbb{R}}) \cong
\mathbb{Z}$ e temos também que
$K_{1}(\mathcal{A}/ \mathcal{K}_{\mathbb{R}}) \cong \mathbb{Z} \oplus \mathbb{Z} \oplus \mathbb{Z}$.
\end{Teo}
\textbf{Demonstração:}

Consideremos a sequência exata curta:
\begin{center}
\begin{picture}(350,15)
\put(10,0){0}
\put(20,3){\vector(1,0){25}}
\put(48,0){\makebox[60pt]{$\mathcal{E}/ \mathcal{K}_{\mathbb{R}}$}}
\put(111,3){\vector(1,0){35}}
\put(149,0){\makebox[60pt]{$\mathcal{A}/ \mathcal{K}_{\mathbb{R}}$}}
\put(212,3){\vector(1,0){35}}
\put(250,0){\makebox[60pt]{$ \mathcal{A}/ \mathcal{E}$}}
\put(313,3){\vector(1,0){25}}
\put(341,0){0}
\put(111,8){\makebox[35pt]{$ \varphi$}}
\put(212,8){\makebox[35pt]{$\psi $}}
\end{picture}
\end{center}
onde $\varphi$ é a inclusão e $\psi: [a]_{\mathcal{K}_{\mathbb{R}}} \longmapsto [a]_{\mathcal{E}}$, $a \in \mathcal{A}$.

E também sua sequência exata cíclica de seis termos para K-teoria\\
\begin{picture}(395,15)
\put(68,0){\makebox[60pt]{$K_{0}(\mathcal{E}/ \mathcal{K}_{\mathbb{R}})$}}
\put(137,3){\vector(1,0){35}}
\put(180,0){\makebox[60pt]{$K_{0}(\mathcal{A}/ \mathcal{K}_{\mathbb{R}})$}}
\put(248,3){\vector(1,0){35}}
\put(290,0){\makebox[60pt]{$K_{0}(\mathcal{A}/ \mathcal{E})$}}
\put(137,8){\makebox[35pt]{$ \varphi_{*}  $}}
\put(248,8){\makebox[35pt]{$ \psi_{*} $}}
\end{picture}

\begin{picture}(350,45)
\put(78,0){\vector(0,1){35}}
\put(300,35){\vector(0,-1){35}}
\put(81,17){$\delta_{1}  $}
\put(303,17){$ \delta_{0} $}
\end{picture}

\begin{picture}(350,15)
\put(48,0){\makebox[60pt]{$K_{1}(\mathcal{A}/ \mathcal{E})$}}
\put(152,3){\vector(-1,0){35}}
\put(160,0){\makebox[60pt]{$K_{1}(\mathcal{A}/ \mathcal{K}_{\mathbb{R}})$}}
\put(263,3){\vector(-1,0){35}}
\put(270,0){\makebox[60pt]{$K_{1}(\mathcal{E}/ \mathcal{K}_{\mathbb{R}})$}}
\put(117,8){\makebox[35pt]{$ \psi_{*} $}}
\put(228,8){\makebox[35pt]{$ \varphi_{*} $}}
\end{picture}

Substituindo as K-teorias que conhecemos temos\\
\begin{picture}(395,15)
\put(68,0){\makebox[60pt]{$\mathbb{Z} \oplus \mathbb{Z}$}}
\put(137,3){\vector(1,0){35}}
\put(180,0){\makebox[60pt]{$K_{0}(\mathcal{A}/ \mathcal{K}_{\mathbb{R}})$}}
\put(248,3){\vector(1,0){35}}
\put(280,0){\makebox[60pt]{$\mathbb{Z} \oplus \mathbb{Z}$}}
\put(137,8){\makebox[35pt]{$ \varphi_{*}  $}}
\put(248,8){\makebox[35pt]{$ \psi_{*} $}}
\end{picture}

\begin{picture}(350,45)
\put(78,0){\vector(0,1){35}}
\put(290,35){\vector(0,-1){35}}
\put(81,17){$\delta_{1}  $}
\put(293,17){$ \delta_{0} $}
\end{picture}

\begin{picture}(350,15)
\put(46,0){\makebox[60pt]{$\mathbb{Z} \oplus \mathbb{Z} \oplus \mathbb{Z} \oplus \mathbb{Z}$}}
\put(152,3){\vector(-1,0){35}}
\put(160,0){\makebox[60pt]{$K_{1}(\mathcal{A}/ \mathcal{K}_{\mathbb{R}})$}}
\put(263,3){\vector(-1,0){35}}
\put(260,0){\makebox[60pt]{$\mathbb{Z} \oplus \mathbb{Z} $}}
\put(117,8){\makebox[35pt]{$ \psi_{*} $}}
\put(228,8){\makebox[35pt]{$ \varphi_{*} $}}
\end{picture} 
 
Precisamos calcular $\delta_{0}$ e $\delta_{1}$ para conhecermos então a
K-teoria de $\mathcal{A}/ \mathcal{K}_{\mathbb{R}}$.

Vamos primeiro calcular $\delta_{0}([[b(D)]_{\mathcal{E}}]_{0})$, onde
$[[b(D)]_{\mathcal{E}}]_{0}$ foi descrito na página 55 equação (3.1). Pela Proposição 12.2.2 de \cite{[R]} sabemos
que precisamos de um elemento $x \in
\mathcal{A}$ tal que $[x]_{\mathcal{K}_{\mathbb{R}}}$ seja auto-adjunto e
$\psi([x]_{\mathcal{K}_{\mathbb{R}}}) =[b(D)]_{\mathcal{E}}$. 

Como $b(D)$ já é auto-adjunto em $\mathcal{A}$ basta então tomarmos $x =
b(D)$. Pela Proposição 12.2.2 de \cite{[R]} temos que
existe um único elemento $u \in \mathcal{U}( \widetilde{\mathcal{E}/
\mathcal{K}_{\mathbb{R}}})$ tal que $\bar{\varphi}(u) = e^{2 \pi i
  [b(D)]_{\mathcal{K}_{\mathbb{R}}}}$. Como $\widetilde{\mathcal{E}/
  \mathcal{K}_{\mathbb{R}}} \subset \mathcal{A}/ \mathcal{K}_{\mathbb{R}}$ temos $ \varphi =
\bar{\varphi}$, logo pela definição de $\varphi$ temos que $u = e^{2 \pi i
  [b(D)]_{\mathcal{K}_{\mathbb{R}}}}$. Portanto o gerador
$[[b(D)]_{\mathcal{E}}]_{0}$ de $K_{0}( \mathcal{A} / \mathcal{E})$ é levado
em $-[e^{2 \pi i [b(D)]_{\mathcal{K}_{\mathbb{R}}}}]_{1} \in K_{1}(
\mathcal{E}/ \mathcal{K}_{\mathbb{R}})$ pela aplicação
$\delta_{0}$. Utilizando o isomorfismo $\gamma$ da Proposição 1.24 e
homotopias verificamos em seguida que $[e^{2 \pi i
  [b(D)]_{\mathcal{K}_{\mathbb{R}}}}]_{1}$ corresponde ao elemento $(-1, -1)$ de $K_{1}(C(M_{SL}, \mathcal{K}_{ \mathbb{Z}}))$.

Usando o isomorfismo da Proposição 1.24 temos que $[e^{2 \pi i
  [b(D)]_{\mathcal{K}_{\mathbb{R}}}}]_{1}$ é levado em $[h]_{1}$ onde $h$ é a função:
\begin{displaymath}
\begin{array}{cccc}
h: & M_{SL}                       & \longrightarrow  & \mathcal{K}_{ \mathbb{Z}}\\
   & (e^{2 \pi i \varphi}, \pm 1) & \longmapsto      & Y_{ \varphi}(e^{2
  \pi i b(M - \varphi)})Y_{- \varphi}.\\
\end{array}
\end{displaymath}

Com respeito ao isomorfismo da Proposição 3.3, sabemos que o elemento que
corresponde a $(1, 1) \in \mathbb{Z} \oplus \mathbb{Z}$ em $K_{1}(C(M_{SL},
\mathcal{K}_{ \mathbb{Z}}))$ é dado pela classe em $K_{1}$ de 
\begin{displaymath}
\begin{array}{cccc}
l: & M_{SL}                       & \longrightarrow & \mathcal{K}_{ \mathbb{Z}}\\
   & (z, \pm 1) & \longmapsto     & 1 + (z - 1) E,\\
\end{array}
\end{displaymath}
onde $E((a_{n})_{n}) = a_{0}e_{0}$.
Seja $a \in CS(\mathbb{R})$ então,\\[0.2cm]
\centerline{$Y_{ \varphi + 1}a(M - ( \varphi + 1))Y_{-( \varphi + 1)} = Y_{
    \varphi}Y_{1}a(M - ( \varphi + 1))Y_{- 1}Y_{-\varphi} = Y_{ \varphi}a(
M-  \varphi )Y_{-\varphi}$}\\[.2cm]
logo $Y_{ \varphi}(e^{2 \pi i b(M - \varphi)})Y_{- \varphi}$ é periódica de
período 1. 
Consideremos então $Y_{ \varphi}(e^{2 \pi i b(M - \varphi)})Y_{- \varphi}$
para $\varphi \in [ - \frac{1}{2}, \frac{1}{2} ]$.

Sem perda de generalidade podemos supor que a função $b$ é constante em $( -
\infty, -1/5]$ e em $[1/5, + \infty)$. Se
$1/5 \leq | \varphi| \leq 1/2$ então $| j - \varphi| \geq \frac{1}{5}$ para
todo $j$. Se $|\varphi| < 1/5$ então $| j - \varphi| \geq \frac{1}{5}$ para
todo $j$ não nulo.

Logo para todo $\varphi \in [-1/2, 1/2]$ temos que
\begin{displaymath}
 e^{2 \pi i b(j - \varphi)}= \left\{ 
 \begin{array}{lll}
   e^{2 \pi i b(- \varphi)} & se & j = 0,\\
   1 & se & j \neq 0.\\
 \end{array}
 \right.
\end{displaymath}
Portanto $e^{2 \pi i b(M - \varphi)} = Id + (e^{2 \pi i b(- \varphi)} - Id)E$
para todo $|\varphi| \leq 1/2$. Seja $g(\varphi) = e^{2 \pi i b( - \varphi)}
-Id$ temos que $g \in C^{\infty}_{c}([ - \frac{1}{5}, \frac{1}{5} ])$ e $e^{2
  \pi i b(M - \varphi)} = g(\varphi)E + Id$ logo $Y_{ \varphi}(e^{2 \pi i b(M
  - \varphi)})Y_{- \varphi} = Y_{ \varphi}(g( \varphi)E + Id)Y_{ -
  \varphi}$, para todo $|\varphi| \leq 1/2$.

Seja $Z_{t} = Y_{h(t)}$ onde $h: [-1/2, 1/2] \rightarrow \mathbb{R}$ com $h
\equiv 0$ em $[-1/4, 1/4]$. Logo $Z_{t} = Id$ se $t \in [ - \frac{1}{4}, \frac{1}{4} ]$.

Seja $v: [0, 1] \longrightarrow \mathcal{U}(C([-1/2, 1/2]), \mathcal{K}_{\mathbb{Z}}))$ onde
$x \longmapsto v_{x}$ definida por $v_{x}(t) = Y_{(1-x)t + xh(t)}(g(t)E +
1)Y_{-[ (1-x)t + xh(t)]}$. Se $|t| \geq 1/5$ temos $g(t)E + Id = Id$. Logo
para todo $x$, $v_{x}$ pode ser encarado como um elemento de $\mathcal{U}(C(S^1,
\mathcal{K}_{\mathbb{Z}}))$ e  $v_{0}(t) = Y_{t}(g(t)E + 1)Y_{-t}$ e $v_{1}(t)
= Y_{h(t)}(g(t)E + 1)Y_{-h(t)} = Z_{t}(g(t)E + 1)Z_{-t}$. Portanto
$Y_{t}(g(t)E + 1)Y_{-t}$ e $Z_{t}(g(t)E + 1)Z_{-t}$ são elementos homotópicos
de $C(S^1,\mathcal{K}_{\mathbb{Z}})$.

Como $Z_{t} = Id$ se $t \in [ - \frac{1}{4}, \frac{1}{4} ]$ e $ g$ tem suporte
contido em $[ - \frac{1}{5}, \frac{1}{5} ]$ temos que $Z_{t}(g(t)E + 1)Z_{-t}
= g(t)E + 1$ para todo $t \in [-1/2, 1/2]$ . Então pela Proposição 8.1.4 de \cite{[R]} e pelos resultados acima temos   

$$[Y_{t}(g(t)E + 1)Y_{-t}]_{1} = [Z_{t}(g(t)E + 1)Z_{-t}]_{1} = [g(t)E + 1]_{1}, $$
onde $F(t)$ denota a função $(e^{2 \pi i t}, \pm 1) \ni M_{SL} \mapsto
F(t) \in \mathcal{K}_{\mathbb{Z}}$.

Usando o que foi visto na página 47 sobre número de rotação temos que  as funções $\bar{z} = (z
\longmapsto \bar{z})$ e $ (e^{2 \pi i \varphi} \longmapsto e^{2 \pi i
  b(-\varphi)}) = \tilde{g}(\varphi)$ são homotópicas em
$\mathcal{U}(C(S^1))$. Seja então $l: [0, 1] \longrightarrow
\mathcal{U}(C(S^1))$ tal que $l(0) = \bar{z}$ e $l(1) = \tilde{g}(\varphi)$ e
seja $v: [0, 1] \longrightarrow \mathcal{U}(C(S^1, \mathcal{K}_{\mathbb{Z}}))$
dada por  
$  x \longmapsto (l(x) - 1 )E + 1$. Temos que $v(0) = ( \tilde{g}(\varphi) -
1)E + 1 = g(\varphi)E + 1$ e $v(1) = ( \bar{z} - 1)E + 1$. E sabemos que $(
\bar{z} - 1)E + 1$ é o $-1$ em $K_{1}(C(S^1, \mathcal{K}_{\mathbb{Z}}))$.(veja
demonstração da Proposição 3.3)

Portanto pela Proposição 8.1.4 de \cite{[R]} temos $[g(\varphi)E + 1]_{1} = [(\bar{z} - 1)E + 1)]_{1}$, logo $\delta_{0}((1,0)) = (1,1)$.

Com uma demonstração análoga podemos mostrar que $\delta_{0}((0, 1)) =
(-1,-1)$. Então temos que $\delta_{0}((x,y)) = (x-y)(1,1)$. Com isto segue que
$\mbox{Ker}\, \delta_{0} = \{(x,y) \in \mathbb{Z} \oplus \mathbb{Z} : x = y \} \cong  \mathbb{Z}$ e $\mbox{Im}\, \delta_{0} = \{(a, a): a \in \mathbb{Z} \} \cong \mathbb{Z}$.

Com isso encerramos o cálculo de $\delta_{0}$ e passaremos agora a calcular a aplicação $ \delta_{1}$.

Já sabemos pela Proposição 3.2 que os quatro geradores de $K_{1}(
\mathcal{A}/ \mathcal{E})$ são $[[A_{1}]_{\mathcal{E}}]_{1}$, $[[A_{2}]_{\mathcal{E}}]_{1}$, $[[A_{3}]_{\mathcal{E}}]_{1}$, $[[A_{4}]_{\mathcal{E}}]_{1}$, onde
$A_{1} = L(M)b(D) + c(D)$, $A_{2} = \tilde{L}(M)b(D) + c(D)$, $A_{3} = b(D) +
L(M)c(D)$,  $A_{4} = b(D) + \tilde{L}(M)c(D)$. 

É obvio que $a = [A_{1}]_{\mathcal{K}_{\mathbb{R}}}$ é uma pré-imagem de $u = [A_{1}]_{\mathcal{E}}$. Como $[A_{1}]_{\mathcal{E}} =
u$ é unitário temos que $b = a^*$ é uma pré-imagem de $u^{-1}$. Pelo Lema 3.4
temos que 
\begin{displaymath}
 \delta_{1}([[A_{1}]_{\mathcal{E}}]_{1}) = \left [ \left(
  \begin{array}{cc}
   2aa^* - (aa^*)^2 & a(2 - a^*a)(1 - a^*a)\\
   (1 - a^*a)a^* & (1 - a^*a)^2\\
  \end{array}
 \right ) \right ]_{0} - \left [ \left (
  \begin{array}{ll}
   1 & 0\\
   0 & 0\\
  \end{array} 
 \right ) \right ]_{0}
\end{displaymath}
onde $\delta_{1}([[A_{1}]_{\mathcal{E}}]_{1}) \in K_{0}(
\mathcal{E}/ \mathcal{K}_{\mathbb{R}}) \;{\buildrel \gamma_{*} \over
  \cong}\; K_{0}(C(M_{SL}, \mathcal{K}_{\mathbb{Z}})) \;{\buildrel \nu_{*} \over
  \cong}\; K_{0}(\mathcal{K}_{\mathbb{Z}}) \oplus
K_{0}(\mathcal{K}_{\mathbb{Z}})$ (veja Proposição 3.3). Seja $ \eta = \nu_{*} \circ \gamma_{*}$ então
considerando $\alpha = \gamma_{A_{1}}(1,1)$, $\alpha^{'}
=\gamma_{A_{1}}(1,-1)$, $\beta = \gamma_{A^{*}_{1}}(1,1)$ e $\beta^{'} =
\gamma_{A_{1}^{*}}(1, -1)$ temos   
\begin{displaymath}
 \eta \circ \delta_{1}([[A_{1}]_{\mathcal{E}}]_{1}) =  \left ( \left [ \left(
  \begin{array}{cc}
   2 \alpha \beta - (\alpha \beta)^2 & \alpha(2 - \beta \alpha)(1 - \beta \alpha)\\
   (1 - \beta \alpha) \beta & (1 - \beta \alpha)^2\\
  \end{array}
 \right ) \right ]_{0} - \left [ \left (
  \begin{array}{ll}
   1 & 0\\
   0 & 0\\
  \end{array} 
 \right ) \right ]_{0},
\right.
\end{displaymath}
\begin{displaymath}
  \left. \left [ \left(
  \begin{array}{cc}
   2 \alpha^{'} \beta^{'} - (\alpha^{'} \beta^{'})^2 & \alpha^{'}(2 - \beta^{'} \alpha^{'})(1 - \beta^{'} \alpha^{'})\\
   (1 - \beta^{'} \alpha^{'}) \beta^{'} & (1 - \beta^{'} \alpha^{'})^2\\
  \end{array}
 \right ) \right ]_{0} - \left [ \left (
  \begin{array}{ll}
   1 & 0\\
   0 & 0\\
  \end{array} 
 \right ) \right ]_{0} \right )
\end{displaymath}
Vamos observar que $\gamma_{A_{1}}(1, \pm 1) \notin \mathcal{K}_{\mathbb{Z}}$,
já que $A_{1} \in \mathcal{A}$ e $A_{1} \notin \mathcal{E}$. Observemos também
que $A^{*}_{1} A_{1} - Id
\in \mathcal{E}$ então $\gamma_{A^{*}_{1} A_{1}}(1, \pm 1) - \gamma_{Id}(1,
\pm 1) \in
\mathcal{K}_{\mathbb{Z}}$. Portanto $[\gamma_{A_{1}}(1,
1)]_{\mathcal{K}_{\mathbb{Z}}}$ e $[\gamma_{A_{1}}(1,
-1)]_{\mathcal{K}_{\mathbb{Z}}}$ são inversíveis em
$\mathcal{L}(L^2( \mathbb{Z}))/ \mathcal{K}_{\mathbb{Z}}$.

Temos então pelo Lema 3.4 que  
\begin{displaymath}
 \bar{\delta}_{1}( [\gamma_{A_{1}}(1,1)]_{1}) =   \left [ \left(
  \begin{array}{cc}
   2 \alpha \beta - (\alpha \beta)^2 & \alpha(2 - \beta \alpha)(1 - \beta \alpha)\\
   (1 - \beta \alpha) \beta & (1 - \beta \alpha)^2\\
  \end{array}
 \right ) \right ]_{0} - \left [ \left (
  \begin{array}{ll}
   1 & 0\\
   0 & 0\\
  \end{array} 
 \right ) \right ]_{0} 
\end{displaymath}
e
\begin{displaymath}
 \bar{\delta}_{1}( [\gamma_{A_{1}}(1,-1)]_{1}) =   \left [ \left(
  \begin{array}{cc}
   2 \alpha^{'} \beta^{'} - (\alpha^{'} \beta^{'})^2 & \alpha^{'}(2 - \beta^{'} \alpha^{'})(1 - \beta^{'} \alpha^{'})\\
   (1 - \beta^{'} \alpha^{'}) \beta^{'} & (1 - \beta^{'} \alpha^{'})^2\\
  \end{array}
 \right ) \right ]_{0} - \left [ \left (
  \begin{array}{ll}
   1 & 0\\
   0 & 0\\
  \end{array} 
 \right ) \right ]_{0} 
\end{displaymath}
onde $\bar{\delta}_{1}$ é o índice da sequência de K-teoria associada a
sequência exata curta abaixo
\begin{center}
\begin{picture}(395,15)
\put(10,0){0}
\put(30,3){\vector(1,0){25}}
\put(48,0){\makebox[60pt]{$\mathcal{K}_{\mathbb{Z}}$}}
\put(101,3){\vector(1,0){35}}
\put(149,0){\makebox[60pt]{$\mathcal{L}(L^2( \mathbb{Z}))$}}
\put(212,3){\vector(1,0){35}}
\put(270,0){\makebox[60pt]{$\mathcal{L}(L^2( \mathbb{Z}))/ \mathcal{K}_{\mathbb{Z}} $}}
\put(353,3){\vector(1,0){25}}
\put(388,0){0}
\put(101,8){\makebox[35pt]{$ i$}}
\put(212,8){\makebox[35pt]{$\pi $}}
\end{picture}.
\end{center}

Mas o índice para esta sequência é o índice de Fredholm por 9.4.2 de \cite{[R]}
. Por 19.1.14 de \cite{[Hormander]} temos que\\ 
\centerline{$ind \gamma_{A_{1}}(1, \pm 1) = Tr((Id - \gamma_{A^{*}_{1}}(1, \pm
  1)  \gamma_{A_{1}}(1, \pm 1) )^N) - Tr((Id - \gamma_{A_{1}}(1, \pm 1)
  \gamma_{A^{*}_{1}}(1, \pm 1) )^N)$}\\
 para algum $N > 0$ natural para o qual se tenha 

\centerline{$(Id - \gamma_{A^{*}_{1}}(1, \pm 1)  \gamma_{A_{1}}(1, \pm 1) )^N$ e $(Id - \gamma_{A_{1}}(1, \pm 1)  \gamma_{A^{*}_{1}}(1, \pm 1) )^N$ trace-class.}

O mesmo raciocínio vale para $A_{2}$, $A_{3}$, $A_{4}$.

Utilizando a Proposição 1.25 temos que \\[0.4cm]
$\gamma_{A_{1}}(1,-1) = Id$ e $\gamma_{A^{*}_{1}}(1,-1) = Id$,\\
$\gamma_{A_{1}}(1,1) = b(M-1)Y_{-1} + Y_{1}c(M-1)Y_{-1}$ e $\gamma_{A^{*}_{1}}(1,1) = Y_{1}b(M-1) + Y_{1}c(M-1)Y_{-1}$,\\
$\gamma_{A_{2}}(1,-1) = b(M-1)Y_{-1} + Y_{1}c(M-1)Y_{-1}$ e $\gamma_{A^{*}_{2}}(1,-1) = Y_{1}b(M-1) + Y_{1}c(M-1)Y_{-1}$,\\
$\gamma_{A_{2}}(1,1) = Id$ e $\gamma_{A^{*}_{2}}(1,1) = Id$,\\
$\gamma_{A_{3}}(1,-1) = Id$ e $\gamma_{A^{*}_{3}}(1,-1) = Id$,\\
$\gamma_{A_{3}}(1,1) = Y_{1}b(M-1)Y_{-1} + c(M-1)Y_{-1}$ e $\gamma_{A^{*}_{3}}(1,1) = Y_{1}b(M-1)Y_{-1} + Y_{1}c(M-1)$,\\
$\gamma_{A_{4}}(1,-1) = Y_{1}b(M-1)Y_{-1} + c(M-1)Y_{-1}$ e $\gamma_{A^{*}_{4}}(1,-1) = Y_{1}b(M-1)Y_{-1} + Y_{1}c(M-1)$,\\
$\gamma_{A_{4}}(1,1) = Id$ e $\gamma_{A^{*}_{4}}(1,1) = Id$.

Construindo as matrizes desses operadores podemos verificar que para $N = 1$ temos que
$(Id - \gamma_{A^{*}_{i}}(1, \pm 1) \gamma_{A_{i}}(1, \pm 1))$ e $(Id -
\gamma_{A_{i}}(1, \pm 1) \gamma_{A^{*}_{i}}(1, \pm 1))$, para $i = 1,..,4$, são
de posto finito.

Apresentamos a seguir as matrizes dos operadores $(Id - \gamma_{A^{*}_{1}}(1,
 1) \gamma_{A_{1}}(1, 1))$ e \\$(Id - \gamma_{A_{1}}(1, 1) \gamma_{A^{*}_{1}}(1,
 1))$.
\begin{displaymath}
\left[ \begin{array}{ccccccccccc}
         0 & 0 & \ldots & \ldots & \ldots & 0 & 0 & 0 & 0 & 0 & \ldots \\
         0 & \ddots & \vdots & \vdots & \vdots & \vdots & \vdots & \vdots & \vdots
         & \vdots & \vdots\\
    \vdots & \vdots & \ddots & \vdots & \vdots & \vdots & \vdots & \vdots &
         \vdots & \vdots & \vdots\\
    \vdots & \vdots & \vdots & \ddots & \vdots & \vdots & \vdots & \vdots &
         \vdots & \vdots & \vdots \\
    \vdots & \vdots & \vdots & \vdots & \ddots & \vdots & \vdots & \vdots &
         \vdots & \vdots & \vdots\\
    \vdots & \vdots & \vdots & \vdots & \vdots & 0 & \vdots & \vdots & \vdots
         & \vdots & \vdots\\
    \vdots & \vdots & \vdots & \vdots & \vdots & \vdots & 0 & -1/2 & 0 &
         \vdots & \vdots\\
    \vdots & \vdots & \vdots & \vdots & \vdots & \vdots & -1/2 & 1/2 & -1/2 &
         0 & \vdots\\
    \vdots & \vdots & \vdots & \vdots & \vdots & \vdots & 0 & -1/2 & 0 &
         \vdots & \vdots\\
    \vdots & \vdots & \vdots & \vdots & \vdots & \vdots & \vdots & 0 & 0 &
         \ddots & \vdots \\
    \vdots & \vdots & \vdots & \vdots & \vdots & \vdots & \vdots & \vdots &
         \vdots & \vdots & \ddots\\
   \end{array}
\right]
\end{displaymath}
 \begin{displaymath}
\left[ \begin{array}{ccccccccccc}
         0 & 0 & \ldots & \ldots & \ldots & 0 & 0 & 0 & 0 & 0 & \ldots \\
         0 & \ddots & \vdots & \vdots & \vdots & \vdots & \vdots & \vdots & \vdots
         & \vdots & \vdots \\
    \vdots & \vdots & \ddots & \vdots & \vdots & \vdots & \vdots & \vdots &
         \vdots & \vdots & \vdots\\
    \vdots & \vdots & \vdots & \ddots & \vdots & \vdots & \vdots & \vdots &
         \vdots & \vdots & \vdots \\
    \vdots & \vdots & \vdots & \vdots & \ddots & \vdots & \vdots & \vdots &
         \vdots & \vdots & \vdots\\
    \vdots & \vdots & \vdots & \vdots & \vdots & 0 & \vdots & \vdots & \vdots
         & \vdots & \vdots\\
    \vdots & \vdots & \vdots & \vdots & \vdots & \vdots & -1/4 & -1/4 & 0 &
         \vdots & \vdots\\
    \vdots & \vdots & \vdots & \vdots & \vdots & \vdots & -1/4 & -1/4 & 0 &
         0 & \vdots\\
    \vdots & \vdots & \vdots & \vdots & \vdots & \vdots & 0 & 0 & 0 &
         \vdots & \vdots\\
    \vdots & \vdots & \vdots & \vdots & \vdots & \vdots & \vdots & 0 & 0 &
         \ddots & \vdots \\
    \vdots & \vdots & \vdots & \vdots & \vdots & \vdots & \vdots & \vdots &
         \vdots & \vdots & \ddots\\
   \end{array}
\right]
\end{displaymath}
 Obtemos então:  \\
$\bar{\delta}_{1}(\gamma_{A_{1}}(1,1)) =  Tr(Id - \gamma_{A^{*}_{1}}(1,1) \gamma_{A_{1}}(1,1)) - Tr(Id - \gamma_{A_{1}}(1,1) \gamma_{A^{*}_{1}}(1,1)) = 1$,\\
$\bar{\delta}_{1}(\gamma_{A_{1}}(1,-1)) = Tr(Id - \gamma_{A^{*}_{1}}(1,-1) \gamma_{A_{1}}(1,-1)) - Tr(Id - \gamma_{A_{1}}(1,-1) \gamma_{A^{*}_{1}}(1,-1)) = 0$,\\
$\bar{\delta}_{1}(\gamma_{A_{2}}(1,1)) =  Tr(Id - \gamma_{A^{*}_{2}}(1,1) \gamma_{A_{2}}(1,1)) - Tr(Id - \gamma_{A_{2}}(1,1) \gamma_{A^{*}_{2}}(1,1)) = 0$,\\
$\bar{\delta}_{1}(\gamma_{A_{2}}(1,-1)) = Tr(Id - \gamma_{A^{*}_{2}}(1,-1) \gamma_{A_{2}}(1,-1)) - Tr(Id - \gamma_{A_{2}}(1,-1) \gamma_{A^{*}_{2}}(1,-1)) = 1$,\\
$\bar{\delta}_{1}(\gamma_{A_{3}}(1,1)) = Tr(Id - \gamma_{A^{*}_{3}}(1,1) \gamma_{A_{3}}(1,1)) - Tr(Id - \gamma_{A_{3}}(1,1) \gamma_{A^{*}_{3}}(1,1)) = -1$,\\
$\bar{\delta}_{1}(\gamma_{A_{3}}(1,-1)) = Tr(Id - \gamma_{A^{*}_{3}}(1,-1) \gamma_{A_{3}}(1,-1)) - Tr(Id - \gamma_{A_{3}}(1,-1) \gamma_{A^{*}_{3}}(1,-1)) = 0$,\\
$\bar{\delta}_{1}(\gamma_{A_{4}}(1,1)) = Tr(Id - \gamma_{A^{*}_{4}}(1,1) \gamma_{A_{4}}(1,1)) - Tr(Id - \gamma_{A_{4}}(1,1) \gamma_{A^{*}_{4}}(1,1)) = 0$,\\
$\bar{\delta}_{1}(\gamma_{A_{4}}(1,-1)) = Tr(Id - \gamma_{A^{*}_{4}}(1,-1) \gamma_{A_{4}}(1,-1)) - Tr(Id - \gamma_{A_{4}}(1,-1) \gamma_{A^{*}_{4}}(1,-1)) = -1$.

Portanto podemos concluir que \\
\centerline{$\eta \circ \delta_{1}([[A_{1}]_{\mathcal{E}}]_{1}) = (1,0)$,}\\
\centerline{$\eta \circ \delta_{1}([[A_{2}]_{\mathcal{E}}]_{1}) = (0,1)$,}\\
\centerline{$\eta \circ \delta_{1}([[A_{3}]_{\mathcal{E}}]_{1}) = (-1,0)$,}\\
\centerline{$\eta \circ \delta_{1}([[A_{4}]_{\mathcal{E}}]_{1}) = (0,-1)$.}

Segue da exatidão da sequência exata cíclica de seis termos da página 59 que $\mbox{Im}\,
\delta_{1}= \mbox{Ker}\, \varphi_{*} \cong \mathbb{Z} \oplus \mathbb{Z}$;
$\mbox{Im}\, \varphi_{*} = \mbox{Ker}\, \psi_{*} = 0$ e $\mbox{Im}\, \psi_{*} =
\mbox{Ker}\, \delta_{0} = \mathbb{Z}(1,1)$ onde $(1,1) \in K_{0}(\mathcal{A}/
\mathcal{E}) = \mathbb{Z}[[b(D)]_{\mathcal{E}}]_{0} \oplus
\mathbb{Z}[[c(D)]_{\mathcal{E}}]_{0}$ (veja página 55, (3.1)). Portanto
$$K_{0}(\mathcal{A}/ \mathcal{K}_{\mathbb{R}}) =
\mathbb{Z}[[Id]_{\mathcal{K}_{\mathbb{R}}}]_{0}.$$

Pelos cálculos acima temos que $\mbox{Ker}\,\delta_{1} = \{(x, y, z, w) \in
\mathbb{Z} \oplus \mathbb{Z} \oplus \mathbb{Z} \oplus \mathbb{Z}: y = w \quad
\mbox{e} \quad x = z \}$ logo $(1, 0, 1, 0)$ e $(0, 1, 0, 1) \in K_{1}(\mathcal{A}/\mathcal{E}) = \mathbb{Z}[[A_{1}]_{\mathcal{E}}]_{1} \oplus
\mathbb{Z}[[A_{2}]_{\mathcal{E}}]_{1} \oplus
\mathbb{Z}[[A_{3}]_{\mathcal{E}}]_{1} \oplus
\mathbb{Z}[[A_{4}]_{\mathcal{E}}]_{1} $ (veja página 56, (3.2)) geram
$\mbox{Ker}\,\delta_{1}$.

Novamente pela exatidão da sequência exata cíclica de seis termos temos que
$\mbox{Im}\, \psi_{*} = \mbox{Ker}\, \delta_{1} \cong \mathbb{Z} \oplus
\mathbb{Z}$. Pelo Teorema do homomorfismo e exatidão da sequência exata
cíclica temos que  $\mbox{Ker}\, \psi_{*} = \mbox{Im}\, \varphi_{*}
\cong \mathbb{Z} \oplus \mathbb{Z} / \mathbb{Z} (1, 1) \cong \mathbb{Z}
(1, 0)$ onde $(1,0) \in K_{1}(\mathcal{E}/ \mathcal{K}_{\mathbb{R}})$.

Podemos então construir a seguinte sequência exata curta 
\begin{center}
\begin{picture}(350,15)
\put(31,0){0}
\put(47,3){\vector(1,0){25}}
\put(57,0){\makebox[60pt]{$\mathbb{Z}$}}
\put(102,3){\vector(1,0){35}}
\put(146,0){\makebox[60pt]{$K_{1}(\mathcal{A}/ \mathcal{K}_{\mathbb{R}})$}}
\put(215,3){\vector(1,0){35}}
\put(250,0){\makebox[60pt]{$ \mathbb{Z} \oplus \mathbb{Z}$}}
\put(307,3){\vector(1,0){25}}
\put(341,0){0}
\put(102,8){\makebox[35pt]{$ i$}}
\put(215,8){\makebox[35pt]{$\psi_{*} $}}
\end{picture}
\end{center}
onde $i$ é a inclusão.

Pela Proposição 1.25 e demonstração do Teorema 3.5 temos que $[b(M) e^{2 \pi i
  c(D)} + c(M)]_{\mathcal{K}_{\mathbb{R}}}$, $b$ e $c$ como na página 54,
  corresponde ao $(1, 0)$ de $K_{1}(\mathcal{E}/
  \mathcal{K}_{\mathbb{R}})$. Então usando a sequência exata acima temos que 
 $$K_{1}(\mathcal{A}/ \mathcal{K}_{\mathbb{R}}) =
 \mathbb{Z}[[L(M)]_{\mathcal{K}_{\mathbb{R}}}]_{1} \oplus
 \mathbb{Z}[[\tilde{L}(M)]_{\mathcal{K}_{\mathbb{R}}}]_{1} \oplus
 \mathbb{Z}[[b(M) e^{2 \pi i  c(D)} + c(M)]_{\mathcal{K}_{\mathbb{R}}}]_{1}.$$
 $\cqd$

Agora já estamos prontos para calcular a K-teoria da $C^*$- álgebra $\mathcal{A}$.
\begin{Teo}
Seja $\pi: \mathcal{A} \longrightarrow  \mathcal{A}/
\mathcal{K}_{\mathbb{R}}$ a projeção canônica. Então $\pi_{*}: K_{i}( \mathcal{A})
\longrightarrow \mbox{Im}\, \pi_{*}$ é um isomorfismo. Logo
$K_{0}(\mathcal{A}) \cong \mathbb{Z}$ e $K_{1}(\mathcal{A}) \cong \mathbb{Z}
\oplus \mathbb{Z}$.
\end{Teo}
\textbf{Demonstração:}

 Consideremos a sequência exata curta
\begin{center}
\begin{picture}(280,15)
\put(15,0){0}
\put(30,3){\vector(1,0){25}}
\put(43,0){\makebox[60pt]{$\mathcal{K}_{\mathbb{R}}$}}
\put(90,3){\vector(1,0){35}}
\put(110,0){\makebox[60pt]{$\mathcal{A}$}}
\put(152,3){\vector(1,0){35}}
\put(181,0){\makebox[60pt]{$\mathcal{A}/ \mathcal{K}_{\mathbb{R}}$}}
\put(230,3){\vector(1,0){25}}
\put(262,0){0}
\put(90,8){\makebox[35pt]{$i $}}
\put(152,8){\makebox[35pt]{$\pi$}}
\end{picture}
\end{center}
e sua sequência exata ciclíca de seis termos para K-teoria

\begin{picture}(350,15)
\put(48,0){\makebox[60pt]{$K_{0}(\mathcal{K}_{\mathbb{R}})$}}
\put(111,3){\vector(1,0){35}}
\put(149,0){\makebox[60pt]{$K_{0}(\mathcal{A})$}}
\put(212,3){\vector(1,0){35}}
\put(250,0){\makebox[60pt]{$K_{0}(\mathcal{A}/ \mathcal{K}_{\mathbb{R}})$}}
\put(111,8){\makebox[35pt]{$ i_{*}$}}
\put(212,8){\makebox[35pt]{$ \pi_{*} $}}
\end{picture}

\begin{picture}(350,45)
\put(78,0){\vector(0,1){35}}
\put(280,35){\vector(0,-1){35}}
\put(81,17){$\delta_{1}  $}
\put(283,17){$ \delta_{0} $}
\end{picture}

\begin{picture}(350,15)
\put(48,0){\makebox[60pt]{$K_{1}(\mathcal{A}/ \mathcal{K}_{\mathbb{R}})$}}
\put(146,3){\vector(-1,0){35}}
\put(149,0){\makebox[60pt]{$K_{1}(\mathcal{A})$}}
\put(247,3){\vector(-1,0){35}}
\put(250,0){\makebox[60pt]{$K_{1}( \mathcal{K}_{\mathbb{R}})$}}
\put(111,8){\makebox[35pt]{$ \pi_{*}$}}
\put(212,8){\makebox[35pt]{$ i_{*} $}}
\end{picture}

Substituindo as K-teorias conhecidas temos

\begin{picture}(350,15)
\put(48,0){\makebox[60pt]{$\mathbb{Z}$}}
\put(111,3){\vector(1,0){35}}
\put(149,0){\makebox[60pt]{$K_{0}(\mathcal{A})$}}
\put(212,3){\vector(1,0){35}}
\put(250,0){\makebox[60pt]{$\mathbb{Z}$}}
\put(111,8){\makebox[35pt]{$ i_{*}$}}
\put(212,8){\makebox[35pt]{$ \pi_{*} $}}
\end{picture}

\begin{picture}(350,45)
\put(78,0){\vector(0,1){35}}
\put(280,35){\vector(0,-1){35}}
\put(81,17){$\delta_{1}  $}
\put(283,17){$ \delta_{0} $}
\end{picture}

\begin{picture}(350,15)
\put(48,0){\makebox[60pt]{$\mathbb{Z} \oplus \mathbb{Z} \oplus \mathbb{Z}$}}
\put(146,3){\vector(-1,0){35}}
\put(149,0){\makebox[60pt]{$K_{1}(\mathcal{A})$}}
\put(247,3){\vector(-1,0){35}}
\put(250,0){\makebox[60pt]{$0$}}
\put(111,8){\makebox[35pt]{$ \pi_{*}$}}
\put(212,8){\makebox[35pt]{$ i_{*} $}}
\end{picture}

Sabemos pela Proposição 9.4.4 de \cite{[R]} que o índice para esta sequência é
o índice de Fredholm. Vamos calcular $\delta_{1}$ nos geradores de
$K_{1}(\mathcal{A} / \mathcal{K}_{\mathbb{R}})$.

Temos que $b(M) e^{2 \pi i c(D)} + c(M) \in \mathcal{C}$, logo $[b(M) e^{2 \pi
  i c(D)} + c(M)]_{\mathcal{K}_{\mathbb{R}}} \in \mathcal{C} /
\mathcal{K}_{\mathbb{R}}$ e como foi visto no Capítulo 2, Teorema 2.13, $\delta_{1}([[T]_{\mathcal{K}_{\mathbb{R}}}]_{1}) = ind(T) = w(\sigma_{T})$ nesta álgebra.

Logo $\delta_{1}([[b(M) e^{2 \pi i c(D)} + c(M)]_{\mathcal{K}_{\mathbb{R}}}]_{1}) =
w(\sigma_{b(M) e^{2 \pi i c(D)} + c(M)}) = -1$. Sabemos
também que $ind(L(M)) = ind(\tilde{L}(M)) = 0$, pois são operadores de
multiplicação. Portanto $\delta_{1}$ é sobrejetora e $\delta_{1}(x, y, z) = -z$.

Pela exatidão da sequência exata cíclica de seis termos temos que $\mbox{Im}\,
i_{*} = \mbox{Ker} \, \pi_{*} = 0$; $\mbox{Im}\,
\pi_{*} = \mbox{Ker} \, \delta_{0} = \mathbb{Z}$; $\mbox{Im}\,
i_{*} = \mbox{Ker} \, \pi_{*} = 0$ e $\mbox{Im}\, \pi_{*} = \mbox{Ker} \,
\delta_{1} \cong \mathbb{Z} \oplus \mathbb{Z}$, logo $$K_{0}(\mathcal{A}) =
\mathbb{Z}[Id]_{0}\hspace{0.5cm} \mbox {e} \hspace{0.5cm} K_{1}(\mathcal{A}) =
\mathbb{Z}[L(M)]_{1} \oplus \mathbb{Z}[\tilde{L}(M)]_{1}.\qquad \cqd $$ 
\chapter{K-teoria de $\mathcal{A}$ por produto cruzado}
\hspace{8mm}Neste capítulo calcularemos, novamente, a K-teoria de
$\mathcal{A}$ utilizando a seguinte estratégia. Primeiro provaremos que a
$\mbox{Im}\, \gamma$, onde $\gamma$ é a aplicação da Proposição 1.25 tem
estrutura de produto cruzado e então usando a sequência exata de
Pimsner-Voiculescu \cite{[PV]} calcularemos sua K-teoria. Utilizando então a sequência exata curta
induzida por $\gamma$
 \begin{center}
\begin{picture}(350,15)
\put(10,0){0}
\put(20,3){\vector(1,0){25}}
\put(48,0){\makebox[60pt]{$J_{0}/ \mathcal{K}_{\mathbb{R}}$}}
\put(111,3){\vector(1,0){35}}
\put(149,0){\makebox[60pt]{$\mathcal{A} / \mathcal{K}_{\mathbb{R}}$}}
\put(212,3){\vector(1,0){35}}
\put(250,0){\makebox[60pt]{$\mbox{Im} \, \gamma$}}
\put(313,3){\vector(1,0){25}}
\put(341,0){0}
\put(111,8){\makebox[35pt]{}}
\put(212,8){\makebox[35pt]{}}
\end{picture}
\end{center}
calculamos a K-teoria de $\mathcal{A} / \mathcal{K}_{\mathbb{R}}$. De posse
desta K-teoria e usando a sequência exata curta abaixo calculamos a K-teoria
de $\mathcal{A}$.
\begin{center}
\begin{picture}(280,15)
\put(15,0){0}
\put(30,3){\vector(1,0){25}}
\put(43,0){\makebox[60pt]{$\mathcal{K}_{\mathbb{R}}$}}
\put(90,3){\vector(1,0){35}}
\put(110,0){\makebox[60pt]{$\mathcal{A}$}}
\put(152,3){\vector(1,0){35}}
\put(181,0){\makebox[60pt]{$\mathcal{A}/ \mathcal{K}_{\mathbb{R}}$}}
\put(230,3){\vector(1,0){25}}
\put(262,0){0}
\put(90,8){\makebox[35pt]{}}
\put(152,8){\makebox[35pt]{}}
\end{picture}.
\end{center}
O que diferencia a estratégia deste capítulo para o anterior é que no Capítulo
3 utilizamos uma sequência exata curta induzida pelo $\sigma$-símbolo para
calcular a K-teoria de $\mathcal{A}/ \mathcal{K}_{\mathbb{R}}$.

Um estudo geral sobre produto cruzado de $C^*$- álgebras pode ser encontrado no Capítulo 7 de \cite{[Pe]}. Aqui, trabalharemos apenas em casos onde a $C^*$- álgebra $A$ tem unidade e o grupo é $\mathbb{Z}$ com a topologia discreta.

Dados uma $C^*$- álgebra $A$ com unidade, o grupo $\mathbb{Z}$ com a topologia
discreta, e uma ação $\alpha$ de $\mathbb{Z}$ em $A$, define-se uma nova $C^*$- álgebra $A \underset{\alpha}{\rtimes} \mathbb{Z}$ como sendo a $C^*$- álgebra envolvente de $l_{1}(\mathbb{Z}, A)$ (seção (2.7) de \cite{[Dx]}), esta vista como uma álgebra de Banach com a involução, o produto e a norma definidos respectivamente por:\\[0.2cm]
\centerline{$y^*(n) = \alpha^{n}(y^{*}_{-n})$,$\quad \alpha \in Aut(A)$,}\\[0.2cm]
\centerline{$y \times z(n) = \dps{\sum_{k \in \mathbb{Z}}} y_{k}
  \alpha^{k}(z_{n-k})$ \hspace{0.2cm} e}\\[0.2cm]
\centerline{$||y||_{1} = \dps{\sum_{n \in \mathbb{Z}}}|y_{n}|$}\\[0.2cm]
para todo $y = (y_{n})_{n \in \mathbb{Z}}$ e $z = (z_{n})_{n \in
  \mathbb{Z}}$ em $l_{1}( \mathbb{Z}, A)$ e $n \in \mathbb{Z}$. Temos que $k( \mathbb{Z}, A)$, funções de $l_{1}( \mathbb{Z}, A)$ que tem suporte finito, é um subconjunto fortemente denso de $l_{1}( \mathbb{Z}, A)$.

Sejam $a \in A$ e $n \in \mathbb{Z}$. Denote por $a \delta_{n}$ o elemento de $l_{1}( \mathbb{Z}, A)$ dado por
\begin{displaymath}
 a \delta_{n} (k) = \left\{ 
 \begin{array}{lll}
   0 & se & k \neq n,\\
   a & se & k = n.\\
 \end{array}
 \right.
\end{displaymath}
Note que $(a \delta_{n})^{*} = \alpha^{-n}(a^*) \delta_{-n}$ e que $(a \delta_{n}) \times (b \delta_{m}) = a \alpha^{n}(b) \delta_{n+m}$ para todo $n, m \in \mathbb{Z}$ e todo $a, b \in A$. E não é difícil mostrar que o conjunto  $\{ x \delta_{0}, 1 \delta_{1}: x \in A \}$ gera $l_{1}(\mathbb{Z},A)$ como álgebra de Banach.
\section{ A estrutura de produto cruzado de $\mathcal{A}^{\diamond}$}

Seja $A = CS( \mathbb{R})$, e consideremos a sequência exata curta
\begin{center}
\begin{picture}(350,15)
\put(10,0){0}
\put(26,3){\vector(1,0){25}}
\put(51,0){\makebox[60pt]{$C( \mathbb{C})$}}
\put(111,3){\vector(1,0){35}}
\put(149,0){\makebox[60pt]{$C([0,1])$}}
\put(212,3){\vector(1,0){35}}
\put(235,0){\makebox[60pt]{$ \mathbb{C}$}}
\put(283,3){\vector(1,0){25}}
\put(323,0){0}
\put(111,8){\makebox[35pt]{$i$}}
\put(212,8){\makebox[35pt]{$\psi $}}
\end{picture}
\end{center}
onde $C(\mathbb{C})$ é o cone da álgebra $\mathbb{C}$, isto é $C( \mathbb{C}) = \{f: [0,1] \rightarrow \mathbb{C}\;:\; f(0) = 0\}$; $\psi(f) = f(0)$ e $i$ é a inclusão.
Aplicando a sequência exata cíclica de seis termos para K-teoria na sequência
acima e usando que $K_{0}(C( \mathbb{C})) = K_{1}(C( \mathbb{C})) = 0$, veja
exemplos 4.1.5 e 8.3 de \cite{[R]}, obtemos $K_{0}(C([0,1])) =
\mathbb{Z}[1]_{0}$ e $K_{1}(C([0,1])) = 0$. Como $A$ é isomorfo a
$C([0,1])$ temos que $K_{1}(A) = 0$ e $K_{0}(A) = \mathbb{Z}[1]_{0}$.

Consideremos $\alpha: \mathbb{Z} \longrightarrow Aut(A)$ onde $n \longmapsto [\alpha(n)a](x) = a(x-n)$. A aplicação $\alpha$ está bem definida e temos que $[ \alpha (n + m) (a)](x) = a(x - n - m) = [ \alpha(n) \circ \alpha(m)(a)](x)$ para todo $n ,m \in \mathbb{Z}$. Logo $\alpha$ é um homomorfismo de grupo.
Como $\mathbb{Z}$ é um grupo discreto temos que $\alpha$ é contínua, logo $\alpha$ é uma ação.

Usando a sequência exata de Pimsner-Voiculescu \cite{[PV]} temos
\newpage
\begin{picture}(350,15)
\put(48,0){\makebox[60pt]{$K_{0}(A)$}}
\put(111,3){\vector(1,0){35}}
\put(149,0){\makebox[60pt]{$K_{0}(A)$}}
\put(212,3){\vector(1,0){35}}
\put(250,0){\makebox[60pt]{$K_{0}(A \underset{\alpha}{\rtimes} \mathbb{Z} )$}}
\put(111,8){\makebox[35pt]{$ id - \alpha_{*}^{-1}$}}
\put(212,8){\makebox[35pt]{$ i_{*} $}}
\end{picture}

\begin{picture}(350,45)
\put(78,0){\vector(0,1){35}}
\put(280,35){\vector(0,-1){35}}
\put(81,17){$\delta^{'}_{1}  $}
\put(283,17){$ \delta^{'}_{0} $}
\end{picture}

\begin{picture}(350,15)
\put(48,0){\makebox[60pt]{$K_{1}(A \underset{\alpha}{\rtimes} \mathbb{Z})$}}
\put(146,3){\vector(-1,0){35}}
\put(149,0){\makebox[60pt]{$K_{1}(A)$}}
\put(247,3){\vector(-1,0){35}}
\put(250,0){\makebox[60pt]{$K_{1}(A)$}}
\put(111,8){\makebox[35pt]{$i_{*}$}}
\put(212,8){\makebox[35pt]{$ id - \alpha_{*}^{-1}$}}
\end{picture} 
 
Como $\alpha(n)(1)(x) = 1(x-n) = 1$, segue que $\alpha_{*}$ é igual a
identidade de $K_{0}(A)$. Usando que $K_{1}(A) = 0$ e que $K_{0}(A) =
\mathbb{Z}[1]_{0}$ ficamos com a seguinte sequência exata

\centerline{$ 0 \quad \longrightarrow \quad K_{1}(A \underset{\alpha}{\rtimes}
  \mathbb{Z})  \quad {\buildrel \delta^{'}_{1} \over \longrightarrow} \quad \mathbb{Z}[1]_{0} \quad {\buildrel 0 \over \longrightarrow} \quad \mathbb{Z}[1]_{0} \quad  {\buildrel i_{*} \over \longrightarrow} \quad K_{0}(A \underset{\alpha}{\rtimes} \mathbb{Z}) \quad  \longrightarrow \quad 0$}

Podemos concluir então que $$K_{1}(A \underset{\alpha}{\rtimes} \mathbb{Z})
\cong \mathbb{Z}\hspace{.5cm} \mbox{e}$$ $$K_{0}(A \underset{\alpha}{\rtimes}
\mathbb{Z}) = \mathbb{Z}[1]_{0}.$$

Se $B$ é uma $C^*$- álgebra e $\alpha \in Aut(B)$ sempre existe um unitário $u
\in B \underset{\alpha}{\rtimes} \mathbb{Z}$ tal que $\alpha(a) = uau^*$ para
todo $a \in B$ (veja, por exemplo, \cite{[Cris]}, pág. 49). Diz-se que um tal
$u$ implementa a ação $\alpha$ e que $\alpha$ é o automorfismo induzido por
$u$. No nosso caso, temos que $\alpha(n)a = \delta_{1}^na(M)(\delta_{1}^*)^n$ para todo
$a \in A$ e $n \in \mathbb{Z}$, isto é $\alpha$ é o automorfismo induzido pelo
unitário $\delta_{1} \in A \underset{\alpha}{\rtimes} \mathbb{Z}$.

Pela equação apresentada no começo da página 102 de \cite{[PV]} temos que
$\delta^{'}_{1}([(1 \otimes 1_{n} - F) + Fx(u^* \otimes 1_{n})F]_{1}) = [F]_{0}$,
onde $F, x \in A \otimes M_{n}$, $n \in \mathbb{N}$; $F$ é uma projeção e $u$
é o unitário que implementa a ação. Consideremos $F = 1_{A} = x$ e $n = 1$, pela equação
acima temos que $\delta^{'}_{1}([\delta_{-1}]_{1}) = [1]_{0}$. E como foi visto acima que $\delta^{'}_{1}$ é
um isomorfismo temos que $$K_{1}(A \underset{\alpha}{\rtimes} \mathbb{Z}) = \mathbb{Z}[\delta_{-1}]_{1}.$$

Consideremos a $C^*$- álgebra $\mathcal{A}^{\diamond}$ definida na Proposição 1.26. Vamos verificar na Proposição abaixo que $\mathcal{A}^{\diamond}$ é isomorfo a $A \underset{\alpha}{\rtimes} \mathbb{Z}$.
\begin{Prop}
Seja $A = CS( \mathbb{R})$, então existe $ \varphi: A \underset{\alpha}{\rtimes} \mathbb{Z} \longrightarrow \mathcal{A}^{\diamond}$ isomorfismo.
\end{Prop}
\textbf{Demonstração:}

Seja $B$ o fecho do conjunto $B_{0} = \{\; \dps{\sum_{j \in F}} a_{j}(M)T_{j}\;:\; a_{j} \in
CS( \mathbb{R})\; \}$, onde $F \subset \mathbb{Z}$ finito.

Daqui até o fim dessa demonstração estaremos considerando $T_{j} u(x) = u(x
- j)$.

Consideremos $f: k(\mathbb{Z},A) \longrightarrow B_{0}$
onde $(a_{j})_{j} \longmapsto \dps{\sum_{j \in F}} a_{j}(M)T_{j}$, $F \subset
\mathbb{Z}$ finito. Vamos verificar que $f$ é um $*$ - homomorfismo contínuo:\\
$\triangleright$ Produto\\
$f((a_{j})_{j})f((b_{j})_{j}) = \dps{\sum_{j}}\dps{\sum_{k}} a_{j}(M)T_{j}b_{k}(M)T_{k}$. Temos que $T_{j}b_{k}(M) = b_{k}(M-j)T_{j}$, então\\
$f((a_{j})_{j})f((b_{j})_{j}) = \dps{\sum_{j}}\dps{\sum_{k}}a_{j}(M)b_{k}(M-j)T_{j+k}$. Fazendo mudanças de variáveis e usando a definição de $\alpha$ temos
$f((a_{j})_{j})f((b_{j})_{j}) = \dps{\sum_{j}}( \dps{\sum_{l}}a_{l} \alpha^{l}(b_{j-l}))T_{j} = f((a_{j} * b_{j})_{j})$.\\[0.2cm]
 $\triangleright$ Involução\\
$f((a_{j})_{j})^* = \dps{\sum_{j}}T_{-j}a^*_{j}(M) = \sum T_{j'}a^*_{-j'}(M) = \sum \alpha^{j'}(a^*_{-j'})T_{j'} = f((a_{j})^*)$.\\[0.2cm]
$\triangleright$ Continuidade\\
$||f((a_{j}))_{j}|| = || \dps{\sum_{j}}a_{j}(M)T_{j}|| \leq \dps{\sum_{j}}||a_{j}(M)|| ||T_{j}|| = \dps{\sum_{j}}||a_{j}(M)|| = ||(a_{j})_{j}||_{1}.$

Como $f$ é contínua, $k(\mathbb{Z},A)$ é denso em $l_{1}(\mathbb{Z},A)$ e $B_{0}$ é denso em $B$ temos que existe $*$ - homomorfismo $\tilde{f}: l_{1}(\mathbb{Z},A) \longrightarrow B$ extensão contínua de $f$.
\begin{displaymath}
 \begin{array}{ccc}
  k(\mathbb{Z},A) & {\buildrel f \over \rightarrow} & B_{0}\\
  \cap & & \cap\\
  l_{1}(\mathbb{Z},A) & {\buildrel \tilde{f} \over \rightarrow} & B\\
 \end{array}
\end{displaymath}

Por 2.7.4 de \cite{[Dx]} temos que existe
$\varphi: A \underset{\alpha}{\rtimes} \mathbb{Z} \longrightarrow B$ $*$ -
homomorfismo contínuo estendendo $\tilde{f}$.
\begin{displaymath}
 \begin{array}{ccc}
  l_{1}(\mathbb{Z},A) & {\buildrel \tilde{f} \over \rightarrow} & B\\
  \cap & \nearrow \varphi & \\
  A \underset{\alpha}{\rtimes} \mathbb{Z} & & \\
 \end{array}
\end{displaymath}
Vamos verificar agora que $\varphi$ é um isomorfismo.
\vspace{0.3cm}\\
$\triangleright$ $\varphi$ é sobrejetora

Como $\mbox{Im} \, \varphi$ é fechada e contém $B_{0}$ que é denso em $B$
temos que $\varphi$ é sobrejetora.

Para provar a injetividade de $\varphi$ vamos utilizar a Proposição 2.9 de
\cite{[Ex]}, que enunciamos abaixo, 
\begin{Prop}
Sejam $C$ e $C^{'}$ $C^*$- álgebras e sejam $\alpha$ e $\alpha^{'}$ ações de $S^1$ em $C$ e $C^{'}$ respectivamente. Suponha que $\psi: C \longrightarrow C^{'}$ é um homomorfismo covariante. Se a restrição de $\psi$ à subálgebra dos pontos fixos $C_{0}$ é injetiva, então $\psi$ é injetiva.
\end{Prop}  

Então vamos definir ações $\gamma$ e $\beta$ de $S^1$ em $A
\underset{\alpha}{\rtimes} \mathbb{Z}$ e $B$ respectivamente tal que
$\beta_{z}(\varphi(x)) = \varphi(\gamma_{z}(x))$, para todo $z \in S^1$ e para
todo $x \in A \underset{\alpha}{\rtimes} \mathbb{Z}$.  

Dado $z \in S^1$ e $a = (a_{j})_{j} \in l_{1}(\mathbb{Z},A)$ vamos definir $\gamma_{z}(a) = (z^{j}a_{j})_{j}$.\\
Temos que $\gamma$ é uma ação de $S^{1}$ em $l_{1}(\mathbb{Z},A)$. Seja $i: l_{1}(\mathbb{Z},A) \longrightarrow 
 A \underset{\alpha}{\rtimes} \mathbb{Z}$ a inclusão. Então $i \circ
 \gamma_{z}: l_{1}(\mathbb{Z},A) \longrightarrow  A \underset{\alpha}{\rtimes}
 \mathbb{Z}$ é um $*$-homomorfismo contínuo. Por 2.7.4 de \cite{[Dx]} temos que $\gamma$ é uma ação de $S^{1}$ em $ A \underset{\alpha}{\rtimes} \mathbb{Z}$.

Dado $z \in S^1$ vamos definir o operador $U$ em $\mathcal{L}_{\mathbb{R}}$
por $(U_{z}f)(t) = z^t f(t)$. Temos que $U$ é um operador unitário. Dado $z
\in S^1$ vamos definir $\beta_{z}(A) =U_{z}A(U_{z})^{-1}$ para todo $A \in
B$. Temos então que $\beta$ é uma ação de $S^1$ em $B$.

Agora vamos verificar que $\varphi$ é um homomorfismo covariante, ou seja, vamos mostrar que $\beta_{z}(\varphi(x)) = \varphi(\gamma_{z}(x))$, para todo $x \in l_{1}(\mathbb{Z},A)$.
Como $\{ x \delta_{0}, 1 \delta_{1}: x \in A \}$ gera $l_{1}(\mathbb{Z},A)$, é suficiente verificarmos a equação nesses geradores de $l_{1}(\mathbb{Z},A)$.
Para $x \delta_{0}$ temos que $\gamma_{z}( x \delta_{0}) = x \delta_{0}$ e $\varphi( x \delta_{0}) = f( x \delta_{0}) = x(M)$, logo\\[.1cm]
\centerline{ $\varphi( \gamma_{z}( x \delta_{0})) = x(M) = \beta_{z}( \varphi( x \delta_{0}))$.}
Para $1 \delta_{1}$ temos $\gamma_{z}( 1 \delta_{1}) = z \delta_{1}$ e $\varphi( 1 \delta_{1}) = f( 1 \delta_{1}) = T_{1}$, logo\\[.1cm]
\centerline{ $\varphi( \gamma_{z}( 1 \delta_{1})) = z T_{1} = \beta_{z}( \varphi( 1 \delta_{1}))$.}

Vamos verificar agora quem é a álgebra dos pontos fixos de $A \underset{\alpha}{\rtimes} \mathbb{Z}$.
Seja $C$ a álgebra dos pontos fixos de $A \underset{\alpha}{\rtimes} \mathbb{Z}$, temos que $C \cap l_{1}( \mathbb{Z}, A ) = A$.\\
De fato, consideremos a inclusão $i: A \longrightarrow l_{1}( \mathbb{Z}, A )$ onde $x \longmapsto x \delta_{0}$. Se $x \in A$ então pela inclusão temos que $x \in l_{1}(\mathbb{Z},A)$ e com isso $\gamma_{z}(x) = x$, logo $x \in C \cap l_{1}(\mathbb{Z},A)$.
Seja $x \in C \cap l_{1}( \mathbb{Z}, A )$, temos que $x = \sum a_{j}
\delta_{j}$ então $\gamma_{z}(x) = \sum z^{j}a_{j} \delta_{j}$ e como $x \in C$ também temos $\gamma_{z}(x) = x$. Logo $\sum
z^{j}a_{j} \delta_{j} = \sum a_{j} \delta_{j}$, para todo $z \in S^1$, o que implica que $a_{j}
\delta_{j} = 0$ ou $z^{j} - 1 = 0 $, para todo $z \in S^1$. O segundo caso acontece se, e somente
se, $j = 0$. Portanto $a_{j} = 0$ para $j \neq 0$, logo $x = a_{0} \delta_{0}$, isto é, $ x \in A$.

Consideremos agora a aplicação $E: A \underset{\alpha}{\rtimes} \mathbb{Z}
\longrightarrow  A \underset{\alpha}{\rtimes} \mathbb{Z}$ definida por\\
\centerline{$E(b) =  \dps{\int_{0}^{2 \pi}} \gamma_{z}(b)\; dz$, $ b \in A \underset{\alpha}{\rtimes} \mathbb{Z}$ e $dz = \frac{d \theta}{2 \pi}$.}

Se $b$ é tal que $ \gamma_{z}(b) = b$ para todo $z \in S^1$ então $E(b) =
\dps{\int_{0}^{2 \pi}} b \;dz = b$, logo $C \subset \mbox{Im} \, E$. 
Seja $b = (b_{j})_{j} \in l_{1}( \mathbb{Z}, A)$ então\\
\centerline{$E(b) = \dps{\int_{0}^{2 \pi}} \gamma_{z}(b) \;dz =  \dps{\int_{0}^{2 \pi}}
(z^{j}b_{j})_{j} \;dz = (b_{j} \dps{\int_{0}^{2 \pi}} z^{j} \;dz)_{j} = b_{0}
\delta_{0} \in A$.}

Então para todo $a \in l_{1}(\mathbb{Z}, A)$ temos que $E(a) \in A$. Como
$l_{1}(\mathbb{Z}, A)$ é uma subálgebra densa em $A \underset{\alpha}{\rtimes}
\mathbb{Z}$, temos a imagem de $E$ contida em $A$. Portanto temos  $C  \subset
\mbox{Im} \, E \subset A = C \cap l_{1}(\mathbb{Z}, A)$, logo $C = A$.

Resta-nos provar que $\varphi$ é injetora na álgebra dos pontos fixos.
Consideremos a inclusão $i: A \longrightarrow A \underset{\alpha}{\rtimes} \mathbb{Z}$ onde $a \longmapsto a \delta_{0}$. Temos que $ \varphi( a \delta_{0}) = f( a \delta_{0}) = a(M)$ e $\varphi( a \delta_{0}) = \varphi( b \delta_{0}) \Leftrightarrow a(M) = b(M) \Leftrightarrow a = b$.\\
Portanto $\varphi$ é injetora na álgebra dos pontos fixos, logo pela
Proposição 4.2 temos que $\varphi$ é injetora. Logo $\varphi: A
\underset{\alpha}{\rtimes} \mathbb{Z} \longrightarrow B$ é um isomorfismo.

Como $\mathcal{A}^{\diamond} = FBF^{-1}$ temos que $A
\underset{\alpha}{\rtimes} \mathbb{Z}$ é isomorfo a
$\mathcal{A}^{\diamond}$.$\cqd$ 

Na página 68 vimos que $K_{0}(A
\underset{\alpha}{\rtimes} \mathbb{Z}) = \mathbb{Z}[1]_{0}$ e que $K_{1}(A
\underset{\alpha}{\rtimes} \mathbb{Z}) = \mathbb{Z}[\delta_{-1}]_{1}$, logo $$
K_{0}(\mathcal{A}^{\diamond}) = \mathbb{Z}[Id]_{0} \hspace{.5cm} \mbox{e}$$
$$K_{1}(\mathcal{A}^{\diamond}) = \mathbb{Z}[e^{i M}]_{1}.$$

\section{A K-teoria de $\mathcal{A}/ \mathcal{K}_{\mathbb{R}}$ e
  $\mathcal{A}$}
\hspace{8mm}Temos que o conjunto
\begin{equation}
\mathcal{A}_{0} = \Big \{ \dps{ \sum_{j= -N}^{N}}a_{j}(M)b_{j}(D)e^{ijM} + K
\mbox{;} \,  N \in \mathbb{N}, \, a_{j} \in CS( \mathbb{R}), \, b_{j} \in
CS(\mathbb{R}), \, K \in \mathcal{K}_{\mathbb{R}} \Big \}
\end{equation}
é denso em $\mathcal{A}$. Dado  $A = \dps{ \sum_{j=
    -N}^{N}}a_{j}(M)b_{j}(D)e^{ijM} \in \mathcal{A}_{0}$, sejam $A^{\pm} =
\dps{ \sum_{j= -N}^{N}}a_{j}(\pm \infty)b_{j}(D)e^{ijM}$. Temos que $A^{+}$ e
$A^{-} \in \mathcal{A}^{\diamond}$ logo pela demonstração da Proposição 1.26 e pela
definição de $\gamma$ temos que $\gamma_{A} = ( WF^{-1}A^{+}(WF^{-1})^{-1},
WF^{-1}A^{-}(WF^{-1})^{-1})$.

Como isso vale para todo $A$ em uma subálgebra densa de $\mathcal{A}$, segue
que a aplicação $$\mathcal{A}_{0} \ni A \longmapsto (A^{+}, A^{-}) \in
\mathcal{A}^{\diamond} \oplus \mathcal{A}^{\diamond}$$ se estende a um
$*$-homomorfismo $$\psi: \mathcal{A} \longrightarrow \mathcal{A}^{\diamond}
\oplus \mathcal{A}^{\diamond}$$ tal que, para todo $A \in \mathcal{A}$, e com
$U =  \left(   
\begin{array}{cc}
   WF^{-1} & 0\\
   0   & WF^{-1}\\
\end{array}
 \right )$ temos  $U \psi(A) U^{-1} = \gamma(A)$. Segue então que $\mbox{Ker}
 \, \psi = \mbox{Ker} \, \gamma = J_{0}$ (veja Proposição 1.27).
\begin{Prop}
$K_{0}(\mathcal{A}/ \mathcal{K}_{\mathbb{R}})$ é isomorfo a
$\mathbb{Z}$ e $ K_{1}(\mathcal{A}/ \mathcal{K}_{\mathbb{R}})$ é isomorfo a $\mathbb{Z} \oplus \mathbb{Z} \oplus \mathbb{Z} $.
\end{Prop}
\textbf{Demonstração:}

Consideremos a seguinte sequência exata curta:
\begin{center}
\begin{picture}(350,15)
\put(16,0){0}
\put(35,3){\vector(1,0){25}}
\put(48,0){\makebox[60pt]{$J_{0}$}}
\put(96,3){\vector(1,0){35}}
\put(119,0){\makebox[60pt]{$\mathcal{A}$}}
\put(167,3){\vector(1,0){35}}
\put(220,0){\makebox[60pt]{$\mathcal{A}^{\diamond} \oplus \mathcal{A}^{\diamond}$}}
\put(283,3){\vector(1,0){25}}
\put(341,0){0}
\put(96,8){\makebox[35pt]{$ \iota$}}
\put(167,8){\makebox[35pt]{$\psi $}}
\end{picture}
\end{center}
onde $\psi(A) = (A^{+}, A^{-})$ com $A$ e $A^{ \pm}$ como acima.

Podemos tomar então a seguinte sequência exata curta
\begin{center}
\begin{picture}(350,15)
\put(10,0){0}
\put(20,3){\vector(1,0){25}}
\put(48,0){\makebox[60pt]{$J_{0}/ \mathcal{K}_{\mathbb{R}}$}}
\put(111,3){\vector(1,0){35}}
\put(149,0){\makebox[60pt]{$\mathcal{A} / \mathcal{K}_{\mathbb{R}}$}}
\put(212,3){\vector(1,0){35}}
\put(250,0){\makebox[60pt]{$\mathcal{A}^{\diamond} \oplus \mathcal{A}^{\diamond}$}}
\put(313,3){\vector(1,0){25}}
\put(341,0){0}
\put(111,8){\makebox[35pt]{$ \iota$}}
\put(212,8){\makebox[35pt]{$\psi $}}
\end{picture}.
\end{center}

Como $J_{0}/ \mathcal{K}_{\mathbb{R}}$ é isomorfo a $C_{0}( \mathbb{R}
\times \{ - \infty, + \infty \})$ (Lema 1.15) e $C_{0}( \mathbb{R}
\times \{ - \infty, + \infty \})$ é isomorfo a $C_{0}( \mathbb{R}) \oplus
C_{0}( \mathbb{R})$ temos que
$K_{1}(J_{0}/ \mathcal{K}_{\mathbb{R}}) \cong \mathbb{Z} \oplus
\mathbb{Z}$ e $K_{0}(J_{0}/ \mathcal{K}_{\mathbb{R}}) = 0$.

Usando o isomorfismo definido no Lema 1.15 temos que $$K_{1}(J_{0}/ \mathcal{K}_{\mathbb{R}}) = \mathbb{Z}[[U]_{\mathcal{K}_{\mathbb{R}}}]_{1} \oplus
\mathbb{Z}[[V]_{\mathcal{K}_{\mathbb{R}}}]_{1}$$ onde $U = e^{2 \pi i b(M)}b(D) + c(D)$ e  $V = e^{2 \pi i b(M)}c(D) + b(D)$ com
$b$ e $c$ as funções definidas na página 54.

Observemos que $U$ e $V$ pertencem a $J_{0} \oplus \mathbb{C} Id \subset
\mathcal{L}(L^2(\mathbb{R}))$, pois (como $b+c \equiv 1$) $Id - U = (e^{2 \pi i
  b(M)} -Id)b(D) = \tilde{b}(M)b(D)$ com $\tilde{b} \in C_{0}(\mathbb{R})$, e
analogamente para $V$. Como, por definição, $K_{1}(J_{0}/
\mathcal{K}_{\mathbb{R}}) = K_{1}(\widetilde{J_{0}/ \mathcal{K}_{\mathbb{R}}})$,
segue que, de fato, $[[U]_{\mathcal{K}_{\mathbb{R}}}]_{1}$ e
$[[V]_{\mathcal{K}_{\mathbb{R}}}]_{1}$ são elementos de $K_{1}(J_{0}/
\mathcal{K}_{\mathbb{R}})$.

Aplicando a sequência exata cíclica de seis termos para K-teoria na sequência
acima obtemos
\newpage
\begin{picture}(395,15)
\put(48,0){\makebox[60pt]{$0$}}
\put(111,3){\vector(1,0){35}}
\put(149,0){\makebox[60pt]{$K_{0}(\mathcal{A}/ \mathcal{K}_{\mathbb{R}})$}}
\put(212,3){\vector(1,0){35}}
\put(250,0){\makebox[60pt]{$\mathbb{Z} \oplus \mathbb{Z}$}}
\put(111,8){\makebox[35pt]{$ \iota_{*}  $}}
\put(212,8){\makebox[35pt]{$ \psi_{*} $}}
\end{picture}

\begin{picture}(350,45)
\put(78,0){\vector(0,1){35}}
\put(280,35){\vector(0,-1){35}}
\put(81,17){$\delta_{1}  $}
\put(283,17){$ \delta_{0} $}
\end{picture}

\begin{picture}(350,15)
\put(48,0){\makebox[60pt]{$\mathbb{Z} \oplus \mathbb{Z}$}}
\put(146,3){\vector(-1,0){35}}
\put(149,0){\makebox[60pt]{$K_{1}(\mathcal{A}/ \mathcal{K}_{\mathbb{R}})$}}
\put(247,3){\vector(-1,0){35}}
\put(250,0){\makebox[60pt]{$\mathbb{Z} \oplus \mathbb{Z}$}}
\put(111,8){\makebox[35pt]{$ \psi_{*} $}}
\put(212,8){\makebox[35pt]{$ \iota_{*} $}}
\end{picture}
 
Logo $\mbox{Im}\, \delta_{1} = 0$ e $\mbox{Ker}\, \delta_{1} = \mathbb{Z} \oplus \mathbb{Z}$.

Seja $\mathcal{C}$ a $C^*$- álgebra do Capítulo 2, neste capítulo vimos que
$\mathcal{C}/ \mathcal{K}_{\mathbb{R}} \cong C(M_{c})$ onde $M_{c} = \{ (x,
\xi) \in [- \infty, + \infty] \times [- \infty, + \infty]\; :\; |x| + | \xi| =
\infty \}$ e que $K_{1}(\mathcal{C} /
\mathcal{K}_{\mathbb{R}}) = \mathbb{Z}[[e^{2 \pi i c(M)} b(D) +
c(D)]_{\mathcal{K}_{\mathbb{R}}}]_{1}$.

Consideremos as inclusões\\[0.2cm]
\centerline{$J_{0}/ \mathcal{K}_{\mathbb{R}} \quad {\buildrel i \over
    \longrightarrow} \quad \mathcal{C}/ \mathcal{K}_{\mathbb{R}} \quad {\buildrel j
    \over \longrightarrow} \quad \mathcal{A}/ \mathcal{K}_{\mathbb{R}}$.}\\[0.2cm]
Temos que $ \iota = j \circ i$.
Pela definição de $\psi$ temos que $\psi_{*}([Id]_{0}) = ([Id]_{0}, [Id]_{0})
$. Logo com respeito ao isomorfismo $K_{0}(\mathcal{A}^{\diamond} \oplus
\mathcal{A}^{\diamond}) \cong \mathbb{Z} \oplus \mathbb{Z}$ podemos escrever
que $(1, 1) \in \mbox{Im}\, \psi_{*} = \mbox{Ker}\, \delta_{0}$ e $\delta_{0}(1, 1) = 0$.

Seja $W = e^{2 \pi i c(M)} c(D) + b(D)$. Calculando o número de rotação de $\sigma_{U}$ e $\sigma_{W}$ temos que
 $w(\sigma_{U}) = w(\sigma_{W}) = -1$. Logo pelo o que vimos na página 47
 sobre o número de rotação temos
 que $i_{*}([[U]_{\mathcal{K}_{\mathbb{R}}}]_{1}) =
 i_{*}([[W]_{\mathcal{K}_{\mathbb{R}}}]_{1}) = [[e^{2 \pi i c(D)} b(M) +
 c(M)]_{\mathcal{K}_{\mathbb{R}}}]_{1}$. Portanto com respeito ao isomorfismo $K_{1}(J_{0}/ \mathcal{K}_{\mathbb{R}}) = \mathbb{Z}[[U]_{\mathcal{K}_{\mathbb{R}}}]_{1} \oplus
\mathbb{Z}[[V]_{\mathcal{K}_{\mathbb{R}}}]_{1}$ temos que $\iota_{*}((1,0) -
 (0,-1)) = j_{*} \circ i_{*}((1,0) - (0,-1)) = 0$. Logo $(1, 1) \in
 \mbox{Ker}\, \iota_{*} = \mbox{Im}\, \delta_{0}$.

Suponhamos que $\delta_{0}(1, 0) = (m, n)$, como $(x, y) = (x-y)(1,0) +
y(1,1)$ temos que $\delta_{0}(x, y) = (x-y)(m, n)$, para todo $(x,y) \in
\mathbb{Z} \oplus \mathbb{Z}$.

Como $(1, 1) \in \mbox{Im}\, \delta_{0}$ temos que existe $(x, y) \in
\mathbb{Z} \oplus \mathbb{Z}$ tal que $(x- y)(m, n) = (1, 1)$, isto é, $(x -
y)m = 1$ e $(x - y)n = 1$. Logo $m = n \neq 0$ e com isto $\mbox{Im}\,
\delta_{0} = \mathbb{Z}(m, n)$ e $\mbox{Ker}\,\delta_{0} = \mathbb{Z}(1,
1)$. Sabemos que $(1,1) \in \mbox{Ker}\, \iota_{*} = \mbox{Im}\,
\delta_{0} \cong \mathbb{Z}$. Logo, $\mbox{Im}\, \delta_{0} = \mbox{Ker}\, \iota_{*} = (1,1)\mathbb{Z}$.

Podemos então concluir pelo Teorema do homomorfismo que $\mbox{Im}\, \iota_{*} \cong \mathbb{Z}(1, 0)$, onde $(1,0) \in K_{1}(J_{0}/ \mathcal{K}_{\mathbb{R}}) = \mathbb{Z}[[U]_{\mathcal{K}_{\mathbb{R}}}]_{1} \oplus
\mathbb{Z}[[V]_{\mathcal{K}_{\mathbb{R}}}]_{1}$.

Pela exatidão da sequência exata ciclíca de seis termos temos que $\psi_{*}: K_{0}(
\mathcal{A}/ \mathcal{K}_{\mathbb{R}}) \longrightarrow \mbox{Im}\, \psi_{*} =
\mathbb{Z}(1,1)$ é um isomorfismo. Logo $K_{0}(\mathcal{A}/
\mathcal{K}_{\mathbb{R}}) = \mathbb{Z}[[Id]_{\mathcal{K}_{\mathbb{R}}}]_{0}$.

Como $\delta_{1} \equiv 0$
temos que $\mbox{Im}\, \psi_{*} \cong \mathbb{Z} \oplus
\mathbb{Z} \cong K_{1}(\mathcal{A}^{\diamond} \oplus
\mathcal{A}^{\diamond})$. Então podemos construir a  seguinte sequência exata curta
\begin{center}
\begin{picture}(350,15)
\put(31,0){0}
\put(44,3){\vector(1,0){25}}
\put(57,0){\makebox[60pt]{$\mathbb{Z}(1,0)$}}
\put(106,3){\vector(1,0){35}}
\put(146,0){\makebox[60pt]{$K_{1}(\mathcal{A}/ \mathcal{K}_{\mathbb{R}})$}}
\put(215,3){\vector(1,0){35}}
\put(250,0){\makebox[60pt]{$ \mathbb{Z} \oplus \mathbb{Z}$}}
\put(307,3){\vector(1,0){25}}
\put(341,0){0}
\put(102,8){\makebox[35pt]{$ i$}}
\put(215,8){\makebox[35pt]{$\psi_{*} $}}
\end{picture}.
\end{center}
Logo $K_{1}(\mathcal{A}/ K_{\mathbb{R}}) = \mathbb{Z}[[A]_{\mathcal{K}_{\mathbb{R}}}]_{1}
\oplus \mathbb{Z}[[B]_{\mathcal{K}_{\mathbb{R}}}]_{1} \oplus \mathbb{Z}[[U]_{\mathcal{K}_{\mathbb{R}}}]_{1}$ onde $A$ e $B$ são
elementos de $\mathcal{A}$ tal que $\psi(A) = (e^{i M}, 1)$ e $\psi(B) = (1,
e^{i M})$. Podemos tomar, por exemplo, $A = L(M)$ e $B = \tilde{L}(M)$, $L$ e
$\tilde{L}$ definidas na página 54.$\cqd$

Agora já é possivel calcular a K-teoria da $C^*$- álgebra $\mathcal{A}$.
\begin{Teo}
Seja $\pi: \mathcal{A} \longrightarrow  \mathcal{A}/
\mathcal{K}_{\mathbb{R}}$ a projeção canônica. Então $\pi_{*}: K_{i}( \mathcal{A})
\longrightarrow \mbox{Im} \, \pi_{*}$ é um isomorfismo.
\end{Teo}
\textbf{Demonstração:}

Consideremos a sequência exata curta
\begin{center}
\begin{picture}(280,15)
\put(15,0){0}
\put(30,3){\vector(1,0){25}}
\put(43,0){\makebox[60pt]{$\mathcal{K}_{\mathbb{R}}$}}
\put(90,3){\vector(1,0){35}}
\put(110,0){\makebox[60pt]{$\mathcal{A}$}}
\put(152,3){\vector(1,0){35}}
\put(181,0){\makebox[60pt]{$\mathcal{A}/ \mathcal{K}_{\mathbb{R}}$}}
\put(230,3){\vector(1,0){25}}
\put(262,0){0}
\put(90,8){\makebox[35pt]{$i $}}
\put(152,8){\makebox[35pt]{$\pi$}}
\end{picture}
\end{center}
e sua sequência exata ciclíca de seis termos para K-teoria

\begin{picture}(350,15)
\put(48,0){\makebox[60pt]{$K_{0}(\mathcal{K}_{\mathbb{R}})$}}
\put(111,3){\vector(1,0){35}}
\put(149,0){\makebox[60pt]{$K_{0}(\mathcal{A})$}}
\put(212,3){\vector(1,0){35}}
\put(250,0){\makebox[60pt]{$K_{0}(\mathcal{A}/ \mathcal{K}_{\mathbb{R}})$}}
\put(111,8){\makebox[35pt]{$ i_{*}$}}
\put(212,8){\makebox[35pt]{$ \pi_{*} $}}
\end{picture}

\begin{picture}(350,45)
\put(78,0){\vector(0,1){35}}
\put(280,35){\vector(0,-1){35}}
\put(81,17){$\delta_{1}  $}
\put(283,17){$ \delta_{0} $}
\end{picture}

\begin{picture}(350,15)
\put(48,0){\makebox[60pt]{$K_{1}(\mathcal{A} / \mathcal{K}_{\mathbb{R}})$}}
\put(146,3){\vector(-1,0){35}}
\put(149,0){\makebox[60pt]{$K_{1}(\mathcal{A})$}}
\put(247,3){\vector(-1,0){35}}
\put(250,0){\makebox[60pt]{$K_{1}(\mathcal{K}_{\mathbb{R}})$}}
\put(111,8){\makebox[35pt]{$ \pi_{*}$}}
\put(212,8){\makebox[35pt]{$ i_{*} $}}
\end{picture}

Substituindo as K-teorias conhecidas temos

\begin{picture}(350,15)
\put(48,0){\makebox[60pt]{$\mathbb{Z}$}}
\put(111,3){\vector(1,0){35}}
\put(149,0){\makebox[60pt]{$K_{0}(\mathcal{A})$}}
\put(212,3){\vector(1,0){35}}
\put(250,0){\makebox[60pt]{$\mathbb{Z}$}}
\put(111,8){\makebox[35pt]{$ i_{*}$}}
\put(212,8){\makebox[35pt]{$ \pi_{*} $}}
\end{picture}

\begin{picture}(350,45)
\put(78,0){\vector(0,1){35}}
\put(280,35){\vector(0,-1){35}}
\put(81,17){$\delta_{1}  $}
\put(283,17){$ \delta_{0} $}
\end{picture}

\begin{picture}(350,15)
\put(48,0){\makebox[60pt]{$\mathbb{Z} \oplus \mathbb{Z} \oplus \mathbb{Z}$}}
\put(146,3){\vector(-1,0){35}}
\put(149,0){\makebox[60pt]{$K_{1}(\mathcal{A})$}}
\put(247,3){\vector(-1,0){35}}
\put(250,0){\makebox[60pt]{$0$}}
\put(111,8){\makebox[35pt]{$ \pi_{*}$}}
\put(212,8){\makebox[35pt]{$ i_{*} $}}
\end{picture}

Sabemos pela Proposição 9.4.4 de \cite{[R]} que o índice para esta sequência é o
índice de Fredholm.

Vamos calcular $\delta_{1}$ nos três geradores de $K_{1}(\mathcal{A} /
\mathcal{K}_{\mathbb{R}})$. Temos que
$\delta_{1}([[L(M)]_{\mathcal{K}_{\mathbb{R}}}]_{1}) =
\delta_{1}([[\tilde{L}(M)]_{\mathcal{K}_{\mathbb{R}}}]_{1}) = 0$ já que eles
são operadores de multiplicação.
 
É fácil ver que $e^{2 \pi i b(M)} b(D) + c(D) \in \mathcal{C}$ logo $[e^{2 \pi
  i b(M)} b(D) + c(D)]_{\mathcal{K}_{\mathbb{R}}} \in \mathcal{C} / \mathcal{K}_{\mathbb{R}}$ e como foi visto no Capítulo
  2, Teorema 2.13, $\delta_{1}([[T]_{\mathcal{K}_{\mathbb{R}}}]_{1}) = ind(T) = w(\sigma_{T})$ nesta álgebra.
Logo $\delta_{1}([[e^{2 \pi
  i b(M)} b(D) + c(D)]_{\mathcal{K}_{\mathbb{R}}}]_{1}) = w(\sigma_{e^{2 \pi
  i b(M)} b(D) + c(D)}) = -1$. Portanto $\delta_{1}$ é sobrejetora e $\delta_{1}(x, y, z) = -z$.

 Usando a sequência exata cíclica de seis termos temos $$K_{0}(\mathcal{A}) =
 \mathbb{Z}[Id]_{0} \hspace{.5cm} \mbox{e}$$ $$K_{1}(\mathcal{A}) =
 \mathbb{Z}[L(M)]_{1} \oplus \mathbb{Z}[\tilde{L}(M)]_{1}.$$ $\cqd$

\end{document}